# Non-dispersive retarded interactions as manifestations of unique symmetry of our space associated with its 3-dimensionality

by

Yakov A. Iosilevskii[1]

## Abstract

It is proved that the inhomogeneous d'Alambertian (wave) differential equation in a Euclidean affine space (EAfS) of any given dimension $n \geq 1$ has a so-called *assembled*, i.e. *non-dispersed*, (*as* contrasted to a so-called *dispersed*), *time-biased* (*retarded* or *advanced*) *solution* if and only if $n=3$ (see Theorem 6.1). To this end, a *special integro-differential formula* (*SIDF*) is deduced for a pair of real-valued *functional forms* (*FF's*) of appropriate kinds that are, for simplicity but without loss of generality, defined and continuous, along with their second-order spatial derivatives, on a closed spherical layer in an adjoined $n$-dimensional ($n$D) Euclidean arithmetical space (EArS) that is restricted by two $(n–1)$D spherical surfaces of radii $r_1$ and $r_2 > r_1$, centered at a given spatial point of observation – a layer that is after all expanded to the whole of the EArS. One of the two FF's is an arbitrary time-biased, and hence *composite-plain* with respect to its spatial variables, real-valued FF that is represented by an appropriate *FF placeholder* (*FFPH*). The latter is parameterized by the double-valued variable '$\lambda$' such that the range of the FFPH comprises retarded FF's at $\lambda=–1$ and advanced FF's at $\lambda=1$. The other FF is a *spherically symmetrical harmonic FF* (*HFF*), which is properly normalized so that if the HFF applies in an $n$D domain that contains its *singular point*, its *associated function* becomes the Green function of an $n$D Poisson (inhomogeneous Laplace) equation. At $n \geq 2$, the SIDF involves two $n$D ($n$-fold) volume Riemann (proper) integrals over the spherical layer and four $(n–1)$D (($n–1$)-fold) surface Riemann integrals over the spherical surfaces. As $r_1 \to 0$ and $r_2 \to \infty$, under certain additional standard assumptions regarding the respective values of a biased FF in the range of the FFPH and of its pertinent derivatives, the SIDF at any given natural number $n \geq 1$ turns into a *boundary-free SIFD* (*BFSIFD*) that expresses the unbiased value of the biased FF at the given temporo-spatial point of observation as the sum of two converging *n*-fold improper integrals over the *n*D EArS (see Theorem 7.1 or equation (1.15) subject to (1.13), (1.14), and (1.16)). One of the two integrals is the negative contraction (convolution) of the HFF and of the biased (retarded or advanced) value of the operand consisting of the d'Alambertian operator followed by an unbiased version of the biased FF. At *n*=3, this integral turns into the known biased (retarded or advanced) solution of the pertinent inhomogeneous d'Alambertian differential equation in the infinite 3D EArS (see Corollary 7.1 or equation (1.19)). The other integral contains the factor $\lambda(3-n)$, while its integrand involves the first-order temporal derivative of the biased FF. This BFSIDF applies also with *n*=1. Thus, the BFSIDF is the forced solution of an inhomogeneous d'Alambertian differential equation in the infinite *n*D EArS if an only if *n*=3. Therefore, the SIDF in general or the BFSIDF in particular can be regarded as a manifestation of *unique symmetry of our 3D Euclidean affine space associated with it 3-dimensionality*, which allows radiation and propagation of *non-dispersive* electromagnetic waves in vacuum and radiation and propagation of *quasi non-dispersive* sound waves in continuous material media – waves that serve as carriers of *detailed*, or *non-dispersive*, *retarded interactions* of material particles. That is to say, the BFSIDF (e.g.) signifies that no *non-dispersed* retarded classical wave interactions and hence no *special relativity* can exist in an *n*D Euclidean real affine space other than a 3D one.

---

[1] Retired from the Israel Ocenographic and Limnological Research Institute, P.O.B. 8030, Haifa 31080, **E-mail:** yi@bezeqint.net. **Phone/Fax**: 972-4-8236071.

## 1. Introduction

For any natural number *n*>3, all pure algebraic properties of an *n-dimensional Euclidean affine space* $\dot{E}_n(R)$ *over the field R of real numbers*, called also an *n*-dimensional Euclidean *real* affine space, i.e. all pertinent axioms and all theorems following from them, are the same as the respective pure algebraic properties of a 3-*dimensional Euclidean real affine space* $\dot{E}_3(R)$. Also, all pure algebraic properties 1- and 2-dimensional Euclidean real affine spaces, $\dot{E}_1(R)$ and $\dot{E}_2(R)$, have analogues among algebraic properties of $\dot{E}_3(R)$ except some minor exceptions. For instance, $\dot{E}_1(R)$ does not have any subspaces of a smaller dimension, whereas $\dot{E}_2(R)$ has only 1-dimensional subspaces. At the same time, a *distinguished hypothetic* $\dot{E}_3(R)$ can be regarded as a *receptacle of nature*, i.e. as a receptacle of matter along with all metamorphoses, which occur to matter in *time* and which are called physical, chemical, biological, etc. processes. *Time* is a hypothetical *non-spatial* one-dimensional continuum that can be regarded as a special version (interpretation) of $\dot{E}_1(R)$, to be denoted by '$\dot{T}(R)$' or briefly by '$\dot{T}$', in which *the above processes* go on in the irreversible direction from past through present to future. It is postulated that, via those processes, $\dot{E}_3(R)$ is united with $\dot{T}(R)$ to form a 4-dimensional pseudo-Euclidean real affine space of index 1 – the *space-time of special theory of relativity*, which is called the *Minkowski space* and which will be denoted by '$\dot{M}_4(R)$'. The presence of gravitating masses in the hypothetic $\dot{E}_3(R)$, along with the inseparable gravitational processes going on in and expanded across $\dot{T}(R)$, changes the known metric properties of $\dot{M}_4(R)$, so that it is replaced by, i.e. as if turns into, the *Riemannian space of general theory of relativity*, $\mathcal{R}$.

Nevertheless, many physical processes that occurs in a relatively small regions of $\dot{E}_3(R)$ during relatively small spans of $\dot{T}(R)$ are satisfactorily described as if they occurred in the manifold being the direct product $\dot{T}(R) \times \dot{E}_3(R)$, whereas many others are satisfactorily described in $\dot{M}_4(R)$. For instance, radiation of electromagnetic waves and their propagation in vacuum and radiation of sound waves and their propagation in material medium are most often regarded as ones that occur in $\dot{T}(R) \times \dot{E}_3(R)$. In this case, the *non-dispersive retarded solution* of an *inhomogeneous wave (d'Alambertian) second-order differential equation* is naturally interpreted as the pertinent wave signal that either accomplishes *non-dispersive* retarded non-local interaction between the source of the wave signal and its responder (receiver) or that carries *non-dispersive* retarded information of the source to an observer through the vacuous or material medium.

There is a great many of definitions of either of the synonymous word-combinations (phrases) "*classical causality principle*" and "*Laplace causality principle*" to be abbreviated as "*causality principle*" unless stated otherwise (see, e.g., Bohm [1957]). Most generally, i.e. with minimum of assumptions, the causality principle can be defined as the *one-directionality of time* (*orthochrony*). *Necessary interrelations* of things, which are hypostases (intrinsic aspects, ways of existence) of the things, are some other causality principles, which are commonly called "*causality laws*". The Green function of an equation of diffusion or of propagation of heat in an affine 3D Euclidean space can be called the causality principle of the respective phenomenon, whereas the retarded Green function or retarded solution of a wave (d'Alambertian) equation in such a space can be called the causality principle of the wave phenomenon or, less explicitly (more generally), a *non-dispersive causality principle* and also, perhaps, a *detailed causality principle*, in analogy with "*principle of detailed* *equilibrium*". In this case, the following general question of philosophical (epistemological) importance may be raised:

«Can the direct product $\dot{T}(R) \times \dot{E}_n(R)$ with $n \neq 3$ serve as a receptacle of some nature that satisfies the classical *non-dispersive* (*detailed*) *causality principle*?»

This general question reduces to the following concrete one:

«Does an *inhomogeneous wave* (*d'Alambertian*) *second-order partial differential equation* in the infinite direct product $\dot{T}(R) \times \dot{E}_n(R)$ with $n \neq 3$ have a *non-dispersive* (*detailed*) *retarded solution*?»

The *categorical* (*unconditional*) *negative answer* to the latter question is given in this paper.

The antonymous qualifiers "*non-dispersive*" ("*detailed*") and "*dispersive*" to "solution of a wave equation" (particularly to "Green function of a wave equation") or to "causality principle" are not, likely, the best ones, but I could not find any better English words to do their duties for expressing the pertinent properties, which will be made explicit below in this section.

To each $\dot{E}_n(R)$ there is an adjoined *n*-dimensional Euclidean abstract linear (vector) space $\hat{E}_n(R)$ over the field $R$ and a surjection from $\dot{E}_n \times \dot{E}_n$ onto $\hat{E}_n$; $\hat{E}_n$ is the underlying set of points of $\dot{E}_n(R)$ and $\hat{E}_n$ is the underlying set of abstract vectors of $\hat{E}_n(R)$. Relative to its any orthonormal basis, $\hat{E}_n(R)$ is isomorphic to the respective adjoined *n*-dimensional Euclidean arithmetical linear (vector) space $\overline{E}_n(R)$ over $R$, whose underlying set $\overline{E}_n$ consists of ordered *n*-tuples of real numbers. The above said applies also in the case of *n*=1, provided that the ordered 1-tuple of a real number is by definition the singleton of that number. Consequently, the above said applies particularly with '$\dot{T}(R)$', '$\hat{T}(R)$', '$\overline{T}(R)$', '$\dot{T}$', '$\hat{T}$', '$\overline{T}$' in place of '$\dot{E}_n(R)$', '$\hat{E}_n(R)$', '$\overline{E}_n(R)$', '$\dot{E}_n$', '$\hat{E}_n$', '$\overline{E}_n$' respectively.

Any classical physical field in $\dot{T}(R) \times \dot{E}_3(R)$ is in fact described by a certain real-valued function, which is defined on $R \times \overline{X}_3$, the understanding being that $\overline{X}_3$ is a certain subset of $\overline{E}_3$, while successive real numbers of $R$ are instants of time indicated by a certain clock. Therefore, in order to answer the second one of the two questions raised above, we have deduced, in analogy with the ordinary second Green's formula for two ordinary (plain, not composite) real-valued functional forms defined on a closed domain of $\overline{E}_3$ (see, e.g., Tikhonov and Samarskii [1990. chapter IV, §2, n°1, pp. 314–320]), an *integro-differential* *formula* (*IDF*) for two [*quasi-*]*time-biased* (*retarded* or *advanced*), and hence *composite-plain with respect to its spatial variables,* real-valued functional forms $\phi^{\langle 1,n \rangle}(x_{*0} + \lambda x'_{[1,n]}, \bar{x}_{[1,n]})$ and $\psi^{\langle 1,n \rangle}(x_{*0} + \lambda x'_{[1,n]}, \bar{x}_{[1,n]})$, which are defined and continuous along with their all second-order spatial derivatives on a given closed subset $\bar{X}_n^c$ of $\bar{E}_n$ ($\bar{X}_n^c \subset \bar{E}_n$), subject to the following notations:

$$x_{*0} + \lambda x'_{[1,n]} \in R, \lambda \in \{-1, 1\}, \tag{1.1}$$

$$\bar{x}_{[1,1]} \stackrel{\equiv}{=} \langle x_1 \rangle = \{x_1\} \in \bar{E}_1, \tag{1.2}$$

$$\bar{x}_{[1,2]} \stackrel{\equiv}{=} \langle x_1, x_2 \rangle \stackrel{\equiv}{=} \{\{x_1\}, \{x_1, x_2\}\} \in \bar{E}_2, \tag{1.3}$$

$$\bar{x}_{[1,n]} \stackrel{\equiv}{=} \langle x_i \rangle_{i \in \omega_{1,n}} \stackrel{\equiv}{=} \langle x_1, x_2, \ldots, x_{n-1}, x_n \rangle \stackrel{\equiv}{=} \langle \bar{x}_{[1,n-1]}, x_n \rangle$$
$$= \underbrace{\langle \langle \ldots \langle x_1, x_2 \rangle, x_3 \rangle, \ldots, x_{n-1} \rangle, x_n \rangle}_{n-1} \in \bar{E}_n, \tag{1.4}$$

$$\bar{x}'_{[1,n]} \stackrel{\equiv}{=} \bar{x}_{[1,n]} - \bar{x}_{*[1,n]} = \langle x_1 - x_{*1}, \ldots, x_n - x_{*n} \rangle \stackrel{\equiv}{=} \langle x'_1, \ldots, x'_n \rangle \in \bar{E}_n, \tag{1.5}$$

$$x'_{[1,n]} \stackrel{\equiv}{=} |\bar{x}'_{[1,n]}| = \sqrt{\sum_{i=1}^{n} x'^2_i} = \sqrt{\sum_{i=1}^{n} (x_i - x_{*i})^2} \geq 0, \tag{1.6}$$

$$x'_{[1,1]} = |x'_1| = |x_1 - x_{*1}|. \tag{1.7}$$

In this case, '$\lambda$' is a parameter that may take on two values: –1 and 1, and $\langle x_{*0}, \bar{x}_{*[1,n]} \rangle \in R \times \bar{E}_n$ is a given (fixed) temporo-spatial point of observation (measuring) of the hypothetical physical fields $\phi^{\langle 1,n \rangle}$ and $\psi^{\langle 1,n \rangle}$, the understanding being that $\bar{x}_{*[1,n]} \notin \bar{X}_n^c$. Therefore, if the functional form $\psi^{\langle 1,n \rangle}(x_0, \bar{x}_{[1,n]})$ (e.g.) is assumed to be defined and continuous along with its *m*th-order derivative (*m*≥0) with respect to the [*quasi-*]time variable '$x_0$' then the properties of differentiability of a time-biased functional form $\psi^{\langle 1,n \rangle}(x_{*0} + \lambda x'_{[1,n]}, \bar{x}_{[1,n]})$ with respect to spatial variables '$x_1$', '$x_2$', …, '$x_n$' up to the order *m* are the same as those of the functional form $\psi^{\langle 1,n \rangle}(x_0, \bar{x}_{[1,n]})$. For the sake of being specific and for convenience, but without loss of generality, if *n*≥2 then the region of integration $\bar{X}_n^c$ is supposed to be a closed (closed-closed) *n*-dimensional spherical layer $\bar{C}_n^{cc}(r_1, r_2, \bar{x}_{*[1,n]})$ that is restricted by two (*n*–1)–dimensional spherical surfaces $\bar{B}_n^b(r_1, \bar{x}_{*[1,n]})$ and $\bar{B}_n^b(r_2, \bar{x}_{*[1,n]})$ of radii $r_1$ and $r_2 > r_1$, which are centered at the point of observation $\bar{x}_{*[1,n]} \in \bar{E}_n$; the superscripts "c", "o", and "b" stand for "closed", "open", and "boundary" respectively. Accordingly,

$$\bar{C}_1^{cc}(r_1, r_2, \bar{x}_{*[1,1]}) \equiv [x_{*1} - r_2, x_{*1} - r_1] \cup [x_{*1} + r_1, x_{*1} + r_2]. \tag{1.9}$$

Instead of $\phi^{\langle 1,n \rangle}(x_{*0} + \lambda x'_{[1,n]}, \bar{x}_{[1,n]})$, we then utilize the *normalized spherically symmetric harmonic functional form* $\eta^{\langle n \rangle}(x'_{[1,n]})$, which is defined as:

$$\eta^{\langle n \rangle}(x'_{[1,n]}) \equiv \begin{cases} a_n x'^{2-n}_{[1,n]} & \text{if } n \neq 2 \text{ (a)} \\ -a_2 \ln x'_{[1,2]} & \text{if } n = 2 \text{ (b)}, \\ a_1 |x'_1| & \text{if } n = 1 \text{ (c)} \end{cases} \tag{1.10}$$

subject to (1.6) and also subject to

$$a_n \equiv \begin{cases} \dfrac{1}{(n-2)S_n(1)} = \dfrac{\Gamma(n/2)}{2(n-2)\pi^{n/2}} & \text{if } n \geq 3 \text{ (a)} \\ \dfrac{1}{S_2(1)} = \dfrac{1}{2\pi} & \text{if } n = 2 \text{ (b)}. \\ -\dfrac{1}{S_1(1)} = -\dfrac{1}{2} & \text{if } n = 1 \text{ (c)} \end{cases} \tag{1.11}$$

$S_n(1)$, defined as

$$S_n(1) = \frac{2\pi^{n/2}}{\Gamma(n/2)} \text{ for each } n \geq 2, \tag{1.12}$$

is the area of an ($n$–1)-dimensional spherical surface in $\bar{E}_n$ of radius 1 (see Appendix A and cf. Sommerville [1958, pp. 135-137]), i.e. the area of $\bar{B}_n^b(1, x_{*[1,n]})$, while $\Gamma$ is the gamma-function. By (1.7), the expression (1.10c) is the instance of (1.10a). It is understood that $\bar{B}_1^b(1, x_{*[1,n]}) = \{-1, 1\}$, so that '$S_1(1)$' is not defined. It is, however, convenient to utilize the formula (1.11a) at $n=1$ with $a_1 \equiv -1/2$ in order to formally set $S_1(1) \equiv 2$, which is in agreement with (1.12) at $n=1$ for $\Gamma(1/2) = \pi^{1/2}$. This definition can be understood as expressing the fact that the boundary $\bar{B}_1^b(1, x_{*[1,1]})$ consists of *two points* –1 and 1.

In no connection with any values of '$a_n$', for each $n \geq 1$:

$$\Delta_{[1,n]} \eta^{\langle n \rangle}(x'_{[1,n]}) \equiv \sum_{i=1}^{n} \frac{\partial^2}{\partial x_i^2} \eta^{\langle n \rangle}(x'_{[1,n]}) = 0 \text{ for each } \bar{x}_{[1,n]} \in \bar{E}_n - \{\bar{x}_{*[1,n]}\}, \tag{1.13}$$

whereas owing to (1.11)

$$\Delta_{[1,n]} \eta^{\langle n \rangle}(x'_{[1,n]}) = -\delta^{\langle n \rangle}(\bar{x}'_{[1,n]}) \text{ for each } \bar{x}_{[1,n]} \in \bar{E}_n, \tag{1.14}$$

where $\delta^{\langle n \rangle}$ is the *n*-ary Dirac $\delta$-function (see Appendix B). Hence, the function $\eta^{\langle n \rangle}$ is the Green function of an *n*-dimensional Poisson (inhomogeneous Laplace) equation in $\overline{E}_{\{n\}}$. The version of the IDF with '$\eta^{\langle n \rangle}(x'_{[1,n]})$' in place of '$\phi^{\langle 1,n \rangle}(x_{*0} + \lambda x'_{[1,n]}, \overline{x}_{[1,n]})$' is called the *special IDF* (*SIDF*).

For each $n \geq 2$, the SIDF involves two *n*-fold volume integrals over $\overline{C}_n^{cc}(r_1, r_2, \overline{x}_{*[1,n]})$, two (*n*–1)-fold surface integrals over $\overline{B}_n^b(r_1, \overline{x}_{*[1,n]})$, and two (*n*–1)-fold surface integrals over $\overline{B}_n^b(r_2, \overline{x}_{*[1,n]})$. Under the general assumption that the functional form $\psi^{\langle 1,n \rangle}(x_0, \overline{x}_{[1,n]})$, and hence $\psi^{\langle 1,n \rangle}(x_{*0} + \lambda x'_{[1,n]}, \overline{x}_{[1,n]})$, is defined and continuous along with their all second-order derivatives on $R \times \overline{C}_n^{cc}(r_1, r_2, \overline{x}_{*[1,n]})$, it is proved the following.

i) As $r_1 \to 0$, one of the two surface integrals over $\overline{B}_{\{n\}}^b(r_1, \overline{x}_{*[1,n]})$ reduces to $\psi^{\langle 1,n \rangle}(x_{*0}, \overline{x}_{*[1,n]})$, while the other one vanishes.

ii) As $r_2 \to \infty$, both surface integrals over $\overline{B}_n^b(r_2, \overline{x}_{*[1,n]})$ vanish.

iii) As $r_1 \to 0$ and $r_2 \to \infty$ both proper (Riemann) volume integrals turn into converging improper volume integrals.

In accordance with the above properties i–iii, for each $n \geq 3$, as $r_1 \to 0$ and $r_2 \to \infty$, the SIDF becomes:

$$\psi^{\langle 1,n \rangle}(x_{*0}, \overline{x}_{*[1,n]})$$
$$= -\int_{\overline{E}_n} \eta^{\langle n \rangle}(x'_{[1,n]}) \left[ \left( \sum_{i=1}^n \frac{\partial^2}{\partial x_i^2} - \frac{\partial^2}{\partial x_0^2} \right) \psi^{\langle 1,n \rangle}(x_0, \overline{x}_{[1,n]}) \right]_{x_0 = x_{*0} + \lambda x'_{[1,n]}} dv_{\{n\}} \quad (1.15)$$
$$- \lambda(n-3) \int_{\overline{E}_n} \varsigma^{\langle n \rangle}(x'_{[1,n]}) \left[ \frac{\partial \psi^{\langle 1,n \rangle}(x_0, \overline{x}_{[1,n]})}{\partial x_0} \right]_{x_0 = x_{*0} + \lambda x'_{[1,n]}} dv_{\{n\}}$$

subject to (1.10) and also subject to

$$\varsigma^{\langle n \rangle}(x'_{[1,n]}) \equiv \begin{cases} a_n x'^{1-n}_{[1,p]} & \text{if } n \neq 2 \text{ (a)} \\ -\dfrac{a_2 (\ln x'_{[1,2]} + 2)}{x'_{[1,2]}} & \text{if } n = 2 \text{ (b)} \\ a_1 & \text{if } n = 1 \text{ (c)} \end{cases} \quad (1.16)$$

Particularly, by (1.10)–(1.12) and (1.16), equation (1.15) reduces to

$$\psi^{\langle 1,3\rangle}(x_{*0},\overline{x}_{*[1,3]})$$
$$=\frac{1}{4\pi}\int_{\overline{E}_3}\frac{1}{x'_{[1,3]}}\left[\left(\sum_{i=1}^{3}\frac{\partial^2}{\partial x_i^2}-\frac{\partial^2}{\partial x_0^2}\right)\psi^{\langle 1,3\rangle}(x_0,\overline{x}_{[1,3]})\right]_{x_0=x_{*0}+\lambda x'_{[1,3]}}dv_{\{3\}} \quad \text{at } n=3, \tag{1.15a}$$

$$\psi^{\langle 1,2\rangle}(x_{*0},\overline{x}_{*[1,2]})=\frac{1}{2\pi}\int_{\overline{E}_2}\ln x'_{[1,2]}\left[\left(\sum_{i=1}^{2}\frac{\partial^2}{\partial x_i^2}-\frac{\partial^2}{\partial x_0^2}\right)\psi^{\langle 1,2\rangle}(x_0,\overline{x}_{[1,2]})\right]_{x_0=x_{*0}+\lambda x'_{[1,2]}}dv_{\{2\}}$$
$$-\frac{\lambda}{2\pi}\int_{\overline{E}_2}\frac{(\ln x'_{[1,2]}+2)}{x'_{[1,2]}}\left[\frac{\partial\psi^{\langle 1,2\rangle}(x_0,\overline{x}_{[1,2]})}{\partial x_0}\right]_{x_0=x_{*0}+\lambda x'_{[1,2]}}dv_{\{2\}} \quad \text{at } n=2, \tag{1.15b}$$

$$\psi^{\langle 1,1\rangle}(x_{*0},\{x_{*1}\})=\frac{1}{2}\int_{-\infty}^{\infty}|x_1|\left[\left(\frac{\partial^2}{\partial x_1^2}-\frac{\partial^2}{\partial x_0^2}\right)\psi^{\langle 1,1\rangle}(x_0,\{x_1\})\right]_{x_0=x_{*0}+\lambda|x_1|}dx_1$$
$$-\lambda\int_{-\infty}^{\infty}\left[\frac{\partial\psi^{\langle 1,1\rangle}(x_{*0},\{x_{*1}\})}{\partial x_0}\right]_{x_0=x_{*0}+\lambda|x_1|}dx_1 \quad \text{at } n=1. \tag{1.15c}$$

Let by definition

$$f^{\langle 1,n\rangle}(x_0,\overline{x}_{[1,n]})\stackrel{\doteq}{=}-\left(\sum_{i=1}^{n}\frac{\partial^2}{\partial x_i^2}-\frac{\partial^2}{\partial x_0^2}\right)\psi^{\langle 1,n\rangle}(x_0,\overline{x}_{[1,n]}). \tag{1.17}$$

Alternatively, one may assume (postulate) that the functional form $f^{\langle 1,n\rangle}(x_0,\overline{x}_{[1,n]})$ denotes a density of sources of the physical field denoted by the functional form $\psi^{\langle 1,n\rangle}(x_0,\overline{x}_{[1,n]})$, so that the latter is a solution of the inhomogeneous wave (d'Alambertian) differential equation:

$$\left(\sum_{i=1}^{n}\frac{\partial^2}{\partial x_i^2}-\frac{\partial^2}{\partial x_0^2}\right)\psi^{\langle 1,n\rangle}(x_0,\overline{x}_{[1,n]})\stackrel{\doteq}{=}-f^{\langle 1,n\rangle}(x_0,\overline{x}_{[1,n]}) \tag{1.18)a}$$

subject to certain additional satisfiable conditions. In either case, the SIDF (1.15) at *n*=3, i.e. at *(!.15a)*, becomes:

$$\psi^{\langle 1,3\rangle}(x_{*0},\overline{x}_{*[1,3]})=\frac{1}{4\pi}\int_{\overline{E}_3}\frac{f^{\langle 1,n\rangle}(x_{*0}+\lambda x'_{[1,3]},\overline{x}_{[1,3]})}{x'_{[1,n]}}dv_{\{3\}}, \tag{1.19}$$

which presents both known forced solutions of the wave equation (1.18) at *n*=3 in $\overline{E}_3$, namely the retarded one at $\lambda=-1$ and the advanced one at $\lambda=1$.

Given a natural number $n\geq 1$, let $G^{\langle 1,n\rangle}$ be the Green function of the equation (1.18) with '=' in place of '≐', which satisfies the equation

$$\left(\sum_{i=1}^{n}\frac{\partial^2}{\partial x_i^2}-\frac{\partial^2}{\partial x_0^2}\right)G^{\langle 1,n\rangle}(x'_0,\overline{x}'_{[1,n]})=-\delta(x'_0)\delta^{\langle n\rangle}(\overline{x}'_{[1,n]}) \tag{1.20}$$

for each $x'_0\stackrel{\doteq}{=}x_0-x_{*0}\in\overline{T}$ and each $\overline{x}'_{[1,n]}\stackrel{\doteq}{=}\overline{x}_{[1,n]}-\overline{x}_{*[1,n]}\in\overline{E}_n$.

Hence,

$$\psi^{\langle 1,n \rangle}(x_{*0}, \bar{x}_{*[1,n]}) = \int\limits_{\bar{E}_n} \int\limits_{-\infty}^{\infty} G^{\langle 1,n \rangle}(x_{*0} - x_0, \bar{x}_{*[1,n]} - \bar{x}_{[1,n]}) f^{\langle 1,n \rangle}(x_0, \bar{x}_{[1,n]}) dx_0 dv_{\{n\}}. \qquad (1.21)$$

Comparison of (1.21) at $n=3$ and (1.19) shows that equation (1.20) at $n=3$ has two solutions:

$$G^{\langle 1,3 \rangle}(x_{*0} - x_0, \bar{x}_{*[1,3]} - \bar{x}_{[1,3]}) = G^{\lambda\langle 1,3 \rangle}(x_{*0} - x_0, \bar{x}_{*[1,3]} - \bar{x}_{[1,3]})$$
$$= \frac{\delta(x_{*0} - x_0 + \lambda x'_{[1,3]})}{4\pi x'_{[1,3]}} = \eta^{(3)}(x'_{[1,3]}) \delta(x_{*0} - x_0 + \lambda x'_{[1,3]}) \text{ with } \lambda \in \{-1, 1\}, \qquad (1.22)$$

so that $G^{-1\langle 1,3 \rangle}$ is the *retarded Green function* and $G^{1\langle 1,3 \rangle}$ is the *advanced Green function* for the wave problem (cf. Morse and Feshbach [1953, section 7.3]).

**Comment 1.1.** Indeed, the pertinent straightforward calculations yield:

$$\sum_{i=1}^{3} \frac{\partial^2}{\partial x_{*i}^2} \frac{\delta(x_{*0} - x_0 + \lambda x'_{[1,3]})}{x'_{[1,3]}} = \left( \sum_{i=1}^{3} \frac{\partial^2 x'^{-1}_{[1,3]}}{\partial x_{*i}^2} \right) \delta(x_{*0} - x_0 + \lambda x'_{[1,3]})$$
$$+ 2\sum_{i=1}^{3} \frac{\partial x'^{-1}_{[1,3]}}{\partial x_{*i}} \frac{\partial \delta(x_{*0} - x_0 + \lambda x'_{[1,3]})}{\partial x_{*i}} + \frac{1}{x'_{[1,3]}} \sum_{i=1}^{3} \frac{\partial^2 \delta(x_{*0} - x_0 + \lambda x'_{[1,3]})}{\partial x_{*i}^2}. \qquad (1.22_1)$$

In this case, by (1.5), (1.6), and (1.10)–(1.14) at $n=3$, it follows that

$$\left( \sum_{i=1}^{3} \frac{\partial^2 x'^{-1}_{[1,3]}}{\partial x_{*i}^2} \right) \delta(x_{*0} - x_0 + \lambda x'_{[1,3]}) = -4\pi \delta^{(3)}(\bar{x}'_{[1,n]}) \delta(x_{*0} - x_0 + \lambda x'_{[1,3]})$$
$$= -4\pi \delta(x_{*0} - x_0) \delta^{(3)}(\bar{x}'_{[1,n]}), \qquad (1.22_2)$$

$$2\sum_{i=1}^{3} \frac{\partial x'^{-1}_{[1,3]}}{\partial x_{*i}} \frac{\partial \delta(x_{*0} - x_0 + \lambda x'_{[1,3]})}{\partial x_{*i}}$$
$$= 2\left( \sum_{i=1}^{3} \frac{\partial x'_{[1,3]}}{\partial x_{*i}} \frac{\partial x'_{[1,3]}}{\partial x_{*i}} \right) \frac{\partial x'^{-1}_{[1,3]}}{\partial x'_{[1,3]}} \frac{\partial \delta(x_{*0} - x_0 + \lambda x'_{[1,3]})}{\partial x'_{[1,3]}} = -\frac{2\lambda}{x'^{2}_{[1,3]}} \frac{\partial \delta(x_{*0} - x_0 + \lambda x'_{[1,3]})}{\partial x_{*0}}, \qquad (1.22_3)$$

$$\sum_{i=1}^{3} \frac{\partial^2 \delta(x_{*0} - x_0 + \lambda x'_{[1,3]})}{\partial x_{*i}^2} = \sum_{i=1}^{3} \frac{\partial}{\partial x_{*i}} \left[ \frac{\partial x'_{[1,3]}}{\partial x_{*i}} \frac{\partial \delta(x_{*0} - x_0 + \lambda x'_{[1,3]})}{\partial x'_{[1,3]}} \right]$$
$$= \lambda \sum_{i=1}^{3} \frac{\partial}{\partial x_{*i}} \left[ \frac{\partial x'_{[1,3]}}{\partial x_{*i}} \frac{\partial \delta(x_{*0} - x_0 + \lambda x'_{[1,3]})}{\partial x_{*0}} \right] = \lambda \left( \sum_{i=1}^{3} \frac{\partial^2 x'_{[1,3]}}{\partial x_{*i}^2} \right) \frac{\partial \delta(x_{*0} - x_0 + \lambda x'_{[1,3]})}{\partial x_{*0}}$$
$$+ \lambda \sum_{i=1}^{3} \frac{\partial x'_{[1,3]}}{\partial x_{*i}} \frac{\partial^2 \delta(x_{*0} - x_0 + \lambda x'_{[1,3]})}{\partial x_{*i} \partial x_{*0}} = \frac{2\lambda}{x'_{[1,3]}} \frac{\partial \delta(x_{*0} - x_0 + \lambda x'_{[1,3]})}{\partial x_{*0}} \qquad (1.22_4)$$
$$+ \lambda \left( \sum_{i=1}^{3} \frac{\partial x'_{[1,3]}}{\partial x_{*i}} \frac{\partial x'_{[1,3]}}{\partial x_{*i}} \right) \frac{\partial^2 \delta(x_{*0} - x_0 + \lambda x'_{[1,3]})}{\partial x_{*0} \partial x'_{[1,3]}}$$
$$= \frac{2\lambda}{x'_{[1,3]}} \frac{\partial \delta(x_{*0} - x_0 + \lambda x'_{[1,3]})}{\partial x_{*0}} + \lambda^2 \frac{\partial^2 \delta(x_{*0} - x_0 + \lambda x'_{[1,3]})}{\partial x_{*0}^2}.$$

By (1.22₁)–(1.22₄), subject to $\lambda^2=1$, it follows from (1.22) that

$$\left(\sum_{i=1}^{3}\frac{\partial^2}{\partial x_{*i}^2} - \frac{\partial^2}{\partial x_{*0}^2}\right)G^{\langle 1,3\rangle}\left(x_{*0} - x_0, \bar{x}_{*[1,n]} - \bar{x}_{[1,n]}\right) = -\delta(x_{*0} - x_0)\delta^{\langle 3\rangle}\left(\bar{x}_{*[1,n]} - \bar{x}_{[1,n]}\right) \quad (1.20a)$$

for each $x_{*0} - x_0 \in \overline{T}$ and each $x_{*[1,n]} - \bar{x}_{[1,n]} \in \overline{E}_n$.

which is the instance of (1.20) at $n=3$ with $\langle x_0, \bar{x}_{[1,3]}\rangle$ and $\langle x_{*0}, \bar{x}_{*[1,3]}\rangle$ exchanged.•

If the source density $f^{\langle 1,3\rangle}(x_0, \bar{x}_{[1,3]})$ is independent of $x_3$ then the *forced solution* $\psi^{\langle 1,3\rangle}(x_0, \bar{x}_{[1,3]})$ of equation (1.18) at $n=3$ is also independent of $x_3$. A «point» source for such a *two-dimensional wave problem* is the *homogeneous linear source*, spreading from $x_3 = -\infty$ to $x_3 = \infty$ along the straight line, parallel to the $x_3$-axis and passing through the point $\bar{x}_{*[1,2]}$ in the $\langle x_1, x_2\rangle$-plane. Therefore, given $\lambda \in \{-1,1\}$, the two-dimensional Green functional form $G^{\lambda\langle 1,2\rangle}(x_{*0} - x_0, \bar{x}_{*[1,2]} - \bar{x}_{[1,2]})$ can be calculated from (1.22) as:

$$G^{\lambda\langle 1,2\rangle}\left(x_{*0} - x_0, \bar{x}_{*[1,2]} - \bar{x}_{[1,2]}\right) = \frac{1}{4\pi}\int_{-\infty}^{\infty}\frac{\delta(x_{*0} - x_0 + \lambda x'_{[1,3]})}{x'_{[1,3]}}dx_3. \quad (1.23)$$

From (1.5) and (1.6) at $n=3$, it follows that

$$x'_{[1,3]} = \sqrt{x'_{[1,2]}{}^2 + (x'_3)^2} \geq 0, \quad (1.24)$$

whence

$$x'_3 = \pm\sqrt{x'_{[1,3]}{}^2 - x'_{[1,2]}{}^2}, \quad (1.25)$$

$$dx'_{[1,3]}/dx'_3 = x'_3/x'_{[1,3]}, \text{ so that } dx'_3/dx'_{[1,3]} = x'_{[1,3]}/x'_3. \quad (1.26)$$

Therefore, (1.23) for each $\lambda \in \{-1,1\}$ can be developed thus:

$$G^{\lambda\langle 1,2\rangle}\left(x_{*0} - x_0, \bar{x}_{*[1,2]} - \bar{x}_{[1,2]}\right) = \frac{1}{4\pi}\int_{-\infty}^{\infty}\frac{\delta(x_{*0} - x_0 + \lambda x'_{[1,3]})}{x'_{[1,3]}}dx_3$$

$$= \frac{1}{2\pi}\int_{0}^{\infty}\frac{\delta(x_{*0} - x_0 + \lambda x'_{[1,3]})}{x'_{[1,3]}}dx'_3 = \frac{1}{2\pi}\int_{x'_{[1,2]}}^{\infty}\frac{\delta(x_{*0} - x_0 + \lambda x'_{[1,3]})}{x'_{[1,3]}}\frac{dx'_3}{dx'_{[1,3]}}dx'_{[1,3]} \quad (1.23_1)$$

$$= \frac{1}{2\pi}\int_{x'_{[1,2]}}^{\infty}\frac{\delta(x_{*0} - x_0 + \lambda x'_{[1,3]})}{x'_3}dx'_{[1,3]} = \frac{1}{2\pi}\int_{x'_{[1,2]}}^{\infty}\frac{\delta(x_{*0} - x_0 + \lambda x'_{[1,3]})}{\sqrt{x'_{[1,3]}{}^2 - x'_{[1,2]}{}^2}}dx'_{[1,3]}..$$

Given $x'_{[1,2]}{}^2 > 0$, the following two statements are true owing to (1.24):

a) There is exactly one $x'_{[1,3]} = x_{*0} - x_0 > 0$ if and only if $x_{*0} - x_0 - x'_{[1,2]} > 0$.

b) There is exactly one $x'_{[1,3]} = x_0 - x_{*0} > 0$ if and only if $x_{*0} - x_0 + x'_{[1,2]} < 0$.

Hence, (1.23₁) yields:

$$G^{-1\langle 1,2\rangle}(x_{*0}-x_0, \overline{x}_{*[1,2]}-\overline{x}_{[1,2]}) = \frac{1}{2\pi}\int_{x'_{[1,2]}}^{\infty}\frac{\delta(x_{*0}-x_0-x'_{[1,3]})}{\sqrt{x'^{2}_{[1,3]}-x'^{2}_{[1,2]}}}dx'_{[1,3]}$$

$$= \frac{1}{2\pi\sqrt{(x_{*0}-x_0)^2-x'^{2}_{[1,2]}}} \cdot \begin{cases} 1 \text{ if } x_{*0}-x_0-x'_{[1,2]}>0 \\ 0 \text{ if } x_{*0}-x_0-x'_{[1,2]}<0 \end{cases} = \frac{H(x_{*0}-x_0-x'_{[1,2]})}{2\pi\sqrt{(x_{*0}-x_0)^2-x'^{2}_{[1,2]}}},$$

(1.23')

$$G^{1\langle 1,2\rangle}(x_{*0}-x_0, \overline{x}_{*[1,2]}-\overline{x}_{[1,2]}) = \frac{1}{2\pi}\int_{x'_{[1,2]}}^{\infty}\frac{\delta(x_{*0}-x_0+x'_{[1,3]})}{\sqrt{x'^{2}_{[1,3]}-x'^{2}_{[1,2]}}}dx'_{[1,3]}$$

$$= \frac{1}{2\pi\sqrt{(x_{*0}-x_0)^2-x'^{2}_{[1,2]}}} \cdot \begin{cases} 1 \text{ if } x_{*0}-x_0+x'_{[1,2]}<0 \\ 0 \text{ if } x_{*0}-x_0+x'_{[1,2]}>0 \end{cases} = \frac{H(x_0-x_{*0}-x'_{[1,2]})}{2\pi\sqrt{(x_{*0}-x_0)^2-x'^{2}_{[1,2]}}},$$

(1.23'')

i.e.

$$G^{\lambda\langle 1,2\rangle}(x_{*0}-x_0, \overline{x}_{*[1,2]}-\overline{x}_{[1,2]}) = \frac{H(\lambda(x_0-x_{*0})-x'_{[1,2]})}{2\pi\sqrt{(x_{*0}-x_0)^2-x'^{2}_{[1,2]}}} \text{ for each } \lambda\in\{-1,1\}, \quad (1.27)$$

where $H$ is the Heaviside function defined as:

$$H(x) \equiv \begin{cases} 1 \text{ if } x>0 \\ 0 \text{ if } x<0 \end{cases}. \quad (1.28)$$

If the source density $f^{\langle 1,3\rangle}(x_0, \overline{x}_{[1,3]})$ is independent of both $x_2$ and $x_3$ then the forced solution $\psi^{\langle 1,3\rangle}(x_0, \overline{x}_{[1,3]})$ of equation (1.18) at $n=3$ is also independent of both $x_2$ and $x_3$. A «point» source for such a *one-dimensional wave problem* is the *homogeneous source*, *spreading over the infinite plane*, parallel to the $\langle x_2, x_3\rangle$-plane and passing through the point $x_{*1}$ of the $x_1$-axis. Therefore, the one-dimensional Green functional forms $G^{-1\langle 1,1\rangle}(x_{*0}-x_0,\{x_{*1}-x_1\})$ and $G^{-1\langle 1,1\rangle}(x_{*0}-x_0,\{x_{*1}-x_1\})$ can be calculated by integrating the two-dimensional ones given by (1.23') and (1.23''), or (1.27), with respect to $x_2$ from $-\infty$ to $\infty$. Hence, given $x'_1\in R$, it follows from (1.23') and (1.23'') that for each $\lambda\in\{-1,1\}$:

$$G^{\lambda\langle 1,1\rangle}(x_{*0}-x_0,\{x_{*1}-x_1\}) = \int_{-\infty}^{\infty}G^{\lambda\langle 1,2\rangle}(x_{*0}-x_0, \overline{x}_{*[1,2]}-\overline{x}_{[1,2]})dx_2$$

$$= \frac{1}{2\pi}\int_{-s}^{s}\frac{dx'_2}{\sqrt{s^2-x'^{2}_2}} \cdot \begin{cases} H(x_{*0}-x_0-x'_{[1,2]}) \text{ if } \lambda=-1 \\ H(x_0-x_{*0}-x'_{[1,2]}) \text{ if } \lambda=1 \end{cases},$$

(1.29)

where

$$s \equiv \sqrt{(x_{*0}-x_0)^2-x'^{2}_1} > 0. \quad (1.30)$$

Since $x'_{[1,2]} = \sqrt{x'^2_1 + x'^2_2} \geq 0$ therefore, given $|x'_1| > 0$, if $x_{*0} - x_0 - x'_{[1,2]} > 0$ then $x_{*0} - x_0 - |x'_{1,}| > 0$, whereas if $x_{*0} - x_0 - |x'_{1,}| < 0$ then $x_{*0} - x_0 - x'_{[1,2]} < 0$, and similarly with '$x_0 - x_{*0}$' in place of '$x_{*0} - x_0$'. Also,

$$\int_{-s}^{s} \frac{dx'_2}{\sqrt{s^2 - x'^2_{,2}}} = 2\int_0^s \frac{dx'_2}{\sqrt{s^2 - x'^2_{,2}}} = 2(\arcsin 1 - \arcsin 0) = \pi. \qquad (1.31)$$

Hence, it follows from (1.29) that

$$G^{\lambda\langle 1,1\rangle}(x_{*0} - x_0, \{x_{*1} - x_1\}) = \frac{1}{2} H(\lambda(x_0 - x_{*0}) - |x_{*1} - x_1|) \text{ for each } \lambda \in \{-1,1\}. \qquad (1.32)$$

In either specific case, $n=2$ or $n=1$, retarded and advanced solutions of equation (1.18) can be calculated with the help of equations (1.20), (1.27), and (1.32) in order to demonstrate their *dispersive* character.

It is noteworthy that in Morse and Feshbach [1953, section 7.3], the Green function of a wave equation in $\overline{T} \times \overline{E}_3$, denoted by '$G$', is related as $G = 4\pi G^{\langle 1,3\rangle}$ to the Green function $G^{\langle 1,3\rangle}$, satisfying equation (1.20) at $n=3$, and that $x=x_1$, $y=x_2$, $z=x_3$, and $u=H$. Also it follows from (1.28) that

$$1 - H(x) = H(-x). \qquad (1.33)$$

With allowance of the above said and also with allowance of some other obvious correspondences between the pertinent notations of this paper and those of Morse and Feshbach (M&F), equation (1.23'), i.e. equation (1.27) at $\lambda=-1$, coincides with equation (7.3.15) of M&F at $c=1$, while equation (1.32) at $\lambda=-1$ coincides with equation (7.3.16) of M&F at $c=1$.

**Comment 1.2.** In analogy with Comment 1.1, it can be demonstrated by straightforward calculations that the Green functional form, given by equation (1.27) or (1.32), satisfies the instance of equation (1.20) at $n=2$ with $\langle x_0, \overline{x}_{[1,2]}\rangle$ and $\langle x_{*0}, \overline{x}_{*[1,2]}\rangle$ exchanged or the instance of equation (1.20) at $n=1$ with $\langle x_0, x_1\rangle$ and $\langle x_{*0}, x_{*1}\rangle$ exchanged, respectively. This is done below in the case of $n=1$.

By (1.2), (1.5)–(1.7), and (1.10)–(1.14) at $n=3$, it follows from (1.32) that

$$\frac{\partial^2 G^{\lambda\langle 1,1\rangle}(x_{*0}-x_0,\{x_{*1}-x_1\})}{\partial x_1^2} = \frac{1}{2}\frac{\partial^2 H(\lambda(x_0-x_{*0})-|x_{*1}-x_1|)}{\partial x_1^2}$$

$$= \frac{1}{2}\frac{\partial}{\partial x_1}\left[\frac{\partial|x_{*1}-x_1|}{\partial x_1}\frac{\partial H(\lambda(x_0-x_{*0})-|x_{*1}-x_1|)}{\partial|x_{*1}-x_1|}\right] \quad (1.32_1)$$

$$= -\frac{1}{2}\frac{\partial}{\partial x_1}\left[\frac{\partial|x_{*1}-x_1|}{\partial x_1}\delta(\lambda(x_0-x_{*0})-|x_{*1}-x_1|)\right].$$

Since

$$\frac{1}{2}\frac{\partial^2|x_{*1}-x_1|}{\partial x_1^2} = \delta(x_{*1}-x_1), \quad \left(\frac{\partial|x_{*1}-x_1|}{\partial x_1}\right)^2 = 1, \quad (1.32_2)$$

therefore the train ($1.32_1$) eloped further thus

$$\frac{\partial^2 G^{\lambda\langle 1,1\rangle}(x_{*0}-x_0,\{x_{*1}-x_1\})}{\partial x_1^2} == -\frac{1}{2}\frac{\partial^2|x_{*1}-x_1|}{\partial x_1^2}\delta(\lambda(x_0-x_{*0})-|x_{*1}-x_1|)$$

$$-\frac{1}{2}\left(\frac{\partial|x_{*1}-x_1|}{\partial x_1}\right)^2\frac{\partial\delta(\lambda(x_0-x_{*0})-|x_{*1}-x_1|)}{\partial|x_{*1}-x_1|} \quad (1.32_3)$$

$$= -\delta(x_{*1}-x_1)\delta(\lambda(x_0-x_{*0})-|x_{*1}-x_1|) - \frac{1}{2}\frac{\partial\delta(\lambda(x_0-x_{*0})-|x_{*1}-x_1|)}{\partial|x_{*1}-x_1|}$$

$$= -\delta(x_{*1}-x_1)\delta(x_0-x_{*0}) - \frac{1}{2}\frac{\partial\delta(\lambda(x_0-x_{*0})-|x_{*1}-x_1|)}{\partial|x_{*1}-x_1|}.$$

At the same time,

$$\frac{\partial^2 G^{\lambda\langle 1,1\rangle}(x_{*0}-x_0,\{x_{*1}-x_1\})}{\partial x_{*0}^2} = \frac{1}{2}\frac{\partial^2 H(\lambda(x_0-x_{*0})-|x_{*1}-x_1|)}{\partial x_{*0}^2}$$

$$= \frac{1}{2}\left[\frac{\partial\lambda(x_0-x_{*0})}{\partial x_{*0}}\right]^2\frac{\partial^2 H(\lambda(x_0-x_{*0})-|x_{*1}-x_1|)}{\partial(|x_{*1}-x_1|)^2} \quad (1.32_4)$$

$$= -\frac{1}{2}\frac{\partial\delta(\lambda(x_0-x_{*0})-|x_{*1}-x_1|)}{\partial|x_{*1}-x_1|}.$$

Hence, for each $\lambda \in \{-1,1\}$:

$$\left(\frac{\partial^2}{\partial x_1^2} - \frac{\partial^2}{\partial x_{*0}^2}\right)G^{\lambda\langle 1,1\rangle}(x_{*0}-x_0,\{x_{*1}-x_1\}) = -\delta(x_{*1}-x_1)\delta(x_0-x_{*0}) \quad (1.20c)$$

for each $x_{*0}-x_0 \in \overline{T}$ and each $x_{*1}-x_1 \in \overline{E}_1$,

as expected.•

**Comment 1.3.** The above method of calculation of $G^{\lambda\langle 1,2\rangle}$ and $G^{\lambda\langle 1,1\rangle}$ in terms of $G^{\lambda\langle 1,3\rangle}$ has the following implications. Given $\lambda \in \{-1,1\}$, equation (1.19) is equivalent to

$$\psi^{\langle 1,3\rangle}(x_{*0},\overline{x}_{*[1,3]}) = \int\limits_{\overline{E}_3}\int\limits_{-\infty}^{\infty} G^{\lambda\langle 1,3\rangle}(x_{*0}-x_0,\overline{x}_{*[1,3]}-\overline{x}_{[1,3]})f^{\langle 1,3\rangle}(x_0,\overline{x}_{[1,3]})dx_0 dv_{\{3\}}$$

$$= \frac{1}{4\pi}\int\limits_{\overline{E}_3}\int\limits_{-\infty}^{\infty}\frac{\delta(x_{*0}-x_0+\lambda x'_{[1,3]})f^{\langle 1,3\rangle}(x_0,\overline{x}'_{[1,3]})}{x'_{[1,n]}}dx_0 dv_{\{3\}} \qquad (1.34)$$

$$= \frac{1}{4\pi}\int\limits_{\overline{E}_3}\frac{f^{\langle 1,3\rangle}(x_{*0}+\lambda x'_{[1,3]},\overline{x}'_{[1,3]})}{x'_{[1,3]}}dv_{\{3\}},$$

which is (1.21) at $n=3$ subject to (1.22). Therefore, (i) if $f^{\langle 1,3\rangle}(x_0,\overline{x}_{[1,3]})$ is independent of $x_3$, i.e. if

$$f^{\langle 1,3\rangle}(x_0,\overline{x}_{[1,3]}) = f^{\langle 1,2\rangle}(x_0,\overline{x}_{[1,2]}), \qquad (1.35)$$

then (1.34) automatically turns into

$$\psi^{\langle 1,2\rangle}(x_{*0},\overline{x}_{*[1,2]}) = \int\limits_{\overline{E}_2}\int\limits_{-\infty}^{\infty} G^{\lambda\langle 1,2\rangle}(x_{*0}-x_0,\overline{x}_{*[1,2]}-\overline{x}_{[1,2]})f^{\langle 1,2\rangle}(x_0,\overline{x}_{[1,2]})dx_0 dv_{\{2\}} \qquad (1.36)$$

subject (1.27); (ii) if $f^{\langle 1,3\rangle}(x_0,\overline{x}_{[1,3]})$ is independent of $x_2$ and $x_3$, i.e. if

$$f^{\langle 1,3\rangle}(x_0,\overline{x}_{[1,3]}) = f^{\langle 1,1\rangle}(x_0,\{x_1\}), \qquad (1.37)$$

then (1.34) automatically turns into

$$\psi^{\langle 1,1\rangle}(x_{*0},\{x_{*1}\}) = \int\limits_{-\infty}^{\infty}\int\limits_{-\infty}^{\infty} G^{\lambda\langle 1,1\rangle}(x_{*0}-x_0,\{x_{*1}-x_1\})f^{\langle 1,1\rangle}(x_0,\{x_1\})dx_0 dx_1 \qquad (1.38)$$

subject (1.32).

Under definition (1.18), equation (1.34) turns into identity (tautology) (1.15a), whereas equations (1.36) and (1.38) turn into the identities:

$$\psi^{\langle 1,2\rangle}(x_{*0},\overline{x}_{*[1,2]})$$
$$= \int\limits_{\overline{E}_2}\int\limits_{-\infty}^{\infty} G^{\lambda\langle 1,2\rangle}(x_{*0}-x_0,\overline{x}_{*[1,2]}-\overline{x}_{[1,2]})\left(\sum_{i=1}^{n}\frac{\partial^2}{\partial x_i^2}-\frac{\partial^2}{\partial x_0^2}\right)\psi^{\langle 1,2\rangle}(x_0,\overline{x}_{[1,2]})dx_0 dv_{\{2\}}, \qquad (1.39)$$

$$\psi^{\langle 1,1\rangle}(x_{*0},\{x_{*1}\})$$
$$= \int\limits_{-\infty}^{\infty}\int\limits_{-\infty}^{\infty} G^{\lambda\langle 1,1\rangle}(x_{*0}-x_0,\{x_{*1}-x_1\})\left(\frac{\partial^2}{\partial x_1^2}-\frac{\partial^2}{\partial x_0^2}\right)\psi^{\langle 1,1\rangle}(x_0,\{x_1\})dx_0 dx_1. \qquad (1.40)$$

In this case, however, the functional forms $\psi^{\langle 1,2\rangle}(x_{*0},\overline{x}_{*[1,2]})$ and $\psi^{\langle 1,1\rangle}(x_{*0},\{x_{*1}\})$, occurring in (1.39) and (1.40), are *semantically different* from the homographic functional forms, occurring in (1.15b) and (1.15c), for the following reason. Equation (1.15) at each given a natural number $n \geq 1$ is deduced from the SIDF for the ball $\overline{B}_n^c(r_2, \overline{x}_{*[1,n]})$ under Axiom 7.1, according to which, given $\lambda \in \{-1,1\}$, given a natural number $n \geq 1$, when $r_2$ becomes large enough ($r_2 \to \infty$), there is a strictly positive real number $r_* \in (r_1, r_2)$ (to be restricted more precisely from below) such that for each $m \in \{0,1,2\}$, there are two strictly positive real numbers: $p_{*n,m} > m$, and $C_{n,m} > 0$ such that

$$\left|\frac{\partial^m \psi^{\langle 1,n \rangle}(x_0, \overline{x}_{[1,n]})}{\partial x_i^m}\right|_{x_0 = x_{*0} + \lambda x'_{[1,n]}} \leq \frac{C_{n,m}}{\left|\overline{x}_{[1,n]} - \overline{x}_{*[1,n]}\right|^{p_{*n,m}}} \tag{1.41}$$

for each $i \in \omega_{0,n}$ and every $\left|\overline{x}_{[1,n]} - \overline{x}_{*[1,n]}\right| > r_*$;

the inequality $\left|\overline{x}_{[1,n]} - \overline{x}_{*[1,n]}\right| > r_*$ is equivalent to the relation $\overline{x}_{[1,n]} \in \overline{E}_n - \overline{B}_n^c(r_*, \overline{x}_{*[1,n]})$. Under definition (1.18), given $\lambda \in \{-1,1\}$, given a natural number $n \geq 1$, it follows from (1.41) that

$$\left|f^{\langle 1,n \rangle}(x_{*0} + \lambda x'_{[1,n]}, \overline{x}_{[1,n]})\right| \leq \frac{(n+1)C_{n,2}}{\left|\overline{x}_{[1,n]} - \overline{x}_{*[1,n]}\right|^{p_{*n,2}}} \text{ for each } \left|\overline{x}_{[1,n]} - \overline{x}_{*[1,n]}\right| > r_*. \tag{1.42}$$

In this case, either one of the conditions (1.35) and (1.37) is incompatible with (1.42) at $n=3$.

On the other hand, the following formal computational procedure is possible. Let us assume that, besides (1.42) at $n=3$,

$$f^{\langle 1,3 \rangle}(x_0, \overline{x}_{[1,3]}) = f^{\langle 1,2 \rangle}(x_0, \overline{x}_{[1,2]}) \text{ for each } x_3 \in [x_{*3} - r_*, x_{*3} + r_*] \tag{1.35'}$$

instead of (1.35) or

$$f^{\langle 1,3 \rangle}(x_0, \overline{x}_{[1,3]}) = f^{\langle 1,1 \rangle}(x_0, \{x_1\})$$
for each $x_2 \in [x_{*2} - r_*, x_{*2} + r_*]$ and each $x_3 \in [x_{*3} - r_*, x_{*3} + r_*]$ \tag{1.37'}

instead of (1.37). In this case, we arrive at (1.34) subject to (1.35') or subject to (1.37') respectively. The limiting transition $r_* \to \infty$ in each one of the two homographic relations (1.34) thus obtained results respectively in (1.36) and (1.38) again, which again turn into relations (1.39) and (1.40) under definition (1.18). This time, however, the functional forms $\psi^{\langle 1,2 \rangle}(x_0, \overline{x}_{*[1,2]})$ and $\psi^{\langle 1,1 \rangle}(x_{*0}, \{x_{*1}\})$, occurring in (1.39) and (1.40), are *semantically the same* as those occurring in (1.15b) and (1.15c). Consequently, the expressions on the right-hand sides of equations (1.15b) and (1.39) or those of equations (1.15c) and (1.40) should *semantically be equal* as well. Whether or not this is the case can be verified only by the pertinent straightforward calculations (transformations), which seem however to be extremely cumbersome and time consuming.●

Given a natural number $n \geq 1$, let

$$G_\omega^{\langle 1,n \rangle}\left(\overline{x}_{*[1,n]} - \overline{x}_{[1,n]}\right) \equiv \int_{-\infty}^{\infty} G^{\langle 1,n \rangle}\left(x_{*0} - x_0, \overline{x}_{*[1,n]} - \overline{x}_{[1,n]}\right) e^{-i\omega(x_{*0} - x_0)} d(x_{*0} - x_0), \quad (1.43)$$

$$\psi_\omega^{\langle 1,n \rangle}\left(\overline{x}_{*[1,n]}\right) \equiv \int_{-\infty}^{\infty} \psi^{\langle 1,n \rangle}\left(x_{*0}, \overline{x}_{*[1,n]}\right) e^{-i\omega x_{*0}} dx_{*0}, \quad (1.44)$$

$$f_\omega^{\langle 1,n \rangle}\left(\overline{x}_{[1,n]}\right) \equiv \int_{-\infty}^{\infty} f^{\langle 1,n \rangle}\left(x_0, \overline{x}_{[1,n]}\right) e^{-i\omega x_0} dx_0, \quad (1.45)$$

so that conversely

$$G^{\langle 1,n \rangle}\left(x_{*0} - x_0, \overline{x}_{*[1,n]} - \overline{x}_{[1,n]}\right) = \frac{1}{2\pi} \int_{-\infty}^{\infty} G_\omega^{\langle 1,n \rangle}\left(\overline{x}_{*[1,n]} - \overline{x}_{[1,n]}\right) e^{i\omega(x_{*0} - x_0)} d\omega, \quad (1.46)$$

$$\psi^{\langle 1,n \rangle}\left(x_{*0}, \overline{x}_{*[1,n]}\right) = \frac{1}{2\pi} \int_{-\infty}^{\infty} \psi_\omega^{\langle 1,n \rangle}\left(\overline{x}_{*[1,n]}\right) e^{i\omega x_{*0}} dx_{*0}, \quad (1.47)$$

$$f^{\langle 1,n \rangle}\left(x_0, \overline{x}_{[1,n]}\right) = \frac{1}{2\pi} \int_{-\infty}^{\infty} f_\omega^{\langle 1,n \rangle}\left(\overline{x}_{[1,n]}\right) e^{i\omega x_{*0}} dx_{*0}. \quad (1.48)$$

Therefore, equation (1.21) becomes

$$\frac{1}{2\pi} \int_{-\infty}^{\infty} \psi_\omega^{\langle 1,n \rangle}\left(\overline{x}_{*[1,n]}\right) e^{i\omega x_{*0}} d\omega$$

$$= \frac{1}{4\pi^2} \int_{\overline{E}_n} \int_{-\infty}^{\infty} \int_{-\infty}^{\infty} \int_{-\infty}^{\infty} G_\omega^{\langle 1,n \rangle}\left(\overline{x}_{*[1,n]} - \overline{x}_{[1,n]}\right) f_{\omega'}^{\langle 1,n \rangle}\left(\overline{x}_{[1,n]}\right) e^{i\omega(x_{*0} - x_0) + i\omega' x_0} d\omega d\omega' dx_0 dv_{\{n\}} \quad (1.49)$$

$$= \frac{1}{2\pi} \int_{\overline{E}_n} \int_{-\infty}^{\infty} G_\omega^{\langle 1,n \rangle}\left(\overline{x}_{*[1,n]} - \overline{x}_{[1,n]}\right) f_\omega^{\langle 1,n \rangle}\left(\overline{x}_{[1,n]}\right) e^{i\omega x_{*0}} d\omega dv_{\{n\}},$$

because

$$\frac{1}{2\pi} \int_{-\infty}^{\infty} e^{i(\omega' - \omega) x_0} dx_0 = \delta(\omega' - \omega). \quad (1.50)$$

Equation (1.49) implies that

$$\psi_\omega^{\langle 1,n \rangle}\left(\overline{x}_{*[1,n]}\right) = \int_{\overline{E}_n} G_\omega^{\langle 1,n \rangle}\left(\overline{x}_{*[1,n]} - \overline{x}_{[1,n]}\right) f_\omega^{\langle 1,n \rangle}\left(\overline{x}_{[1,n]}\right) dv_{\{n\}}. \quad (1.51)$$

By (1.22), equation (1.43) at $n=3$ becomes

$$G_\omega^{\langle 1,3 \rangle}\left(\overline{x}_{*[1,3]} - \overline{x}_{[1,3]}\right) \equiv \frac{1}{4\pi} \int_{-\infty}^{\infty} \frac{\delta\left(x_{*0} - x_0 + \lambda x'_{[1,3]}\right)}{\left|\overline{x}_{*[1,3]} - \overline{x}_{[1,3]}\right|} e^{-i\omega(x_{*0} - x_0)} d(x_{*0} - x_0)$$

$$= \frac{1}{4\pi \left|\overline{x}_{*[1,3]} - \overline{x}_{[1,3]}\right|} e^{i\lambda\omega \left|\overline{x}_{*[1,3]} - \overline{x}_{[1,3]}\right|}. \quad (1.43a)$$

Therefore, equation (1.51) at $n=3$ turns into

$$\psi_\omega^{\langle 1,n\rangle}(\bar{x}_{*[1,3]}) = \frac{1}{4\pi}\int_{\bar{E}_3}\frac{1}{|\bar{x}_{*[1,3]} - \bar{x}_{[1,3]}|}e^{i\lambda\omega|\bar{x}_{*[1,3]} - \bar{x}_{[1,3]}|}dv_{\{3\}}. \tag{1.51a}$$

By way of comparison with (1.43a), let us calculate $G_\omega^{\lambda\langle 1,1\rangle}(x_{*1} - x_1)$ for each $\lambda \in \{-1,1\}$, which is defined by the pertinent variant of (1.43) at $n=1$ subject to (1.32), i.e. by

$$G_\omega^{\lambda\langle 1,1\rangle}(x_{*1} - x_1) \equiv \int_{-\infty}^{\infty} G^{\lambda\langle 1,1\rangle}(x_{*0} - x_0, x_{*1} - x_1)e^{-i\omega(x_{*0}-x_0)}d(x_{*0} - x_0)$$
$$= \frac{1}{2}\int_{-\infty}^{\infty} H(\lambda(x_0 - x_{*0}) - |x_{*1} - x_1|)e^{-i\omega(x_{*0}-x_0)}d(x_{*0} - x_0). \tag{1.52}$$

Setting $\xi \equiv \lambda(x_0 - x_{*0})$, so that $x_{*0} - x_0 = -\lambda\xi$ and $d(x_{*0} - x_0) = -\lambda d\xi$, equation (1.52) can be developed thus:

$$G_\omega^{\lambda\langle 1,1\rangle}(x_{*1} - x_1) = -\frac{\lambda}{2}\int_{\lambda\infty}^{-\lambda\infty} H(\xi - |x_{*1} - x_1|)e^{i\lambda\omega\xi}d\xi = \frac{1}{2}\int_{-\infty}^{\infty} H(\xi - |x_{*1} - x_1|)e^{i\lambda\omega\xi}d\xi. \tag{1.53}$$

In order to calculate the last integral, we shall replace the Heaviside function $H$, defined by (1.18), with the function $H^{(\gamma)}$, defined as:

$$H^{(\gamma)}(x) \equiv \begin{cases} e^{-\gamma x} & \text{if } x > 0 \\ 0 & \text{if } x < 0 \end{cases}, \tag{1.28a}$$

the understanding being that $\gamma$ should be taken to $+0$ in the final result. That is to say, we conventionally assume that

$$G_\omega^{\lambda\langle 1,1\rangle}(x_{*1} - x_1) = \frac{1}{2}\lim_{\gamma\to +0}\int_{-\infty}^{\infty} H^{(\gamma)}(\xi - |x_{*1} - x_1|)e^{i\lambda\omega\xi}d\xi$$
$$= \frac{1}{2}\lim_{\gamma\to +0}\int_{|x_{*1}-x_1|}^{\infty} e^{i\lambda\omega\xi - \gamma(\xi - |x_{*1}-x_1|)}d\xi = -\frac{1}{2}e^{i\lambda\omega|x_{*1}-x_1|}\lim_{\gamma\to +0}\frac{1}{i\lambda\omega - \gamma} = \frac{i\lambda}{2}e^{i\lambda\omega|x_{*1}-x_1|}\frac{1}{\omega + i\lambda 0}. \tag{1.53_1}$$

In this case, for any continuous functional form '$f(\omega)$' of a complex variable '$\omega$', we have, e.g., that:

$$\int_{-\infty}^{\infty}\frac{f(\omega)}{\omega - i0}d\omega = \lim_{\rho\to +0}\left[\int_{-\infty}^{-\rho}\frac{f(\omega)}{\omega}d\omega + \int_{\pi}^{0}\frac{f(\rho e^{i\varphi})}{\rho e^{i\varphi}}i\rho e^{i\varphi}d\varphi + \int_{\rho}^{\infty}\frac{f(\omega)}{\omega}d\omega\right]$$
$$= if(0)\int_{\pi}^{0}d\varphi + P\int_{-\infty}^{\infty}\frac{f(\omega)}{\omega}d\omega = -i\pi f(0) + P\int_{-\infty}^{\infty}\frac{f(\omega)}{\omega}d\omega, \tag{1.55}$$

so that

$$\frac{1}{\omega - i0} = -i\pi\delta(\omega) + P\frac{1}{\omega}, \quad \frac{1}{\omega + i0} = \left(\frac{1}{\omega - i0}\right)^* = i\pi\delta(\omega) + P\frac{1}{\omega} \tag{1.56}$$

or equivalently

$$\frac{1}{\omega + i\lambda 0} = i\lambda\pi\delta(\omega) + P\frac{1}{\omega} \text{ for each } \lambda \in \{-1,1\}. \tag{1.57}$$

Therefore, (1.53$_1$) can be developed further thus:

$$G_\omega^{\lambda\langle 1,1\rangle}(x_{*1} - x_1) = \frac{i\lambda}{2} e^{i\lambda\omega|x_{*1}-x_1|}\left[i\lambda\pi\delta(\omega) + P\frac{1}{\omega}\right]$$
$$= \frac{1}{2} e^{i\lambda\omega|x_{*1}-x_1|}\left[-\pi\delta(\omega) + i\lambda P\frac{1}{\omega}\right] \text{ for each } \lambda \in \{-1,1\}. \tag{1.58}$$

Hence, given $\lambda \in \{-1,1\}$, (1.51) at $n=1$ becomes

$$\psi_\omega^{\langle 1,1\rangle}(x_{*1}) = \int_{-\infty}^{\infty} G_\omega^{\lambda\langle 1,1\rangle}(x_{*1} - x_1) f_\omega^{\langle 1,1\rangle}(x_1) dx_1$$
$$= \frac{1}{2}\int_{-\infty}^{\infty}\left[-\pi\delta(\omega) + i\lambda P\frac{1}{\omega}\right] e^{i\lambda\omega|x_{*1}-x_1|} f_\omega^{\langle 1,1\rangle}(x_1) dx_1 \tag{1.51c}$$
$$= \frac{1}{2}\left[-\pi\delta(\omega)\int_{-\infty}^{\infty} f_0^{\langle 1,1\rangle}(x_1) dx_1 + i\lambda P\frac{1}{\omega}\int_{-\infty}^{\infty} e^{i\lambda\omega|x_{*1}-x_1|} f_\omega^{\langle 1,1\rangle}(x_1) dx_1\right].$$

In accordance with (1.51a), equation (1.19) is called a forced *non-dispersed* solution of the inhomogeneous wave equation (1.18) at $n=3$ in $\overline{T} \times \overline{E}_3$, – a *retarded* one if $\lambda=-1$ and an *advanced* one if $\lambda=1$. By contrast, in accordance with (1.51c), equation (1.38) is called a forced *dispersed* solution of the inhomogeneous wave equation (1.18) at $n=1$ in $\overline{T} \times \overline{E}_1$, – but again a *retarded* one if $\lambda=-1$ and an *advanced* one if $\lambda=1$. Since the second item on the right-hand side of equation (1.15) vanishes if and only if $n=3$, the inhomogeneous wave equation (1.18) at $n\neq 3$ has *no non-dispersive solution*, either retarded or advanced. In other words, *no o non-dispersive causality principle* is satisfied in this general case. Consequently, no special relativity exists in any affine manifold $\dot{T}(R) \times \dot{E}_n(R)$ at $n\neq 3$, so that the physical properties of the affine manifold $\dot{T}(R) \times \dot{E}_3(R)$ are unique. Therefore, *a three-dimensional affine Euclidean space is the only possible receptacle of nature*. In this connection, it is noteworthy that some physical phenomena and particularly some wave motions are treated (approximated) as one-dimensional ones and some as two-dimensional ones. However, all such phenomena occur in the three-dimensional affine Euclidean space, and neither in a one-dimensional one nor in a two-dimensional one. Thus, the SIDF (1.15) is a manifestation of unique symmetry of $\dot{E}_3(R)$, which is associated with its 3-dimensionality and owing to which $\dot{E}_3(R)$ is coordinated with $\dot{T}(R)$ so as to serve a receptacle of nature.

Two different definitions of an affine space can be found, e.g., in Mac Lane and Birkhoff [1967, p. 420] and in Schwartz [1967, Part I, Chapter III, §1]. All necessary properties of $\dot{E}_n(R)$ are explicated in subsections 2.1 and 2.2 of this exposition. A new convenient concise rigorous system of axioms of all relevant algebraic systems underlying the SIDF at any $n \geq 1$, including a system of axioms of an *n*-dimensional Euclidean affine space will be given separately in the second part of this work to be published elsewhere The *axiomatic algebraico-functional theory* (*AAFT*) *of* $\dot{E}_n(R)$ that is developed in that part of the work can be regarded generally as an *underlying discipline of differential and integral calculi* and hence as a formal interface between *any hypothetical physical processes in* $\dot{E}_n(R) \times \dot{T}(R)$ and *their mathematical descriptions in* $\overline{E}_n(R) \times \overline{T}(R)$.

## 2. Basic notation

### 2.1. General definitions

**Definition 2.1.** 1) The signs $\stackrel{.}{=}$ and $\stackrel{\cdot}{=}$ are indiscriminately called the asymmetric, or one-sided, equality signs by definition or, discriminately, the rightward equality sign by definition and the leftward equality sign by definition respectively. A binary figure, in which either sign $\stackrel{.}{=}$ or $\stackrel{\cdot}{=}$ is used assertively, is called a formal binary asymmetric synonymic definition (FBASD). In making a FBASD, at the head of an arrow we shall write the material definiens – the substantive (substance-valued expression), which is already known either from a previous definition or from another source; at the base of the arrow we shall write the material definiendum – the new substantive, which is being introduced by the definition and which is designed to be used instead of or interchangeably with the definiens in the scope of the FBASD. Therefore, the sign $\stackrel{.}{=}$ is rendered into ordinary language thus: "is to stand as a synonym for" or straightforwardly "is the synonymous definiendum of", and $\stackrel{\cdot}{=}$ thus: "can be used instead of interchangeably with" or straightforwardly "is the synonymous definiens of". The [material] definiendum and [material] definiens of a FBASD are indiscriminately called the terms of the definition. Neither the definiendum nor the definiens of an FBASD should involve any function symbols, particularly any outermost (enclosing) quotation marks, that are not their constituent parts and that are therefore used but not mentioned with the following proviso. If it is necessary to indicate the integrity of the definiendum or of the definiens then that term of the definition can be enclosed in square brackets as metalinguistic punctuation marks, which do not, by definition, belong to the enclosed term and which are therefore used but not mentioned. In the scope of a FBASD, which does not include the FBASD itself, tokens of the terms of the FBASD can be related by the ordinary reflexive, symmetric, and transitive sign of equality =. In contrast to =, either sign $\overset{\rightharpoonup}{=}$ or $\overset{\leftharpoonup}{=}$ is transitive, but not reflexive and not symmetric.

2) In order to state formally that two old or two new substantives are to be used interchangeably (synonymously), we shall write the substantives, without any quotation marks that are not their constituent parts, in either order on both sides of the two-sided sign $\overset{\leftrightarrow}{=}$. Such a relation is called a *formal binary symmetric synonymic definition* (*FBSSD*), whereas the sign $\overset{\leftrightarrow}{=}$ is called the *symmetric*, or *two-sided*, *equality sign by definition*. In this case, $\overset{\leftrightarrow}{=}$ is read as "*is to be concurrent to*" or, alternatively, "— $\overset{\leftrightarrow}{=}$ …" is read as "— *and … are to be concurrent*" or as "— *and … are to be synonyms*", where alike ellipses should be replaced alike. In the scope of an FBSSD, tokens of the terms of the FBSSD can be related by the ordinary sign of equality =.

3) In stating synonymic definitions of substantives, the arrows →, ←, and ↔ can be used instead of $\overset{\rightharpoonup}{=}$, $\overset{\leftharpoonup}{=}$, and $\overset{\leftrightarrow}{=}$ respectively, the understanding being that the arrows are general definition signs, which can apply to relations and not only to substantives.●

**Definition 2.2**. $\{x|x \in A \text{ and } P(x)\}$ is the class (or particularly set) of elements of the class (set) *A* having the property (predicate) *P*.●

**Definition 2.3**. If *A* and *B* are sets then *A−B*, called the *difference* between *A* and *B* and also the *relative complement* of *B* in *A*, is the set of all those elements of *A*, which are not elements of *B*, i.e. $A - B \overset{\leftrightarrow}{=} \{x|x \in A \text{ and } x \notin B\}$.●

**Definition 2.4.** 1) '$\omega_0$' *denotes*, i.e. $\omega_0$ *is*, the set of all *natural numbers* from 0 to infinity.

2) Given $n \in \omega_0$,
$$\omega_n \overset{\leftrightarrow}{=} \{i|i \in \omega_0 \text{ and } i \geq n\}, \tag{2.1}$$
i.e. '$\omega_1$', '$\omega_2$', etc *denote* the sets of natural numbers from 1, 2, etc respectively to infinity.

3) Given $m \in \omega_0$, given $n \in \omega_m$,
$$\omega_{m,n} \overset{\leftrightarrow}{=} \{i|i \in \omega_0 \text{ and } n \geq i \geq m\}, \tag{2.2}$$
i.e. '$\omega_{m,n}$' *denotes* the set of natural numbers from a given number *m* to another given number *n* subject to *n*≥*m*. It is understood that
$$\omega_{m,m} = \{m\}, \; \omega_{m,\infty} = \omega_m, \; \omega_{m,n} = \varnothing \text{ if } m \geq n; \tag{2.3}$$

$\varnothing$ is the empty set.•

**Definition 2.5.** 1) '$I_{-\infty,\infty}$' denotes, i.e. $I_{-\infty,\infty}$ is, the set of all *natural integers* (*natural integral numbers*) – strictly positive, strictly negative, and zero.

2) Given $n \in I_{-\infty,\infty}$,

$$I_{n,\infty} = I_{\infty,n} \equiv \{i | i \in I_{-\infty,\infty} \text{ and } i \geq n\}, \tag{2.4}$$

$$I_{-\infty,n} = I_{n,-\infty} \equiv \{i | i \in I_{-\infty,\infty} \text{ and } i \leq n\}, \tag{2.5}$$

i.e. $I_{n,\infty}$ or $I_{\infty,n}$ is the set of all natural integers greater than or equal to $n$, and $I_{-\infty,n}$ or $I_{n,-\infty}$ is the set of all natural integers less than or equal to $n$.

3) Given $m \in I_{-\infty,\infty}$, given $n \in I_{m,\infty}$,

$$I_{m,n} \equiv \{i | i \in I_{-\infty,\infty} \text{ and } n \geq i \geq m\}, \tag{2.6}$$

i.e. $I_{m,n}$ is the set of all natural integers that are greater than or equal to $m$ and less than or equal to $n$.•

**Comment 2.1.** Definitions 2.4(1) and 2.5(1) are *explicative* ones. A theory of natural integers in particular, and a theory of any numbers (as rational, real, or complex ones) in general can consistently be deduced from the five Peano axioms, which are, in turn, theorems of an axiomatic set theory (see, e.g., Halmos [1960, pp. 46–53], Burrill [1967], Feferman [1964]).•

**Definition 2.6.**

1) $\boldsymbol{R}$ is the field of real numbers.

2) $R$ is the underlying set of $\boldsymbol{R}$, so that $R = (-\infty, \infty)$.

3) $\dot{\boldsymbol{E}}_n(\boldsymbol{R})$, or briefly $\dot{\boldsymbol{E}}_n$, is an $n$-dimensional Euclidean affine (point) space over $\boldsymbol{R}$, called also an Euclidean real affine $n$-space.

4) $\dot{E}_n$ is the underlying set of points of $\dot{\boldsymbol{E}}_n(\boldsymbol{R})$.

5) $\hat{\boldsymbol{E}}_n(\boldsymbol{R})$, or briefly $\hat{\boldsymbol{E}}_n$, is an $n$-dimensional Euclidean linear, or abstract vector, space over $\boldsymbol{R}$, adjoint of $\dot{\boldsymbol{E}}_n(\boldsymbol{R})$, called also an Euclidean real linear, or abstract vector, $n$-space.

6) $\hat{E}_n$ is the underlying set of real abstract vectors of $\hat{\boldsymbol{E}}_n(\boldsymbol{R})$.

7) $\overline{\boldsymbol{E}}_n(\boldsymbol{R})$, or briefly $\overline{\boldsymbol{E}}_n$, is an $n$-dimensional Euclidean arithmetical vector space over $\boldsymbol{R}$, isomorphic to $\hat{\boldsymbol{E}}_n(\boldsymbol{R})$, called also an Euclidean real arithmetical vector $n$-space.

8) $\overline{E}_n$ is the underlying set of real arithmetical vectors of $\overline{E}_n(\mathbf{R})$, so that

$$\overline{E}_1 \triangleq R^{1\times} \triangleq \{\langle x_1\rangle | x_1 \in R\} = \{\{x_1\} | x_1 \in R\} \neq R, \tag{2.7}$$

$$\overline{E}_n \triangleq R^{n\times} \triangleq \underbrace{R \times R \times ... \times R}_{n \text{ times } R} \triangleq R^{(n-1)\times} \times R$$
$$\triangleq \underbrace{[[...[R \times R] \times R] \times ...] \times R]}_{n-1} \times R \text{ for each } n \in \omega_2. \tag{2.8}$$

9) $\dot{x}_{\{n\}}$, $\dot{y}_{\{n\}}$, and $\dot{z}_{\{n\}}$ (e.g.) or briefly (whenever there is no danger of confusion) $\dot{x}$, $\dot{y}$, and $\dot{z}$ are arbitrary points in $\dot{E}_n$.

10) $\hat{x}_{\{n\}}$, $\hat{y}_{\{n\}}$, and $\hat{z}_{\{n\}}$ (e.g.) or briefly (whenever there is no danger of confusion) $\hat{x}$, $\hat{y}$, and $\hat{z}$ are arbitrary abstract vectors in $\hat{E}_n$.

11) $\overline{x}_{[1,n]}$, $\overline{y}_{[1,n]}$, and $\overline{z}_{[1,n]}$ (e.g.) or briefly (whenever there is no danger of confusion) $\overline{x}$, $\overline{y}$, and $\overline{z}$ are arbitrary arithmetical vectors in $\overline{E}_n$, so that, e.g.,

$$\overline{x}_{[1,1]} \triangleq \langle x_1 \rangle = \{x_1\} \in \overline{E}_1, \tag{2.9}$$

$$\overline{x}_{[1,2]} \triangleq \langle x_1, x_2 \rangle \triangleq \{\{x_1\},\{x_1,x_2\}\} \in \overline{E}_2, \tag{2.10}$$

$$\overline{x}_{[1,n]} \triangleq \langle x_i \rangle_{i \in \omega_{1,n}} \triangleq \langle x_1, x_2, ..., x_{n-1}, x_n \rangle \triangleq \langle \overline{x}_{[1,n-1]}, x_n \rangle$$
$$= \underbrace{\langle\langle ..\langle x_1, x_2 \rangle, x_3 \rangle, ..., x_{n-1} \rangle, x_n \rangle}_{n-1} \in \overline{E}_n. \tag{2.11}$$

An ordered *n*-tuple standing for an *n*-dimensional arithmetical vector can be written either with angle brackets or with round ones. However, if $x_1$ and $x_2$ are real numbers then the symbol '$(x_1, x_2)$' is ambiguous, for it may stand either for the ordered pair of those numbers in that order or for the open interval $(x_1, x_2)$. Therefore, in denoting ordered multiples, we shall use angle brackets and round brackets interchangeably, while in most general conceptual statements preference will be given to the former.

12) $\dot{X}_n$ (e.g.) is an arbitrary connected subset of $\dot{E}_n$.

13) $\hat{X}_n$ (e.g.) is an arbitrary connected subset of $\hat{E}_n$.

14) $\overline{X}_n$ (e.g.) is an arbitrary connected subset of $\overline{E}_n$.

15) $\dot{0}_{\{n\}}$ is a given (fixed) reference point in $\dot{E}_{\{n\}}$.

16) $\hat{0}_{\{n\}}$ is the null vector in $\hat{E}_n$.

17) $\overline{0}_{[1,n]}$ is the null vector in $\overline{E}_n$, so that

$$\overline{0}_{[1,1]} \equiv \langle 0 \rangle = \{0\}, \tag{2.12}$$

$$\overline{0}_{[1,n]} \equiv \langle \underbrace{0,...,0}_{n \text{ zeros}} \rangle \equiv \langle \overline{0}_{[1,n-1]},0 \rangle \equiv \underbrace{\langle \langle ... \langle 0,0 \rangle, 0 \rangle,... \rangle, 0 \rangle}_{n-1}. \tag{2.13}$$

18) • is the binary function (operation) of *inner*, or *scalar*, *multiplication* on $\hat{E}_n$.

19) $\overline{\hat{e}}_{[1,n]} \equiv \langle \hat{e}_1,...,\hat{e}_n \rangle \in \hat{E}_n^{n\times}$ is an *orthonormal* (*normal orthogonal*) *basis* in $\hat{E}_n(\mathbf{R})$, so that

$$\hat{e}_i \bullet \hat{e}_j = \delta_{ij} \text{ for each } i \in \omega_{1,n} \text{ and for each } j \in \omega_{1,n}, \tag{2.14}$$

where '$\delta_{ij}$' is the Kronecker delta-symbol.

20) $\overline{\overline{e}}_{[1,n][1,n]} \equiv \langle \overline{e}_{1[1,n]},...,\overline{e}_{n[1,n]} \rangle \in \overline{E}_{\{n\}}^{n\times}$ is an *orthonormal basis* in $\overline{E}_n(\mathbf{R})$, so that

$$\overline{e}_{i[1,n]} \bullet \overline{e}_{j[1,n]} = \delta_{ij} \text{ for each } i \in \omega_{1,n} \text{ and each } j \in \omega_{1,n}. \tag{2.15}$$

It is understood that

$$\overline{e}_{1[1,n]} \equiv \langle 1,0,0,...,0 \rangle, \overline{e}_{2[1,n]} \equiv \langle 0,1,0,...,0 \rangle, ..., \overline{e}_{n[1,n]} \equiv \langle 0,0,...,0,1 \rangle, \tag{2.16}$$

or, equivalently,

$$\overline{e}_{i[1,n]} \equiv \langle \delta_{i1},...,\delta_{in} \rangle \text{ for each } i \in \omega_{1,n}. \tag{2.17}$$

21) $c_{\{n\}} \equiv \langle \dot{0}_{\{n\}}, \overline{\hat{e}}_{[1,n]} \rangle \in \dot{E}_{\{n\}} \times \hat{E}_{\{n\}}^{n\times}$ is an orthonormal coordinate system in $\dot{E}_n(\mathbf{R})$. •

**Comment 2.2.** There are a number of kinds of connectedness of sets. In this treatise, we shall follow to the following definition of disconnected and connected sets (cf. Definition 8-29 in Apostol [1963, p. 178]). A set $S$ is said to be *disconnected* if there are two sets $S_1$ and $S_2$ such that $S=S_1 \cup S_2$ and $S_1 \cap S_2 = \varnothing$. A set $S$ is said to be *connected* if it is not disconnected.•

**Definition 2.7.** 1) Given $\dot{0}_{\{n\}} \in \dot{E}_n$, there are two mutually inverse bijections $\hat{V}_{\dot{0}_{\{n\}}} : \dot{E}_n \to \hat{E}_n$ and $\hat{V}_{\dot{0}_{\{n\}}}^{-1} : \hat{E}_n \to \dot{E}_n$, i.e. for each $\dot{x}_{\{n\}} \in \dot{E}_n$ there is exactly one $\hat{x}_{\{n\}} \in \hat{E}_n$ such that

$$\hat{x}_{\{n\}} = \hat{V}_{\dot{0}_{\{n\}}}(\dot{x}_{\{n\}}) \tag{2.18}$$

and conversely for each $\hat{x}_{\{n\}} \in \hat{E}_n$ there is exactly one $\dot{x}_{\{n\}} \in \dot{E}_n$ such that

$$\dot{x}_{\{n\}} = \hat{V}_{\dot{0}_{\{n\}}}^{-1}(\hat{x}_{\{n\}}). \tag{2.19}$$

The bijections $\hat{V}_{\dot{0}_{\{n\}}}$ and $\hat{V}_{\dot{0}_{\{n\}}}^{-1}$ are called the *abstract vectorization of* $\dot{E}_n(\mathbf{R})$ *relative to the point* $\dot{0}_{\{n\}}$ and the *pointillage of* $\hat{E}_n(\mathbf{R})$ *relative to the point* $\dot{0}_{\{n\}}$.

2) Given $\bar{\hat{e}}_{[1,n]} \equiv \langle \hat{e}_1,...,\hat{e}_n \rangle \in \hat{E}_n^{n\times}$, there are two mutually inverse bijections $\overline{C}_{\bar{\hat{e}}_{[1,n]}} : \hat{E}_n \to \overline{E}_n$ and $\overline{C}_{\bar{\hat{e}}_{[1,n]}}^{-1} : \overline{E}_n \to \hat{E}_n$, i.e. for each $\hat{x}_{\{n\}} \in \hat{E}_n$ there is exactly one $\bar{x}_{[1,n]} \in \overline{E}_n$ such that

$$\bar{x}_{[1,n]} = \overline{C}_{\bar{\hat{e}}_{[1,n]}}(\hat{x}_{\{n\}}) \tag{2.20}$$

and conversely for each $\bar{x}_{[1,n]} \in \overline{E}_n$ there is exactly one $\hat{x}_{\{n\}} \in \hat{E}_n$ such that

$$\hat{x}_{\{n\}} = \overline{C}_{\bar{\hat{e}}_{[1,n]}}^{-1}(\bar{x}_{[1,n]}). \tag{2.21}$$

In thic case, $\overline{C}_{\bar{\hat{e}}_{[1,n]}}$ is an *isomorphism from $\hat{E}_n(R)$ onto $\overline{E}_n(R)$*, which is called the *coordinatization of $\hat{E}_n(R)$ relative to the basis $\bar{\hat{e}}_{[1,n]}$*, while $\overline{C}_{\bar{\hat{e}}_{[1,n]}}^{-1}$ is an *isomorphism from $\overline{E}_n(R)$ onto $\hat{E}_n(R)$*, which is called the *abstract vectoriization of $\overline{E}_n(R)$ relative to the basis $\bar{\hat{e}}_{[1,n]}$*. The coordinates $x_1, x_2,..., x_{n-1}, x_n$ of the ordered *n*-tuple $\bar{x}_{[1,n]} \in \overline{E}_n$ are alternatively called the coordinates of the abstract vector $\hat{x}_{\{n\}} \in \hat{E}_n$ relative to the basis vectors $\hat{e}_1, \hat{e}_2,...., \hat{e}_{n-1}, \hat{e}_n$.

3) Given $c_{\{n\}} \equiv \langle \dot{0}_{\{n\}}, \bar{\hat{e}}_{[1,n]} \rangle$, there are two composite mutually inverse bijections $\overline{K}_{c_{\{n\}}} : \dot{E}_n \to \overline{E}_n$ and $\overline{K}_{c_{\{n\}}}^{-1} : \overline{E}_n \to \dot{E}_n$, subject to

$$\overline{K}_{c_{\{n\}}} \equiv \overline{C}_{\bar{\hat{e}}_{[1,n]}} \circ \hat{V}_{\dot{0}_{\{n\}}} \text{ and } \overline{K}_{c_{\{n\}}}^{-1} \equiv \hat{V}_{\dot{0}_{\{n\}}}^{-1} \circ \overline{C}_{\bar{\hat{e}}_{[1,n]}}^{-1}, \tag{2.22}$$

i.e. for each $\dot{x}_{\{n\}} \in \dot{E}_n$ there is exactly one $\bar{x}_{[1,n]} \in \overline{E}_n$ such that

$$\bar{x}_{[1,n]} = \overline{K}_{c_{\{n\}}}(\dot{x}_{\{n\}}) = \overline{C}_{\bar{\hat{e}}_{[1,n]}}(\hat{V}_{\dot{0}_{\{n\}}}(\dot{x}_{\{n\}})), \tag{2.23}$$

and, conversely, for each $\bar{x}_{[1,n]} \in \overline{E}_n$ there is exactly one $\dot{x}_{\{n\}} \in \dot{E}_n$ such that

$$\dot{x}_{\{n\}} = \overline{K}_{c_{\{n\}}}^{-1}(\bar{x}_{[1,n]}) = \hat{V}_{\dot{0}_{\{n\}}}^{-1}(\overline{C}_{\bar{\hat{e}}_{[1,n]}}^{-1}(\bar{x}_{[1,n]})). \tag{2.24}\bullet$$

**Comment 2.3.** It follows from Definition 2.7(3) that given $c_{\{n\}} \equiv \langle \dot{0}_{\{n\}}, \bar{\hat{e}}_{[1,n]} \rangle \in \dot{E}_n \times \hat{E}_n^{n\times}$, for each $\dot{X}_n \subseteq \dot{E}_n$ there is exactly one $\overline{X}_n \subseteq \overline{E}_n$ such that

$$\overline{X}_{\{n\}} = \overline{K}_{c_{\{n\}}}(\dot{X}_n) \equiv \{\bar{x}_{[1,n]} | \bar{x}_{[1,n]} = \overline{K}_{c_{\{n\}}}(\dot{x}_{\{n\}}) \text{ and } \dot{x}_{\{n\}} \in \dot{X}_n \}, \tag{2.25}$$

whence particularly

$$\overline{E}_n = \overline{K}_{c_{\{n\}}}(\dot{E}_n), \overline{E}_1 = \overline{K}_{c_{\{1\}}}(\dot{E}_1). \tag{2.26}\bullet$$

**Definition 2.8.** 1) A hypothetical time continuum that is associated with $\dot{E}_n(R)$ is a special interpretation of $\dot{E}_1(R)$, which will be denoted by '$\dot{T}(R)$' or briefly by '$\dot{T}$', so that the pertinent interpretations of

$$\hat{E}_1,\ \overline{E}_1,\ \dot{E}_1,\ \hat{E}_1,\ \overline{E}_1,\ \dot{0}_{\{1\}}\ \hat{0}_{\{1\}}\ \overline{0}_{[1,1]},\ \hat{e}_1,\ \overline{e}_{1[1,1]},\ c_{\{1\}} \qquad (2.27)$$

will be denoted respectively by

$$'\hat{T}',\ '\overline{T}',\ '\dot{T}',\ '\hat{T}',\ '\overline{T}',\ '\dot{\theta}',\ '\hat{\theta}',\ '\overline{\theta}',\ '\hat{\varepsilon}',\ '\overline{\varepsilon}',\ '\omega', \qquad (2.28)$$

the understanding being that

$$\overline{\theta} = \langle 0 \rangle = \{0\},\ \overline{\varepsilon} = \langle 1 \rangle = \{1\},\ \omega = \langle \dot{\theta}, \langle \overline{\varepsilon} \rangle \rangle. \qquad (2.29)$$

Mnemonically, the letter 'ω' as used above is the first letter of the Greek noun '*ωρολόγιον*' \orológion\ meaning *a clock*. At the same time, the pertinent interpretations of

$$\dot{x}_{\{1\}},\ \dot{y}_{\{1\}},\ \dot{z}_{\{1\}},\ \hat{x}_{\{1'\}},\ \hat{y}_{\{1\}},\ \hat{z}_{\{1\}},\ \overline{x}_{[1,1]},\ \overline{y}_{[1,1]},\ \overline{z}_{[1,1]} \qquad (2.30)$$

(e.g.) will be denoted respectively by

$$'\dot{\xi}',\ '\dot{\eta}',\ '\dot{\zeta}',\ '\hat{\xi}',\ '\hat{\eta}',\ '\hat{\zeta}',\ '\overline{\xi}',\ '\overline{\eta}',\ '\overline{\zeta}', \qquad (2.31)$$

the understanding being that

$$\overline{\xi} = \langle \xi \rangle \equiv \langle x_0 \rangle,\ \overline{\eta} = \langle \eta \rangle \equiv \langle y_0 \rangle,\ \overline{\zeta} = \langle \zeta \rangle \equiv \langle z_0 \rangle. \qquad (2.32)$$

## 2.2. Real-valued functions defined on $\dot{T} \times \dot{X}_n$, on $\overline{T} \times \overline{X}_n$, or on $R \times \overline{X}_n$

**Comment 2.4.** If exists, a hypothetical measurable time-dependent physical field occurring in an *n*-dimensional Euclidean real affine space $\dot{E}_n(R)$ should be described by one or more real-valued functions such as $\Psi_{\dot{T} \times \dot{E}_n}$, which is defined on the direct product $\dot{T} \times \dot{X}_n$, where $\dot{X}_n$ is a certain connected subset of $\dot{E}_n$. In what follows, it is shown that given coordinate systems $\omega$ in $\dot{T}$ and $c_{\{n\}}$ in $\dot{E}_n$, a function $\Psi_{\dot{T} \times \dot{E}_n}$ can be reduced to a certain real-valued function $\Psi_{\overline{T} \times \overline{E}_n}$ defined on the direct product $\overline{T} \times \overline{X}_n$, where $\overline{X}_n$ is the pertinent subset of $\overline{E}_n$. •

**Definition 2.9.** Given an Euclidean real affine space $\dot{E}_n(R)$ of a given dimension $n \in \omega_1$, let $\Psi_{\dot{E} \times \dot{E}_n}$ be a *real-valued* function from $\dot{T} \times \dot{X}_n$, subject to $\dot{X}_n \subseteq \dot{E}_n$, to $Y \subseteq R = (-\infty, \infty)$, i.e. a function whose *domain of departure* ($D_{\mathrm{dp}}$), *domain of arrival* ($D_\mathrm{a}$), *domain of definition* ($D_{\mathrm{df}}$), and *domain of variation* ($D_\mathrm{v}$) are $\dot{T} \times \dot{E}_n$, $R$, $\dot{T} \times \dot{X}_n$, and $Y$ respectively. Thus, symbolically,

$$\Psi_{\dot{T}\times\dot{E}_n}:\dot{T}\times\dot{X}_n\to Y, \tag{2.33}$$

so that

$$D_{\mathrm{dp}}\!\left(\Psi_{\dot{T}\times\dot{E}_n}\right)=\dot{T}\times\dot{E}_n,\ D_{\mathrm{a}}\!\left(\Psi_{\dot{T}\times\dot{E}_n}\right)=R,$$
$$D_{\mathrm{df}}\!\left(\Psi_{\dot{T}\times\dot{E}_n}\right)=\dot{T}\times\dot{X}_n\subseteq\dot{T}\times\dot{E}_n,\ D_{\mathrm{v}}\!\left(\Psi_{\dot{T}\times\dot{E}_n}\right)=Y\subseteq R. \tag{2.34}$$

Consequently, for each $\langle\dot{\xi},\dot{x}_{\{n\}}\rangle\in\dot{T}\times\dot{X}_n$ there is exactly one $y\in Y$ such that

$$y=\Psi_{\dot{T}\times\dot{E}_n}\!\left(\dot{\xi},\dot{x}_{\{n\}}\right). \tag{2.35}\bullet$$

**Definition 2.10.** By items 1 and 3 of Definition 2.7 and by Comment 2.3, it follows from Definition 2.9 that there is a composite function

$$\Psi_{\overline{T}\times\overline{E}_n}:\overline{T}\times\overline{X}_n\to Y \tag{2.36}$$

subject to

$$\Psi_{\overline{T}\times\overline{E}_n}\rightleftharpoons\Psi_{\dot{T}\times\dot{E}_n}\circ\left(\overline{K}_\omega^{\,-1}\times\overline{K}_{c_{\{n\}}}^{\,-1}\right), \tag{2.37}$$

so that

$$D_{\mathrm{dp}}\!\left(\Psi_{\overline{T}\times\overline{E}_n}\right)=\overline{T}\times\overline{E}_n,\ D_{\mathrm{a}}\!\left(\Psi_{\overline{T}\times\overline{E}_n}\right)=R,$$
$$D_{\mathrm{df}}\!\left(\Psi_{\overline{T}\times\overline{E}_n}\right)=\overline{T}\times\overline{X}_n\subseteq\overline{T}\times\overline{E}_n,\ D_{\mathrm{v}}\!\left(\Psi_{\overline{T}\times\overline{E}_n}\right)=Y\subseteq R. \tag{2.38}$$

Consequently, for each $\langle\overline{\xi},\overline{x}_{[1,n]}\rangle\in\overline{T}\times\overline{X}_n$ there is exactly one $y\in Y$ such that

$$\begin{aligned}y=\Psi_{\overline{T}\times\overline{E}_n}\!\left(\overline{\xi},\overline{x}_{[1,n]}\right)&=\Psi_{\dot{T}\times\dot{E}_n}\!\left(\left(\overline{K}_\omega^{\,-1}\times\overline{K}_{c_{\{n\}}}^{\,-1}\right)\!\left(\overline{\xi},\overline{x}_{[1,n]}\right)\right)\\ &=\Psi_{\dot{T}\times\dot{E}_n}\!\left(\overline{K}_\omega^{\,-1}\!\left(\overline{\xi}\right),\overline{K}_{c_{\{n\}}}^{\,-1}\!\left(\overline{x}_{[1,n]}\right)\right)\end{aligned} \tag{2.39}\bullet$$

**Definition 2.11.** For each $\langle\langle x_0\rangle,\overline{x}_{[1,n]}\rangle\in\overline{T}\times\overline{X}_n$, and hence for $\langle x_0,\overline{x}_{[1,n]}\rangle\in R\times\overline{X}_n$, i.e. for each $x_0\in R$ and each $\overline{x}_{[1,n]}\in\overline{X}_n$,

$$\psi^{\langle 1,n\rangle}\!\left(x_0,\overline{x}_{[1,n]}\right)\rightleftharpoons\Psi_{R\times\overline{E}_n}\!\left(x_0,\overline{x}_{[1,n]}\right)\rightleftharpoons\Psi_{\overline{T}\times\overline{E}_n}\!\left(\langle x_0\rangle,\overline{x}_{[1,n]}\right), \tag{2.40}$$

so that

$$D_{\mathrm{dp}}\!\left(\psi^{\langle 1,n\rangle}\right)=R\times\overline{E}_n,\ D_{\mathrm{a}}\!\left(\psi^{\langle 1,n\rangle}\right)=R,$$
$$D_{\mathrm{df}}\!\left(\psi^{\langle 1,n\rangle}\right)=R\times\overline{X}_n\subseteq R\times\overline{E}_n,\ D_{\mathrm{v}}\!\left(\psi^{\langle 1,n\rangle}\right)=Y\subseteq R; \tag{2.41}$$

the superscript '$\langle 1,n\rangle$' on '$\psi$' stands for *the weight of the function* $\psi^{\langle 1,n\rangle}$, i.e. for the number of independent real-valued variables as '$x_0$', '$x_1$', …, '$x_n$', to which the functional variable (operator) '$\psi^{\langle 1,n\rangle}$' can apply with the understanding that the last $n$ variables are components of an arithmetical vector. Consequently, for each $\langle x_0, \bar{x}_{[1,n]} \rangle \in R \times \bar{X}_n$ there is exactly one $y \in Y$ such that

$$y = \psi^{\langle 1,n \rangle}(x_0, \bar{x}_{[1,n]}). \qquad (2.42)\bullet$$

## 2.3. An *n*-dimensional spheres (balls), spherical surfaces, and spherical layers

**Definition 2.12.** Given $n \in \omega_1$, given $\bar{x}_{*[1,n]} \in \bar{E}_n$, for each $\bar{x}_{[1,n]} \in \bar{E}_n$,

$$\bar{x}'_{[1,n]} \triangleq \bar{x}_{[1,n]} - \bar{x}_{*[1,n]} = \langle x_1 - x_{*1}, ..., x_n - x_{*n} \rangle \triangleq \langle x'_1, ..., x'_n \rangle \in \bar{E}_n, \qquad (2.43)$$

$$x'_{[1,n]} \triangleq |\bar{x}'_{[1,n]}| = \sqrt{\sum_{i=1}^{n} x_i'^2} = \sqrt{\sum_{i=1}^{n} (x_i - x_{*i})^2} \geq 0. \qquad (2.44)$$

At *n*=1, these definitions become

$$\bar{x}'_{[1,1]} \triangleq \bar{x}_{[1,1]} - \bar{x}_{*[1,1]} = \langle x_1 - x_{*1} \rangle = \{x_1 - x_{*1}\} \in \bar{E}_1, \qquad (2.45)$$

$$x'_{[1,1]} \triangleq |\bar{x}'_{[1,1]}| = |x'_1| = |x_1 - x_{*1}| = \sqrt{(x_1 - x_{*1})^2} \geq 0. \qquad (2.46)$$

**Definition 2.13.** Given $n \in \omega_1$, given $\bar{x}_{*[1,n]} \in \bar{E}_n$, given $r \in (0, +\infty)$,

$$\bar{B}_n^o(r, \bar{x}_{*[1,n]}) \triangleq \{\bar{x}_{[1,n]} | \bar{x}_{[1,n]} \in \bar{E}_n \text{ and } |\bar{x}_{[1,n]} - \bar{x}_{*[1,n]}| < r\}, \qquad (2.47)$$

$$\bar{B}_n^c(r, \bar{x}_{*[1,n]}) \triangleq \{\bar{x}_{[1,n]} | \bar{x}_{[1,n]} \in \bar{E}_n \text{ and } |\bar{x}_{[1,n]} - \bar{x}_{*[1,n]}| \leq r\}, \qquad (2.48)$$

$$\bar{B}_n^b(r, \bar{x}_{*[1,n]}) \triangleq \{\bar{x}_{[1,n]} | \bar{x}_{[1,n]} \in \bar{E}_n \text{ and } |\bar{x}_{[1,n]} - \bar{x}_{*[1,n]}| = r\}. \qquad (2.49)$$

The sets $\bar{B}_n^o(r, \bar{x}_{*[1,n]})$, $\bar{B}_n^c(r, \bar{x}_{*[1,n]})$, and $\bar{B}_n^b(r, \bar{x}_{*[1,n]})$ are called respectively the *open sphere* (or the *spherical neighborhood*), the *closed sphere*, and *the spherical surface, in* $\bar{E}_n$ of radius *r* centered at $\bar{x}_{*[1,n]}$. The word "*ball*" can be used interchangeably with the word "*sphere*".•

**Definition 2.14.** Given $n \in \omega_1$, given $\bar{x}_{*[1,n]} \in \bar{E}_n$, given $r_1 \in (0, +\infty)$, given $r_2 \in (r_1, +\infty)$,

$$\begin{aligned}
\bar{C}_n^{oo}(r_1, r_2, \bar{x}_{*[1,n]}) &\triangleq \{\bar{x}_{[1,n]} | \bar{x}_{[1,n]} \in \bar{E}_n \text{ and } r_1 < |\bar{x}_{[1,n]} - \bar{x}_{*[1,n]}| < r_2\} \\
\bar{C}_n^{co}(r_1, r_2, \bar{x}_{*[1,n]}) &\triangleq \{\bar{x}_{[1,n]} | \bar{x}_{[1,n]} \in \bar{E}_n \text{ and } r_1 \leq |\bar{x}_{[1,n]} - \bar{x}_{*[1,n]}| < r_2\} \\
\bar{C}_n^{oc}(r_1, r_2, \bar{x}_{*[1,n]}) &\triangleq \{\bar{x}_{[1,n]} | \bar{x}_{[1,n]} \in \bar{E}_n \text{ and } r_1 < |\bar{x}_{[1,n]} - \bar{x}_{*[1,n]}| \leq r_2\} \\
\bar{C}_n^{cc}(r_1, r_2, \bar{x}_{*[1,n]}) &\triangleq \{\bar{x}_{[1,n]} | \bar{x}_{[1,n]} \in \bar{E}_n \text{ and } r_1 \leq |\bar{x}_{[1,n]} - \bar{x}_{*[1,n]}| \leq r_2\}
\end{aligned} \qquad (2.50)$$

subject to (2.44), and

$$\bar{C}_n^b(r_1, r_2, \bar{x}_{*[1,n]}) \triangleq \bar{B}_n^b(r_1, \bar{x}_{*[1,n]}) \cup \bar{B}_n^b(r_2, \bar{x}_{*[1,n]}). \qquad (2.51)$$

The set $\overline{C}_n^{oo}(r_1, r_2, \overline{x}_{*[1,n]})$, $\overline{C}_n^{co}(r_1, r_2, \overline{x}_{*[1,n]})$, $\overline{C}_n^{oc}(r_1, r_2, \overline{x}_{*[1,n]})$, or $\overline{C}_n^{cc}(r_1, r_2, \overline{x}_{*[1,n]})$ is called respectively the *open* (*open-open*), *closed-open*, *open-closed*, or *closed* (*closed-closed*) *spherical layer enclosed between the concentric spherical surfaces* $\overline{B}_n^b(r_1, \overline{x}_{*[1,n]})$ and $\overline{B}_n^b(r_2, \overline{x}_{*[1,n]})$, whereas the set $\overline{C}_n^b(r_1, r_2, \overline{x}_{*[1,n]})$ is called the total boundary surface of each one of the above four layers.•

**Corollary 2.1.**

$$\overline{B}_n^o(r, \overline{x}_{*[1,n]}) \cup \overline{B}_n^b(r, \overline{x}_{*[1,n]}) = \overline{B}_n^c(r, \overline{x}_{*[1,n]}), \tag{2.52}$$

$$\overline{B}_n^o(r, \overline{x}_{*[1,n]}) \cap \overline{B}_n^b(r, \overline{x}_{*[1,n]}) = \varnothing, \tag{2.53}$$

$$\begin{aligned}
\overline{C}_n^{oo}(r_1, r_2, \overline{x}_{*[1,n]}) &= \overline{B}_n^o(r_2, \overline{x}_{*[1,n]}) - \overline{B}_n^c(r_1, \overline{x}_{*[1,n]}), \\
\overline{C}_n^{co}(r_1, r_2, \overline{x}_{*[1,n]}) &= \overline{B}_n^o(r_2, \overline{x}_{*[1,n]}) - \overline{B}_n^o(r_1, \overline{x}_{*[1,n]}), \\
\overline{C}_n^{oc}(r_1, r_2, \overline{x}_{*[1,n]}) &= \overline{B}_n^c(r_2, \overline{x}_{*[1,n]}) - \overline{B}_n^c(r_1, \overline{x}_{*[1,n]}), \\
\overline{C}_n^{cc}(r_1, r_2, \overline{x}_{*[1,n]}) &= \overline{B}_n^c(r_2, \overline{x}_{*[1,n]}) - \overline{B}_n^o(r_1, \overline{x}_{*[1,n]}).
\end{aligned} \tag{2.54}$$

**Proof:** The corollary follows from Definitions 2.13 and 2.14.•

**Corollary 2.2.**

$$\overline{B}_1^o(r, \overline{x}_{*[1,1]}) = (x_{*1} - r, x_{*1} + r), \ \overline{B}_1^c(r, \overline{x}_{*[1,1]}) = [x_{*1} - r, x_{*1} + r], \\ \overline{B}_1^b(r, \overline{x}_{*[1,1]}) = \{x_{*1} - r, x_{*1} + r\}, \tag{2.55}$$

$$\begin{aligned}
\overline{C}_1^{oo}(r_1, r_2, \overline{x}_{*[1,1]}) &= (x_{*1} - r_2, x_{*1} - r_1) \cup (x_{*1} + r_1, x_{*1} + r_2), \\
\overline{C}_1^{co}(r_1, r_2, \overline{x}_{*[1,1]}) &= (x_{*1} - r_2, x_{*1} - r_1] \cup [x_{*1} + r_1, x_{*1} + r_2), \\
\overline{C}_1^{oc}(r_1, r_2, \overline{x}_{*[1,1]}) &= [x_{*1} - r_2, x_{*1} - r_1) \cup (x_{*1} + r_1, x_{*1} + r_2], \\
\overline{C}_1^{cc}(r_1, r_2, \overline{x}_{*[1,1]}) &= [x_{*1} - r_2, x_{*1} - r_1] \cup [x_{*1} + r_1, x_{*1} + r_2],
\end{aligned} \tag{2.56}$$

the understanding being that

$$\overline{x}_{*[1,1]} \equiv \langle x_{*1} \rangle \equiv \{x_{*1}\}. \tag{2.57}$$

**Proof:** The corollary follows from Definitions 2.13 and 2.14 at $n = 1$.•

**Convention 2.1** (*Abbreviative*). Given $n \in \omega_1$, given $\overline{x}_{*[1,n]} \in \overline{E}_n$, given $r_1 \in (0, +\infty)$, given $r_2 \in (r_1, +\infty)$,

$$\overline{X}_n^o = \overline{C}_n^{oo}(r_1, r_2, \overline{x}_{*[1,n]}), \ \overline{X}_n^c = \overline{C}_n^{cc}(r_1, r_2, \overline{x}_{*[1,n]}), \ \overline{X}_n^b \equiv \overline{C}_n^{bb}(r_1, r_2, \overline{x}_{*[1,n]}), \tag{2.58}$$

unless stated otherwise.•

# 3. Functional forms depending on a biased time-like variable and their differential and integro-differential properties

## 3.1. A basic real-valued functional form depending on one time-like variable and on n space-like real variables specified

**Hypothesis 3.1.** 1) Given $n \in \omega_1$, '$\psi^{\langle 1,n \rangle}(x_0, \bar{x}_{[1,n]})$' subject to (2.40) *is* a real-valued functional form, which along with its all second-order partial derivatives with respect to '$x_0$', '$x_1$', ..., '$x_n$' is defined and continuous for each $x_0 \in R$ and each $\bar{x}_{[1,n]} \in \bar{B}_n^c(r_2, \bar{x}_{*[1,n]})$ subject to (2.48), so that symbolically

$$\psi^{\langle 1,n \rangle} \in \mathcal{C}_2\left(R \times \bar{B}_n^c(r_2, \bar{x}_{*[1,n]})\right), \tag{3.1}$$

where $\psi^{\langle 1,n \rangle}$ is the *full associated function of the functional form* '$\psi^{\langle 1,n \rangle}(x_0, \bar{x}_{[1,n]})$'.•

**Comment 3.1.** Here, and generally throughout the treatise, a *quotation* that is formed by enclosing a logograph in *light-faced single quotation marks*, ' ' is a constant that denotes by default a certain *class of homolographic* (*photographic, congruent or proportional*) *isotokens* (*graphic tokens*) of its interior and that is therefore called a *homolographic*, or *photographic*, *autonymous quotation* (*HAQ*). By contrast, a *quotation* that is formed by enclosing a *logographic placrholder* in *bold-faced single quotation marks*, ' ' is a placeholder itself, whose range is a certain class of HAQ's and which is therefore called a *quasi-homolographic*, or *quasi-photographic*, *autonymous quotation* (*QHAQ*). It is understood that, once the interior of a QHAQ is replaced by a concrete logograph of the range of its interior, the bold-faced single quotation marks of the QHAQ should be replaced with light-faced ones, so that the QHAQ is replaced with the respective HAQ of its range. Thus, in accordance with the pertinent phraseology, in which single quotation marks are used but not mentioned, '$\psi^{\langle 1,n \rangle}(x_0, \bar{x}_{[1,n]})$' *is* a placeholder for any functional form of its range, whereas '$\psi^{\langle 1,n \rangle}(x_0, \bar{x}_{[1,n]})$' or equivocally $\psi^{\langle 1,n \rangle}(x_0, \bar{x}_{[1,n]})$ *is* a functional form of that range.•

**Comment 3.2.** For the reasons that will become evident in section 5, the 1-vector-valued variable '$\bar{\xi}$' and also its only real-valued component '$\xi$' or '$x_0$' are said to be *time-like*, while the $n$-vector-valued variable '$\bar{x}_{[1,n]}$' and also any its component '$x_i$' with $i \in \omega_{1,n}$ are said is said to be *spatial*. At the same time, in the general case, $\bar{T} \times \bar{X}_n^c$ is not assumed to be a part of an $(n+1)$-dimensional pseudo-Euclidean space of index 1, i.e. no pseudo-Euclidean metrics is defined for ordered pairs such as $\langle \bar{\xi}, \bar{x}_{[1,n]} \rangle$ being elements $\bar{T} \times \bar{X}_n^c$. It is only in the case of $n=3$ that $\bar{T} \times \bar{X}_3^c$ can after be postulated to be a part of an $(n+1)$-dimensional pseudo-Euclidean space of index 1. However, such an assumption is not necessary for deducing the desired integro-differential theorems.•

### 3.2. A biased functional form and its basic differential properties

**Definition 3.1.** Given $n \in \omega_1$, given $x_{*0} \in (-\infty, \infty)$, given $\bar{x}_{*[1,n]} \in \bar{E}_n$, given $\bar{x}_{[1,n]} \in \bar{E}_n$, given $\lambda \in \{-1, 0, 1\}$,

$$x_0 \equiv x_{*0} + \lambda x'_{[1,n]} \equiv x_{*0} + \lambda |\bar{x}_{[1,n]} - \bar{x}_{*[1,n]}| \in (-\infty, \infty), \quad (3.2)$$

subject to Definition 2.12. Consequently, the functional form '$x_{*0} + \lambda |\bar{x}_{[1,n]} - \bar{x}_{*[1,n]}|$' and a *composite functional form* '$\psi^{\langle 1,n \rangle}(x_{*0} + \lambda |\bar{x}_{[1,n]} - \bar{x}_{*[1,n]}|, \bar{x}_{[1,n]})$' subject to Definition 2.11, i.e. any functional form in the range of the placeholder '$\psi^{\langle 1,n \rangle}(x_{*0} + \lambda |\bar{x}_{[1,n]} - \bar{x}_{*[1,n]}|, \bar{x}_{[1,n]})$', are said to be *retarded* or *negatively biased* if $\lambda = -1$, *advanced* or *positively biased* if $\lambda = 1$, and *unbiased* if $\lambda = 0$.

**Comment 3.3.** By definition (2.49), it follows from (3.3) that

$$x_{*0} + \lambda x'_{[1,n]} = x_{*0} + \lambda r \text{ if } \bar{x}_{[1,n]} \in \bar{B}_n^b(r, \bar{x}_{*[1,n]}). \quad (3.3)\bullet$$

**Comment 3.4.** The terminology that is introduced in Definition 3.1 applies also to the composite functional form '$\psi^{\lambda, \langle 1,n,n \rangle}(x_{*0}, \bar{x}_{[1,n]}, \bar{x}_{*[1,n]})$', which is defined as:

$$\begin{aligned}\psi^{\lambda, \langle 1,n,n \rangle}(x_{*0}, \bar{x}_{[1,n]}, \bar{x}_{*[1,n]}) &\equiv \psi^{\langle 1,n \rangle}(x_{*0} + \lambda x'_{[1,n]}, \bar{x}_{[1,n]}) \\ &\equiv \psi^{\langle 1,n \rangle}(x_{*0} + \lambda |\bar{x}_{[1,n]} - \bar{x}_{*[1,n]}|, \bar{x}_{[1,n]}),\end{aligned} \quad (3.4)$$

In this case, if $\lambda \in \{-1, 1\}$ then the additional superscript '$\langle 1,n,n \rangle$' on '$\psi$' indicates that the functional variable (operator) '$\psi^{\lambda, \langle 1,n,n \rangle}$' applies simultaneously to $2n+1$ independent real-valued variables, one of which, as '$x_{*0}$', is a time or time-like scalar, while two groups of $n$ variables in each group, as '$x_1$', …, '$x_n$' and '$x_{*1}$', …, '$x_{*n}$' form two $n$-dimensional arithmetical vectors. On the other hand, if $\lambda = 0$ then

$$\psi^{0, \langle 1,n,n \rangle}(x_{*0}, \bar{x}_{[1,n]}, \bar{x}_{*[1,n]}) = \psi^{\langle 1,n \rangle}(x_{*0}, \bar{x}_{[1,n]}). \quad (3.4a)\bullet$$

**Lemma 3.1.** Given $n \in \omega_1$, given $x_{*0} \in R$, given $\bar{x}_{*[1,n]} \in \bar{E}_n$, given $\lambda \in \{-1, 1\}$, for each $\bar{x}_{[1,n]} \in \bar{E}_n - \{\bar{x}_{*[1,n]}\}$, for each $i \in \omega_{1,n}$, for each $j \in \omega_{1,n}$:

$$\frac{\partial (x_{*0} + \lambda x'_{[1,n]})}{\partial x_i} = \lambda \frac{\partial x'_{[1,n]}}{\partial x_i}, \tag{3.5}$$

$$\frac{\partial x'_{[1,n]}}{\partial x_i} = \frac{x'_i}{x'_{[1,n]}}, \quad \frac{\partial^2 x'_{[1,n]}}{\partial x_i \partial x_j} = \frac{{x'_{[1,n]}}^2 \delta_{ij} - x'_i x'_j}{{x'_{[1,n]}}^3}. \tag{3.6}$$

**Proof:** The lemma follows from definition (3.2) subject to (2.43) and (2.44).•

**Definition 3.2.** When *appropriate and convenient*, we shall use the following notations for the associated functions of partial derivatives of the functional forms $\psi^{\langle 1,n \rangle}(x_{*0}, \bar{x}_{[1,n]})$:

$$\psi^{\langle 1,n \rangle}_{x_0} \equiv \frac{\partial \psi^{\langle 1,n \rangle}}{\partial x_0}, \psi^{\langle 1,n \rangle}_{x_0 x_0} \equiv \frac{\partial^2 \psi^{\langle 1,n \rangle}}{\partial x_0^2}, \tag{3.7}$$

$$\psi^{\langle 1,n \rangle}_{x_i} \equiv \frac{\partial \psi^{\langle 1,n \rangle}}{\partial x_i}, \psi^{\langle 1,n \rangle}_{x_j} \equiv \frac{\partial \psi^{\langle 1,n \rangle}}{\partial x_j}, \psi^{\langle 1,n \rangle}_{x_i x_0} \equiv \frac{\partial^2 \psi^{\langle 1,n \rangle}}{\partial x_i \partial x_0}, \psi^{\langle 1,n \rangle}_{x_0 x_i} \equiv \frac{\partial^2 \psi^{\langle 1,n \rangle}}{\partial x_0 \partial x_i},$$

$$\psi^{\langle 1,n \rangle}_{x_i x_i} \equiv \frac{\partial^2 \psi^{\langle 1,n \rangle}}{\partial x_i^2}, \psi^{\langle 1,n \rangle}_{x_i x_j} \equiv \frac{\partial^2 \psi^{\langle 1,n \rangle}}{\partial x_i \partial x_j}, \psi^{\langle 1,n \rangle}_{x_j x_j} \equiv \frac{\partial^2 \psi^{\langle 1,n \rangle}}{\partial x_j^2} \tag{3.8}$$

for each $i \in \omega_{1,n}$ and each $j \in \omega_{1,n}$,

the understanding being that, e.g.,

$$\psi^{\langle 1,n \rangle}_{x_0}(x_{*0} + \lambda x'_{[1,n]}, \bar{x}_{[1,n]}) \equiv \left[ \frac{\partial \psi^{\langle 1,n \rangle}(x_0, \bar{x}_{[1,n]})}{\partial x_0} \right]_{x_0 = x_{*0} + \lambda x'_{[1,n]}},$$

$$\psi^{\langle 1,n \rangle}_{x_0 x_0}(x_{*0} + \lambda x'_{[1,n]}, \bar{x}_{[1,n]}) \equiv \left[ \frac{\partial^2 \psi^{\langle 1,n \rangle}(x_0, \bar{x}_{[1,n]})}{\partial x_0^2} \right]_{x_0 = x_{*0} + \lambda x'_{[1,n]}},$$

$$\psi^{\langle 1,n \rangle}_{x_i}(x_{*0} + \lambda x'_{[1,n]}, \bar{x}_{[1,n]}) \equiv \left[ \frac{\partial \psi^{\langle 1,n \rangle}(x_0, \bar{x}_{[1,n]})}{\partial x_i} \right]_{x_0 = x_{*0} + \lambda x'_{[1,n]}}, \tag{3.9}$$

$$\psi^{\langle 1,n \rangle}_{x_0 x_i}(x_{*0} + \lambda x'_{[1,n]}, \bar{x}_{[1,n]}) \equiv \left[ \frac{\partial^2 \psi^{\langle 1,n \rangle}(x_0, \bar{x}_{[1,n]})}{\partial x_0 \partial x_i} \right]_{x_0 = x_{*0} + \lambda x'_{[1,n]}},$$

$$\psi^{\langle 1,n \rangle}_{x_i x_j}(x_{*0} + \lambda x'_{[1,n]}, \bar{x}_{[1,n]}) \equiv \left[ \frac{\partial^2 \psi^{\langle 1,n \rangle}(x_0, \bar{x}_{[1,n]})}{\partial x_i \partial x_j} \right]_{x_0 = x_{*0} + \lambda x'_{[1,n]}},$$

etc.•

**Convention 3.1.** When there is no danger of confusion, '$\psi^{\langle 1,n \rangle}$' can be abbreviated as '$\psi$' and '$\psi^{\lambda,\langle 1,n,n \rangle}$' as '$\psi^{(\lambda)}$'.•

**Convention 3.2.** Every subsequent equality or train of equalities that involves some of the variables '$n$', '$\lambda$', '$\bar{x}_{*0}$', '$\bar{x}_{*[1,n]}$', and '$\bar{x}_{[1,n]}$' is supposed to be preceded by the respective ones of the qualifiers: 'given $n \in \omega_1$', 'given $\lambda \in \{-1,1\}$', 'given $\bar{x}_{*0} \in R$', 'given $\bar{x}_{*[1,n]} \in \bar{E}_n$', and by the quantifier 'for each $\bar{x}_{[1,n]} \in \bar{X}_n^c$', unless stated otherwise.•

**Lemma 3.2:**

$$\frac{\partial^2 \psi^{\langle 1,n \rangle}\left(x_{*0} + \lambda x'_{[1,n]}, \bar{x}_{[1,n]}\right)}{\partial x_i \partial x_j}$$

$$= \psi^{\langle 1,n \rangle}_{x_i x_j}\left(x_{*0} + \lambda x'_{[1,n]}, \bar{x}_{[1,n]}\right) - \frac{x'_i x'_j}{{x'_{[1,n]}}^2} \psi^{\langle 1,n \rangle}_{x_0 x_0}\left(x_{*0} + \lambda x'_{[1,n]}, \bar{x}_{[1,n]}\right)$$

$$+ \lambda \frac{\partial}{\partial x_i} \left[ \frac{x'_j}{x'_{[1,n]}} \psi^{\langle 1,n \rangle}_{x_0}\left(x_{*0} + \lambda x'_{[1,n]}, \bar{x}_{[1,n]}\right) \right] \quad (3.10)$$

$$+ \lambda \frac{\partial}{\partial x_j} \left[ \frac{x'_i}{x'_{[1,n]}} \psi^{\langle 1,n \rangle}_{x_0}\left(x_{*0} + \lambda x'_{[1,n]}, \bar{x}_{[1,n]}\right) \right]$$

$$- \lambda \frac{{x'_{[1,n]}}^2 \delta_{ij} - x'_i x'_j}{{x'_{[1,n]}}^3} \psi^{\langle 1,n \rangle}_{x_0}\left(x_{*0} + \lambda x'_{[1,n]}, \bar{x}_{[1,n]}\right).$$

**Proof:** By the chain rule of differentiation of a composite functional form, it follows that

$$\frac{\partial \psi^{\langle 1,n \rangle}\left(x_{*0} + \lambda x'_{[1,n]}, \bar{x}_{[1,n]}\right)}{\partial x_j}$$

$$= \lambda \frac{\partial x'_{[1,n]}}{\partial x_j} \psi^{\langle 1,n \rangle}_{x_0}\left(x_{*0} + \lambda x'_{[1,n]}, \bar{x}_{[1,n]}\right) + \psi^{\langle 1,n \rangle}_{x_j}\left(x_{*0} + \lambda x'_{[1,n]}, \bar{x}_{[1,n]}\right). \quad (3.10_1)$$

Differentiation of both sides of the identity (3.10$_1$) with respect to $x_i$ yields:

$$\frac{\partial^2 \psi^{\langle 1,n \rangle}\left(x_{*0} + \lambda x'_{[1,n]}, \bar{x}_{[1,n]}\right)}{\partial x_i \partial x_j}$$

$$= \lambda \frac{\partial}{\partial x_i}\left[ \frac{\partial x'_{[1,n]}}{\partial x_j} \psi^{\langle 1,n \rangle}_{x_0}\left(x_{*0} + \lambda x'_{[1,n]}, \bar{x}_{[1,n]}\right) \right] + \frac{\partial}{\partial x_i} \psi^{\langle 1,n \rangle}_{x_j}\left(x_{*0} + \lambda x'_{[1,n]}, \bar{x}_{[1,n]}\right). \quad (3.10_2)$$

By the same chain rule, the last term in (3.10$_2$) can be developed thus:

$$\frac{\partial}{\partial x_i} \psi^{\langle 1,n \rangle}_{x_j}\left(x_{*0} + \lambda x'_{[1,n]}, \bar{x}_{[1,n]}\right) = \lambda \frac{\partial x'_{[1,n]}}{\partial x_i} \psi^{\langle 1,n \rangle}_{x_0 x_j}\left(x_{*0} + \lambda x'_{[1,n]}, \bar{x}_{[1,n]}\right)$$

$$+ \psi^{\langle 1,n \rangle}_{x_i x_j}\left(x_{*0} + \lambda x'_{[1,n]}, \bar{x}_{[1,n]}\right). \quad (3.10_3)$$

Analogously,

$$\frac{\partial}{\partial x_j} \psi^{\langle 1,n \rangle}_{x_0}\left(x_{*0} + \lambda x'_{[1,n]}, \bar{x}_{[1,n]}\right) = \lambda \frac{\partial x'_{[1,n]}}{\partial x_j} \psi^{\langle 1,n \rangle}_{x_0 x_0}\left(x_{*0} + \lambda x'_{[1,n]}, \bar{x}_{[1,n]}\right)$$

$$+ \psi^{\langle 1,n \rangle}_{x_j x_0}\left(x_{*0} + \lambda x'_{[1,n]}, \bar{x}_{[1,n]}\right), \quad (3.10_4)$$

whence

$$\psi_{x_j x_0}^{\langle 1,n \rangle}(x_{*0} + \lambda x'_{[1,n]}, \bar{x}_{[1,n]})$$
$$= \frac{\partial}{\partial x_j} \psi_{x_0}^{\langle 1,n \rangle}(x_{*0} + \lambda x'_{[1,n]}, \bar{x}_{[1,n]}) - \lambda \frac{\partial x'_{[1,n]}}{\partial x_j} \psi_{x_0 x_0}^{\langle 1,n \rangle}(x_{*0} + \lambda x'_{[1,n]}, \bar{x}_{[1,n]}).$$
(3.10$_5$)]

Since

$$\psi_{x_0 x_j}^{\langle 1,n \rangle}(x_{*0} + \lambda x'_{[1,n]}, \bar{x}_{[1,n]}) = \psi_{x_j x_0}^{\langle 1,n \rangle}(x_{*0} + \lambda x'_{[1,n]}, \bar{x}_{[1,n]}),$$
(3.10$_6$)

therefore equation (3.10$_3$) can be developed by (3.10$_5$) thus:

$$\frac{\partial}{\partial x_i} \psi_{x_j}^{\langle 1,n \rangle}(x_{*0} + \lambda x'_{[1,n]}, \bar{x}_{[1,n]}) = \lambda \frac{\partial x'_{[1,n]}}{\partial x_i} \frac{\partial}{\partial x_j} \psi_{x_0}^{\langle 1,n \rangle}(x_{*0} + \lambda x'_{[1,n]}, \bar{x}_{[1,n]})$$
$$- \frac{\partial x'_{[1,n]}}{\partial x_i} \frac{\partial x'_{[1,n]}}{\partial x_j} \psi_{x_0 x_0}^{\langle 1,n \rangle}(x_{*0} + \lambda x'_{[1,n]}, \bar{x}_{[1,n]}) + \psi_{x_i x_j}^{\langle 1,n \rangle}(x_{*0} + \lambda x'_{[1,n]}, \bar{x}_{[1,n]}),$$
(3.10$_7$)

for $\lambda^2 = 1$. This can immediately be rewritten as

$$\frac{\partial}{\partial x_i} \psi_{x_j}^{\langle 1,n \rangle}(x_{*0} + \lambda x'_{[1,n]}, \bar{x}_{[1,n]}) = \lambda \frac{\partial}{\partial x_j}\left[\frac{\partial x'_{[1,n]}}{\partial x_i} \psi_{x_0}^{\langle 1,n \rangle}(x_{*0} + \lambda x'_{[1,n]}, \bar{x}_{[1,n]})\right]$$
$$- \lambda \frac{\partial^2 x'_{[1,n]}}{\partial x_i \partial x_j} \psi_{x_0}^{\langle 1,n \rangle}(x_{*0} + \lambda x'_{[1,n]}, \bar{x}_{[1,n]})$$
$$+ \psi_{x_i x_j}^{\langle 1,n \rangle}(x_{*0} + \lambda x'_{[1,n]}, \bar{x}_{[1,n]}) - \frac{\partial x'_{[1,n]}}{\partial x_i} \frac{\partial x'_{[1,n]}}{\partial x_j} \psi_{x_0 x_0}^{\langle 1,n \rangle}(x_{*0} + \lambda x'_{[1,n]}, \bar{x}_{[1,n]}).$$
(3.10$_8$)

Substitution of (3.10$_8$) into (3.10$_2$) yields

$$\frac{\partial^2 \psi^{\langle 1,n \rangle}(x_{*0} + \lambda x'_{[1,n]}, \bar{x}_{[1,n]})}{\partial x_i \partial x_j}$$
$$= \psi_{x_i x_j}^{\langle 1,n \rangle}(x_{*0} + \lambda x'_{[1,n]}, \bar{x}_{[1,n]}) - \frac{\partial x'_{[1,n]}}{\partial x_i} \frac{\partial x'_{[1,n]}}{\partial x_j} \psi_{x_0 x_0}^{\langle 1,n \rangle}(x_{*0} + \lambda x'_{[1,n]}, \bar{x}_{[1,n]})$$
$$+ \lambda \frac{\partial}{\partial x_i}\left[\frac{\partial x'_{[1,n]}}{\partial x_j} \psi_{x_0}^{\langle 1,n \rangle}(x_{*0} + \lambda x'_{[1,n]}, \bar{x}_{[1,n]})\right] + \lambda \frac{\partial}{\partial x_j}\left[\frac{\partial x'_{[1,n]}}{\partial x_i} \psi_{x_0}^{\langle 1,n \rangle}(x_{*0} + \lambda x'_{[1,n]}, \bar{x}_{[1,n]})\right]$$
$$- \lambda \frac{\partial^2 x'_{[1,n]}}{\partial x_i \partial x_j} \psi_{x_0}^{\langle 1,n \rangle}(x_{*0} + \lambda x'_{[1,n]}, \bar{x}_{[1,n]}).$$
(3.10$_9$)

By (3.6), equation (3.10$_9$) reduces to (3.10). QED.•

**Comment 3.5.** Differentiation of both sides of the identity (3.10$_1$) with respect to $x_0$ yields:

$$\frac{\partial^2 \psi^{\langle 1,n \rangle}\left(x_{*0} + \lambda x'_{[1,n]}, \bar{x}_{[1,n]}\right)}{\partial x_0 \partial x_j}$$

$$= \lambda \frac{\partial x'_{[1,n]}}{\partial x_j} \psi^{\langle 1,n \rangle}_{x_0 x_0}\left(x_{*0} + \lambda x'_{[1,n]}, \bar{x}_{[1,n]}\right) + \psi^{\langle 1,n \rangle}_{x_0 x_j}\left(x_{*0} + \lambda x'_{[1,n]}, \bar{x}_{[1,n]}\right).$$

(3.10$_{10}$)

By (3.10$_6$), it follows from comparison of (3.10$_4$) and (3.10$_{10}$) that:

$$\frac{\partial^2 \psi^{\langle 1,n \rangle}\left(x_{*0} + \lambda x'_{[1,n]}, \bar{x}_{[1,n]}\right)}{\partial x_j \partial x_0} = \frac{\partial^2 \psi^{\langle 1,n \rangle}\left(x_{*0} + \lambda x'_{[1,n]}, \bar{x}_{[1,n]}\right)}{\partial x_0 \partial x_j}$$

(3.10$_{11}$)•

**Definition 3.3.**

$$\Delta_{[1,n]} \equiv \sum_{i=1}^{n} \nabla_i^2 = \sum_{i=1}^{n} \frac{\partial^2}{\partial x_i^2},$$

(3.11)

$$D_{[0,n]} \equiv \Delta_{[1,n]} - \frac{\partial^2}{\partial x_0^2} = \sum_{i=1}^{n} \frac{\partial^2}{\partial x_i^2} - \frac{\partial^2}{\partial x_0^2}.$$

(3.12)

$\Delta_{[1,n]}$ is called the *Laplacian operator* in $\bar{E}_n$, and $D_{[0,n]}$ is called the *d'Alambertian operator* in $\bar{T} \times \bar{E}_n$.•

**Lemma 3.3.**

$$\left[D_{[0,n]} \psi^{\langle 1,n \rangle}\left(x_0, \bar{x}_{[1,n]}\right)\right]_{x_0 = x_{*0} + \lambda x'_{[1,n]}}$$

$$= \psi^{\langle 1,n \rangle}_{x_i x_i}\left(x_{*0} + \lambda x'_{[1,n]}, \bar{x}_{[1,n]}\right) - \psi^{\langle 1,n \rangle}_{x_0 x_0}\left(x_{*0} + \lambda x'_{[1,n]}, \bar{x}_{[1,n]}\right)$$

$$= \sum_{i=1}^{n} \frac{\partial}{\partial x_i}\left[\frac{\partial \psi^{\langle 1,n \rangle}\left(x_{*0} + \lambda x'_{[1,n]}, \bar{x}_{[1,n]}\right)}{\partial x_i} - \frac{2\lambda x'_i}{x'_{[1,n]}} \psi^{\langle 1,n \rangle}_{x_0}\left(x_{*0} + \lambda x'_{[1,n]}, \bar{x}_{[1,n]}\right)\right]$$

$$+ \frac{\lambda(n-1)}{x'_{[1,n]}} \psi^{\langle 1,n \rangle}_{x_0}\left(x_{*0} + \lambda x'_{[1,n]}, \bar{x}_{[1,n]}\right),$$

(3.13)

subject to (2.43)–(2.46) and (3.12).

**Proof:** Summation of both sides of equation (3.10) at $j = i$ over all values of '$i$' from 1 to $n$ yields

$$\sum_{i=1}^{n} \frac{\partial^2 \psi^{\langle 1,n \rangle}\left(x_{*0} + \lambda x'_{[1,n]}, \bar{x}_{[1,n]}\right)}{\partial x_i^2}$$

$$= \sum_{i=1}^{n} \psi^{\langle 1,n \rangle}_{x_i x_i}\left(x_{*0} + \lambda x'_{[1,n]}, \bar{x}_{[1,n]}\right) - \psi^{\langle 1,n \rangle}_{x_0 x_0}\left(x_{*0} + \lambda x'_{[1,n]}, \bar{x}_{[1,n]}\right)$$

$$+ 2\lambda \sum_{i=1}^{n} \frac{\partial}{\partial x_i}\left[\frac{x'_i}{x'_{[1,n]}} \psi^{v\langle 1,n \rangle}_{x_0}\left(x_{*0} + \lambda x'_{[1,n]}, \bar{x}_{[1,n]}\right)\right]$$

$$- \frac{\lambda(n-1)}{x'_{[1,n]}} \psi^{\langle 1,n \rangle}_{x_0}\left(x_{*0} + \lambda x'_{[1,n]}, \bar{x}_{[1,n]}\right),$$

(3.13$_1$)

which is equivalent to (3.13). QED.•

**Comment 3.6.** Identity (3.13) has the following peculiarities.

1) The expression on the right-hand side of identity (3.13) is the sum of two terms, of which the first is the *divergence* of a certain expression, while the second term involves no differentiation with respect to spatial coordinates at all.

2) By (3.12), the functional form on the left-hand side of identity (3.13) is regular at each $\langle x_{*0}, \bar{x}_{[1,n]} \rangle \in R \times \bar{X}_n^c$. At the same time, the last term on the right-hand side of identity (3.13) unlimitedly increases as $\bar{x}_{[1,n]}$ tends to $\bar{x}_{*[1,n]}$. Therefore, a term of the same absolute value and of the opposite sign must be present in the first term on the right-hand side of identity (3.13).

3) Since

$$\sum_{i=1}^{n} \frac{\partial}{\partial x_i} \frac{x_i'}{x_{[1,n]}'} = \sum_{i=1}^{n} \left[ \frac{1}{x_{[1,n]}'} \frac{\partial x_i'}{\partial x_i} + \frac{\partial}{\partial x_i} \frac{1}{x_{[1,n]}'} \right] = \sum_{i=1}^{n} \left[ \frac{1}{x_{[1,n]}'} + x_i' \frac{\partial}{\partial x_i} \frac{1}{x_{[1,n]}'} \right]$$

$$= \frac{1}{x_{[1,n]}'} \sum_{i=1}^{n} 1 - \frac{1}{x_{[1,n]}'^3} \sum_{i=1}^{n} x_i' x_i' = \frac{n-1}{x_{[1,n]}'}, \tag{3.14}$$

therefore (3.13) can alternatively be written as:

$$\left[ D_{[0,n]} \psi^{\langle 1,n \rangle}(x_0, \bar{x}_{[1,n]}) \right]_{x_0 = x_{*0} + \lambda x_{[1,n]}'}$$

$$= \sum_{i=1}^{n} \frac{\partial^2 \psi^{\langle 1,n \rangle}(x_{*0} + \lambda x_{[1,n]}', \bar{x}_{[1,n]})}{\partial x_i^2} - \frac{2\lambda}{x_{[1,n]}'} \sum_{i=1}^{n} x_i' \frac{\partial \psi_{x_0}^{\langle 1,n \rangle}(x_{*0} + \lambda x_{[1,n]}', \bar{x}_{[1,n]})}{\partial x_i} \tag{3.13$_2$}$$

$$- \frac{\lambda(n-1)}{x_{[1,n]}'} \psi_{x_0}^{\langle 1,n \rangle}(x_{*0} + \lambda x_{[1,n]}', \bar{x}_{[1,n]}).$$

Unlike identity (3.13), the expression on the right-hand side of the identity (3.13$_2$) involves two items not being divergences of functions. Therefore, the latter identity will not be used in the sequel.•

### 3.3. Two basic differential tautologies for a pair of biased functional forms

**Hypothesis 3.2.** Hypothesis 7.1 applies with '$\phi$' in place of '$\psi$'. That is to say, $\phi^{\langle 1,n \rangle}$ is another real-valued function such that $\phi^{\langle 1,n \rangle} \in C_2(R \times \bar{B}_n^c(r_2, \bar{x}_{*[1,n]}))$.•

**Convention 3.3.** 1) We shall use the following abbreviations:

$$\psi^{\langle 1,n \rangle} \equiv \psi^{\langle 1,n \rangle}(x_{*0} + \lambda x_{[1,n]}', \bar{x}_{[1,n]}), \phi^{\langle 1,n \rangle} \equiv \phi^{\langle 1,n \rangle}(x_{*0} + \lambda x_{[1,n]}', \bar{x}_{[1,n]}),$$

$$\psi_{\oplus}^{\langle 1,n \rangle} \equiv \psi_{\oplus}^{\langle 1,n \rangle}(x_{*0} + \lambda x_{[1,n]}', \bar{x}_{[1,n]}), \phi_{\oplus}^{\langle 1,n \rangle} \equiv \phi_{\oplus}^{\langle 1,n \rangle}(x_{*0} + \lambda x_{[1,n]}', \bar{x}_{[1,n]}),$$

where '⊕' is an ellipsis, the occurrences of which in a definiens and in the respective definiendum can be replaced alike by any one of the subscripts:

$$\text{'}x_0\text{', '}x_0 x_0\text{', '}x_i\text{', '}x_j\text{', '}x_i x_0\text{', '}x_0 x_i\text{', '}x_i x_i\text{', '}x_i x_j\text{', '}x_j x_j\text{', '}x_0 x_i\text{',} \tag{3.16}$$

in accordance with definition (3.8).

2) '$\tau_n$' is a placeholder of the list of parameters '$\bar{x}_{*0}, \lambda, \bar{x}_{*[1,n]}$', i.e. $\tau_n$ *is* that list (without quotation marks), – the parameters, which are determined the biased functional form $x_{*0} + \lambda \left| \bar{x}_{[1,n]} - \bar{x}_{*[1,n]} \right|$. ●

**Lemma 3.4:** *The first (asymmetric) differential tautology for a pair of biased functional forms.*

$$\sum_{\alpha=0}^{3} K^{(\alpha)}\left(\phi^{\langle 1,n \rangle}, \psi^{\langle 1,n \rangle}, \tau_n, \bar{x}_{[1,n]}\right) = 0, \tag{3.17}$$

where

$$K^{(0)}\left(\phi^{\langle 1,n \rangle}, \psi^{\langle 1,n \rangle}, \tau_n, \bar{x}_{[1,n]}\right) \equiv \phi^{\langle 1,n \rangle} \left[ \sum_{i=1}^{n} \psi^{\langle 1,n \rangle}_{x_i x_i} - \psi^{\langle 1,n \rangle}_{x_0 x_0} \right], \tag{3.18}$$

$$K^{(1)}\left(\phi^{\langle 1,n \rangle}, \psi^{\langle 1,n \rangle}, \tau_n, \bar{x}_{[1,n]}\right) \equiv -\frac{\lambda}{x'_{[1,n]}} \left[ (n-1)\phi^{\langle 1,n \rangle} + 2\sum_{i=1}^{n} x'_i \frac{\partial \phi^{\langle 1,n \rangle}}{\partial x_i} \right] \psi^{\langle 1,n \rangle}_{x_0}, \tag{3.19}$$

$$K^{(2)}\left(\phi^{\langle 1,n \rangle}, \psi^{\langle 1,n \rangle}, \tau_n, \bar{x}_{[1,n]}\right) \equiv \sum_{i=1}^{n} \frac{\partial \phi^{\langle 1,n \rangle}}{\partial x_i} \frac{\partial \psi^{\langle 1,n \rangle}}{\partial x_i}, \tag{3.20}$$

$$K^{(3)}\left(\phi^{\langle 1,n \rangle}, \psi^{\langle 1,n \rangle}, \tau_n, \bar{x}_{[1,n]}\right) \equiv -\sum_{i=1}^{n} \frac{\partial}{\partial x_i} L_i\left(\phi^{\langle 1,n \rangle}, \psi^{\langle 1,n \rangle}, \tau, \bar{x}_{[1,n]}\right) \tag{3.21}$$

subject to

$$L_i\left(\phi^{\langle 1,n \rangle}, \psi^{\langle 1,n \rangle}, \tau_n, \bar{x}_{[1,n]}\right) \equiv \phi^{\langle 1,n \rangle} \left[ \psi^{\langle 1,n \rangle}_{x_i} - \frac{\lambda x'_i}{x'_{[1,n]}} \psi^{\langle 1,n \rangle}_{x_0} \right] \text{ for each } i \in \omega_{1,n}. \tag{3.22}$$

**Proof:** Multiplying both sides of the identity (3.13) by $\phi^{\langle 1,n \rangle}\left(x_{*0} + \lambda x'_{[1,n]}, \bar{x}_{[1,n]}\right)$ and then making use of (3.18) and (3.19) subject to (3.15), one can, after simple reduction, write the result in the form

$$K^{(0)}\left(\phi^{\langle 1,n\rangle},\psi^{\langle 1,n\rangle},\tau_n,\overline{x}_{[1,n]}\right)=-K^{(1)}\left(\phi^{\langle 1,n\rangle},\psi^{\langle 1,n\rangle},\tau_n,\overline{x}_{[1,n]}\right)$$

$$-\frac{2\lambda}{x'_{[1,n]}}\psi_{x_0}^{\langle 1,n\rangle}\sum_{i=1}^{n}x'_i\frac{\partial\phi^{\langle 1,n\rangle}}{\partial x_i}+\phi^{\langle 1,n\rangle}\sum_{i=1}^{n}\frac{\partial}{\partial x_i}\left[\frac{\partial\psi^{\langle 1,n\rangle}}{\partial x_i}-\frac{2\lambda x'_i}{x'_{[1,n]}}\psi_{x_0}^{\langle 1,n\rangle}\right]$$

$$=-K^{(1)}\left(\phi^{\langle 1,n\rangle},\psi^{\langle 1,n\rangle},\tau,\overline{x}_{[1,n]}\right)-\frac{2\lambda}{x'_{[1,n]}}\psi_{x_0}^{\langle 1,n\rangle}\sum_{i=1}^{n}x'_i\frac{\partial\phi^{\langle 1,n\rangle}}{\partial x_i} \quad (3.23)$$

$$+\sum_{i=1}^{n}\frac{\partial}{\partial x_i}\left[\phi^{\langle 1,n\rangle}\left[\frac{\partial\psi^{\langle 1,n\rangle}}{\partial x_i}-\frac{2\lambda x'_i}{x'_{[1,n]}}\psi_{x_0}^{\langle 1,n\rangle}\right]\right]-\sum_{i=1}^{n}\frac{\partial\phi^{\langle 1,n\rangle}}{\partial x_i}\frac{\partial\psi^{\langle 1,n\rangle}}{\partial x_i}.$$

Under the definition

$$L_i\left(\phi^{\langle 1,n\rangle},\psi^{\langle 1,n\rangle},\tau_n,\overline{x}_{[1,n]}\right)\doteq\phi^{\langle 1,n\rangle}\left[\frac{\partial\psi^{\langle 1,n\rangle}}{\partial x_i}-\frac{2\lambda x'_i}{x'_{[1,n]}}\psi_{x_0}^{\langle 1,n\rangle}\right]\text{ for each }i\in\omega_{1,n}, \quad (3.24)$$

the final result in (3.23) reduces to

$$K^{(0)}\left(\phi^{\langle 1,n\rangle},\psi^{\langle 1,n\rangle},\tau_n,\overline{x}_{[1,n]}\right)=-\sum_{\alpha=1}^{3}K^{(\alpha)}\left(\phi^{\langle 1,n\rangle},\psi^{\langle 1,n\rangle},\tau_n,\overline{x}_{[1,n]}\right) \quad (3.17_1)$$

subject to (3.18)–(3.21). However, by the chain rule of differentiation of a composite functional form, it follows that

$$\frac{\partial\psi^{\langle 1,n\rangle}\left(x_{*0}+\lambda x'_{[1,n]},\overline{x}_{[1,n]}\right)}{\partial x_i}$$
$$=\lambda\frac{\partial x'_{[1,n]}}{\partial x_i}\psi_{x_0}^{\langle 1,n\rangle}\left(x_{*0}+\lambda x'_{[1,n]},\overline{x}_{[1,n]}\right)+\psi_{x_i}^{\langle 1,n\rangle}\left(x_{*0}+\lambda x'_{[1,n]},\overline{x}_{[1,n]}\right) \quad (3.24_1)$$

(cf. $(3.10_1)$), i.e. briefly

$$\frac{\partial\psi^{\langle 1,n\rangle}}{\partial x_i}=\lambda\frac{\partial x'_{[1,n]}}{\partial x_i}\psi_{x_0}^{\langle 1,n\rangle}+\psi_{x_i}^{\langle 1,n\rangle}=\frac{\lambda x'_i}{x'_{[1,n]}}\psi_{x_0}^{\langle 1,n\rangle}+\psi_{x_i}^{\langle 1,n\rangle}, \quad (3.24_2)$$

so that (3.24) reduces to (3.22). Hence, $(3.17_1)$ is equivalent to (3.17).•

**Comment 3.7.** Like (3.24), the expression on the right-hand side of either one of equations (3.19) and (3.20) can be developed with the help of $(3.24_2)$ or with the help of the variant of $(3.24_2)$ with 'ϕ' in place of 'ψ'.•

**Lemma 3.5:** *The second (symmetric) differential tautology for a pair of biased functional forms.*

$$\sum_{\alpha\in\{0,1,3\}}K^{(\alpha)}\left(\phi^{\langle 1,n\rangle},\psi^{\langle 1,n\rangle},\tau_n,\overline{x}_{[1,n]}\right)=-K^{(2)}\left(\phi^{\langle 1,n\rangle},\psi^{\langle 1,n\rangle},\tau_n,\overline{x}_{[1,n]}\right)$$
$$=-K^{(2)}\left(\psi^{\langle 1,n\rangle},\phi^{\langle 1,n\rangle},\tau_n,\overline{x}_{[1,n]}\right)=\sum_{\alpha\in\{0,1,3\}}K^{(\alpha)}\left(\psi^{\langle 1,n\rangle},\phi^{\langle 1,n\rangle},\tau_n,\overline{x}_{[1,n]}\right) \quad (3.25)$$

**Proof:** Under the hypothesis of the lemma, the variant of Lemma 3.4 with '$\phi$' and '$\psi$' exchanged is valid. Therefore, besides (3.17),

$$\sum_{\alpha=0}^{3} K^{(\alpha)}\left(\psi^{\langle 1,n\rangle}, \phi^{\langle 1,n\rangle}, \tau_n, \overline{x}_{[1,n]}\right) = 0. \tag{3.26}$$

At the same time, it is seen from (3.20) that

$$K^{(2)}\left(\phi^{\langle 1,n\rangle}, \psi^{\langle 1,n\rangle}, \tau_n, \overline{x}_{[1,n]}\right) = K^{(2)}\left(\psi^{\langle 1,n\rangle}, \phi^{\langle 1,n\rangle}, \tau_n, \overline{x}_{[1,n]}\right). \tag{3.27}$$

Identity (3.23) immediately follows from (3.17) and (3.26) by (3.27).•

### 3.4. Two basic Integro-differential tautologies for biased functional forms

**Theorem 3.1:** *The first (asymmetric) basic integro-differential tautology for a pair of biased functional forms.* 1) For each $n \in \omega_1$:

$$\sum_{\alpha=0}^{3} \Omega^{(\alpha)}\left(\phi^{\langle 1,n\rangle}, \psi^{\langle 1,n\rangle}, \tau_n, \overline{X}_n^{\mathrm{c}}\right) = 0, \tag{3.28}$$

where for each $\alpha \in \omega_{0,3}$:

$$\Omega^{(\alpha)}\left(\phi^{\langle 1,n\rangle}, \psi^{\langle 1,n\rangle}, \tau_n, \overline{X}_n^{\mathrm{c}}\right) \doteq \int_{\overline{X}_n^{\mathrm{c}}} K^{(\alpha)}\left(\phi^{\langle 1,n\rangle}, \psi^{\langle 1,n\rangle}, \tau_n, \overline{x}_{[1,n]}\right) dv_{\{n\}} \tag{3.29}$$

subject to

$$dv_{\{n\}} \doteq dx_1 \cdots dx_n \tag{3.30}$$

and also subject to (2.58).

2) For each $n \in \omega_2$:

$$\begin{aligned}\Omega^{(3)}\left(\phi^{\langle 1,n\rangle}, \psi^{\langle 1,n\rangle}, \tau_n, \overline{X}_{\{n\}}^{\mathrm{c}}\right) &= -\int_{\overline{X}_n^{\mathrm{c}}} \sum_{i=1}^{n} \frac{\partial}{\partial x_i} L_i\left(\phi^{\langle 1,n\rangle}, \psi^{\langle 1,n\rangle}, \tau_n, \overline{x}_{[1,n]}\right) \\ &= -\Sigma\left(\phi^{\langle 1,n\rangle}, \psi^{\langle 1,n\rangle}, \tau_n, \overline{X}_n^{\mathrm{b}}\right),\end{aligned} \tag{3.31}$$

where, by the Ostrogradsky-Gauss theorem (see, e.g., Budak and Fomin [1978, pp. 218–224]),

$$\Sigma\left(\phi^{\langle 1,n\rangle}, \psi^{\langle 1,n\rangle}, \tau_n, \overline{X}_n^{\mathrm{b}}\right) \doteq \int_{\overline{X}_n^{\mathrm{b}}} \sum_{i=1}^{n} L_i\left(\phi^{\langle 1,n\rangle}, \psi^{\langle 1,n\rangle}, \tau_n, \overline{x}_{[1,n]}\right) n_i\left(\overline{x}_{[1,n]}\right) ds_{\{n\}} \tag{3.32}$$

subject to (3.22) or (3.24); $ds_{\{n\}}$ is an element of the $(n-1)$-dimensional area of the boundary surface $\overline{X}_n^{\mathrm{b}}$ at the integration point $\overline{x}_{[1,n]} \in \overline{X}_n^{\mathrm{b}}$, whereas the ordered *n*-tuple $\overline{n}_{[1,n]}\left(\overline{x}_{[1,n]}\right)$, defined as

$$\overline{n}_{[1,n]}\left(\overline{x}_{[1,n]}\right) \doteq \left\langle n_1\left(\overline{x}_{[1,n]}\right), \ldots, n_n\left(\overline{x}_{[1,n]}\right)\right\rangle, \tag{3.33}$$

is the unit normal vector of the surface $\overline{X}_n^b$ at the integration point *outward* with respect to $\overline{X}_n^c$.

a) In general, the spherical surface $B_n^b(r, \overline{x}_{*[1,n]})$ of a radius $r$ is said to be *positive with respect to its interior* $B_n^o(r, \overline{x}_{*[1,n]})$ and is alternatively denoted by '$B_n^{b+}(r, \overline{x}_{*[1,n]})$',, and it is said to be *negative* and is alternatively denoted by '$B_n^{b-}(r, \overline{x}_{*[1,n]})$' with respect to its exterior $\overline{E}_n - \overline{B}_n^c(r, \overline{x}_{*[1,n]})$ in the sense that all its unit outward vectors are *centrifugal* in the former case and *centripetal* in the latter case; i.e.

$$n_i(\overline{x}_{[1,n]}) = \begin{cases} -\dfrac{x_i'}{x_{[1,n]}'} = -\dfrac{x_i - x_{*i}}{r} & \text{if } \overline{x}_{[1,n]} \in B_n^{b-}(r, \overline{x}_{*[1,n]}) \text{ (a)} \\ \dfrac{x_i'}{x_{[1,n]}'} = \dfrac{x_i - x_{*i}}{r} & \text{if } \overline{x}_{[1,n]} \in B_n^{b+}(r, \overline{x}_{*[1,n]}) \text{ (b)} \end{cases}. \quad (3.34)$$

Accordingly, the boundary spherical surfaces $B_n^b(r_1, \overline{x}_{*[1,n]})$ and $B_n^b(r_2, \overline{x}_{*[1,n]})$ of $\overline{C}_n^{cc}(r_1, r_2, \overline{x}_{*[1,n]})$ are alternatively denoted by '$B_n^{b-}(r_1, \overline{x}_{*[1,n]})$' and '$B_n^{b+}(r_2, \overline{x}_{*[1,n]})$', because

$$n_i(\overline{x}_{[1,n]}) = \begin{cases} -\dfrac{x_i'}{x_{[1,n]}'} = -\dfrac{x_i - x_{*i}}{r_1} & \text{if } \overline{x}_{[1,n]} \in B_n^{b-}(r_1, \overline{x}_{*[1,n]}) \text{ (a)} \\ \dfrac{x_i'}{x_{[1,n]}'} = \dfrac{x_i - x_{*i}}{r_2} & \text{if } \overline{x}_{[1,n]} \in B_n^{b+}(r_2, \overline{x}_{*[1,n]}) \text{ (b)} \end{cases}. \quad (3.35)$$

b) By (2.58) and (3.35), it follows from (3.32) that

$$\Sigma(\phi^{\langle 1,n \rangle}, \psi^{\langle 1,n \rangle}, \tau_n, \overline{X}_n^b) = \Sigma(\phi^{\langle 1,n \rangle}, \psi^{\langle 1,n \rangle}, \tau_n, \overline{C}_n^{cc}(r_1, r_2, \overline{x}_{*[1,n]}))$$
$$= \Sigma(\phi^{\langle 1,n \rangle}, \psi^{\langle 1,n \rangle}, \tau_n, \overline{B}_n^{b-}(r_1, \overline{x}_{*[1,n]})) + \Sigma(\phi^{\langle 1,n \rangle}, \psi^{\langle 1,n \rangle}, \tau_n, B_n^{b+}(r_2, \overline{x}_{*[1,n]})), \quad (3.36)$$

where

$$\Sigma(\phi^{\langle 1,n \rangle}, \psi^{\langle 1,n \rangle}, \tau_n, \overline{B}_n^{b-}(r_1, \overline{x}_{*[1,n]}))$$
$$= \int_{\overline{B}_n^{b-}(r_1, \overline{x}_{*[1,n]})} \sum_{i=1}^{n} L_i(\phi^{\langle 1,n \rangle}, \psi^{\langle 1,n \rangle}, \tau_n, \overline{x}_{[1,n]}) n_i(\overline{x}_{[1,n]}) ds_{\{n\}} \quad (3.37)$$
$$= -\dfrac{1}{r_1} \int_{\overline{B}_n^{b}(r_1, \overline{x}_{*[1,n]})} \sum_{i=1}^{n} L_i(\phi^{\langle 1,n \rangle}, \psi^{\langle 1,n \rangle}, \tau_n, \overline{x}_{[1,n]})(x_i - x_{*i}) ds_{\{n\}},$$

$$\Sigma(\phi^{\langle 1,n \rangle}, \psi^{\langle 1,n \rangle}, \tau_n, \overline{B}_n^{b+}(r_2, \overline{x}_{*[1,n]}))$$
$$= \int_{\overline{B}_n^{b+}(r_2, \overline{x}_{*[1,n]})} \sum_{i=1}^{n} L_i(\phi^{\langle 1,n \rangle}, \psi^{\langle 1,n \rangle}, \tau_n, \overline{x}_{[1,n]}) n_i(\overline{x}_{[1,n]}) ds_{\{n\}} \quad (3.38)$$
$$= -\dfrac{1}{r_2} \int_{\overline{B}_n^{b}(r_2, \overline{x}_{*[1,n]})} \sum_{i=1}^{n} L_i(\phi^{\langle 1,n \rangle}, \psi^{\langle 1,n \rangle}, \tau_n, \overline{x}_{[1,n]})(x_i - x_{*i}) ds_{\{n\}}.$$

3) For $n=1$:

$$\begin{aligned}
&\Omega^{(3)}\left(\phi^{\langle 1,1\rangle},\psi^{\langle 1,1\rangle},\tau_1,\overline{X}_1^c\right) \\
&= \Omega^{(3)}\left(\phi^{\langle 1,1\rangle},\psi^{\langle 1,1\rangle},\tau_1,[x_{*1}-r_2,x_{*1}+r_2]\setminus(x_{*1}-r_1,x_{*1}+r_1)\right) \\
&= \Omega^{(3)}\left(\phi^{\langle 1,1\rangle},\psi^{\langle 1,1\rangle},\tau_1,[x_{*1}-r_2,x_{*1}+r_2]\right) \\
&\quad - \Omega^{(3)}\left(\phi^{\langle 1,1\rangle},\psi^{\langle 1,1\rangle},\tau_1,[x_{*1}+r_1,x_{*1}+r_1]\right),
\end{aligned} \quad (3.39)$$

where

$$\begin{aligned}
\Omega^{(3)}\left(\phi^{\langle 1,1\rangle},\psi^{\langle 1,1\rangle},\tau_1,[x_{*1}-r_1,x_{*1}+v]\right) &\equiv -\int_{x_{*1}-r_1}^{x_{*1}+r_1}\frac{\partial}{\partial x_1}L_1\left(\phi^{\langle 1,1\rangle},\psi^{\langle 1,1\rangle},\tau_1,x_1\right)dx_1 \\
&= L_1\left(\phi^{\langle 1,1\rangle},\psi^{\langle 1,1\rangle},\tau_1,x_{*1}-r_1\right)-L_1\left(\phi^{\langle 1,1\rangle},\psi^{\langle 1,1\rangle},\tau_1,x_{*1}+r_1\right),
\end{aligned} \quad (3.40)$$

$$\begin{aligned}
\Omega^{(3)}\left(\phi^{\langle 1,1\rangle},\psi^{\langle 1,1\rangle},\tau_1,[x_{*1}-r_2,x_{*1}+r_2]\right) &\equiv -\int_{x_{*1}-r_2}^{x_{*1}+r_2}\frac{\partial}{\partial x_1}L_1\left(\phi^{\langle 1,1\rangle},\psi^{\langle 1,1\rangle},\tau_1,x_1\right)dx_1 \\
&= L_1\left(\phi^{\langle 1,1\rangle},\psi^{\langle 1,1\rangle},\tau_1,x_{*1}-r_2\right)-L_1\left(\phi^{\langle 1,1\rangle},\psi^{\langle 1,1\rangle},\tau_1,x_{*1}+r_2\right).
\end{aligned} \quad (3.41)$$

**Proof:** 1) Item 1 of the theorem follows straightforwardly from Lemma 3.4 by Hypothesis 3.1.

2) The statement of item 2 of the theorem includes its proof.

3) At $n=1$ and $\alpha=3$, definition (3.29) subject to (3.21) (with allowance for (3.22)) can straightforwardly be developed thus:

$$\begin{aligned}
&\Omega^{(3)}\left(\phi^{\langle 1,1\rangle},\psi^{\langle 1,1\rangle},\tau_1,\overline{X}_n^c\right) \\
&= \Omega^{(3)}\left(\phi^{\langle 1,1\rangle},\psi^{\langle 1,1\rangle},\tau_1,[x_{*1}-r_2,x_{*1}-r_1]\cup[x_{*1}+r_1,x_{*1}+r_2]\right) \\
&= \Omega^{(3)}\left(\phi^{\langle 1,1\rangle},\psi^{\langle 1,1\rangle},\tau_1,[x_{*1}-r_2,x_{*1}-r_1]\right) \\
&\quad + \Omega^{(3)}\left(\phi^{\langle 1,1\rangle},\psi^{\langle 1,1\rangle},\tau_1,[x_{*1}+r_1,x_{*1}+r_2]\right)
\end{aligned} \quad (3.42)$$

subject to

$$\begin{aligned}
\Omega^{(3)}\left(\phi^{\langle 1,1\rangle},\psi^{\langle 1,1\rangle},\tau_1,[x_{*1}-r_2,x_{*1}-r_1]\right) &\equiv -\int_{x_{*1}-r_2}^{x_{*1}-r_1}\frac{\partial}{\partial x_1}L_1\left(\phi^{\langle 1,1\rangle},\psi^{\langle 1,1\rangle},\tau_1,x_1\right)dx_1 \\
&= L_1\left(\phi^{\langle 1,1\rangle},\psi^{\langle 1,1\rangle},\tau_1,x_{*1}-r_2\right)-L_1\left(\phi^{\langle 1,1\rangle},\psi^{\langle 1,1\rangle},\tau_1,x_{*1}-r_1\right),
\end{aligned} \quad (3.43)$$

$$\begin{aligned}
\Omega^{(3)}\left(\phi^{\langle 1,1\rangle},\psi^{\langle 1,1\rangle},\tau_1,[x_{*1}+r_1,x_{*1}+r_2]\right) &\equiv -\int_{x_{*1}+r_1}^{x_{*1}+r_2}\frac{\partial}{\partial x_1}L_1\left(\phi^{\langle 1,1\rangle},\psi^{\langle 1,1\rangle},\tau_1,x_1\right)dx_1 \\
&= L_1\left(\phi^{\langle 1,1\rangle},\psi^{\langle 1,1\rangle},\tau_1,x_{*1}+r_1\right)-L_1\left(\phi^{\langle 1,1\rangle},\psi^{\langle 1,1\rangle},\tau_1,x_{*1}+r_2\right),
\end{aligned} \quad (3.44)$$

which is equivalent to (3.39)–(3.41).•

**Comment 3.8.** For the equality (3.31) subject to (3.32) be valid, it is sufficient that the domain $\overline{X}_n^c$ should be *regular in each coordinate direction* (see, e.g., Budak and Fomin [1978, pp. 90, 219]). Owing to (2.58), this sufficient condition is satisfied.•

**Theorem 3.2:** *The second (symmetric) basic integro-differential tautology for a pair of biased functional forms.*

$$\sum_{\alpha \in \{0,1,3\}} \Omega^{(\alpha)}\left(\phi^{\langle 1,n \rangle}, \psi^{\langle 1,n \rangle}, \tau_n, \overline{X}_n^c\right) = \sum_{\alpha \in \{0,1,3\}} \Omega^{(\alpha)}\left(\psi^{\langle 1,n \rangle}, \phi^{\langle 1,n \rangle}, \tau_n, \overline{X}_n^c\right) \quad (3.45)$$

subject to the notations of Theorem 3.1 and of the variant of Theorem 3.1 with '$\phi^{\langle 1,n \rangle}$' and '$\psi^{\langle 1,n \rangle}$' exchanged.

**Proof:** From the tautology (3.28) subject to (3.29) and from the variant of that tautology with '$\phi^{\langle 1,n \rangle}$' and '$\psi^{\langle 1,n \rangle}$', it follows by (3.25) that

$$\sum_{\alpha \in \{0,1,3\}} \Omega^{(\alpha)}\left(\phi^{\langle 1,n \rangle}, \psi^{\langle 1,n \rangle}, \tau_n, \overline{X}_n^c\right) = -\Omega^{(2)}\left(\phi^{\langle 1,n \rangle}, \psi^{\langle 1,n \rangle}, \tau_n, \overline{X}_n^c\right)$$
$$= -\Omega^{(2)}\left(\psi^{\langle 1,n \rangle}, \phi^{\langle 1,n \rangle}, \tau_n, \overline{X}_n^c\right) = \sum_{\alpha \in \{0,1,3\}} \Omega^{(\alpha)}\left(\psi^{\langle 1,n \rangle}, \phi^{\langle 1,n \rangle}, \tau_n, \overline{X}_n^c\right), \quad (3.46)$$

which immediately implies the tautology (3.45).•

# 4. Special basic integro-differential tautologies for a single biased functional form

## 4.1. A spherically symmetric harmonic function in $\overline{E}_{\{n\}}$

**Lemma 4.1.** Given $n \in \omega_1$, given $p \in I_{-\infty,\infty}$, for each $\overline{x}_{[1,n]} \in \overline{E}_n - \{\overline{0}_{[1,n]}\}$:

$$\Delta_{[1,n]} x_{[1,n]}^{\ p} = p(p+n-2) x_{[1,n]}^{\ p-2}, \quad (4.1)$$

subject to

$$x_{[1,n]} \triangleq \left|\overline{x}_{[1,n]}\right| = \sqrt{\sum_{i=1}^{n} x_i^{\ 2}} > 0 \quad (4.2)$$

(cf. (2.43) and (2.44)) and also subject to (3.11).

**Proof:** It follows from (3.11) and (4.2) that for each $i \in \omega_{1,n}$ and each $j \in \omega_{1,n}$:

$$\nabla_i x_{[1,n]}^{\ p} = p x_{[1,n]}^{\ p-1} \nabla_i x_{[1,n]} = p x_{[1,n]}^{\ p-2} x_i, \quad (4.3)$$

$$\nabla_i \nabla_j x_{[1,n]}^{\ p} = p \nabla_j \left[x_{[1,n]}^{\ p-2} x_i\right] = p \left[(p-2) x_{[1,n]}^{\ -2} x_i x_j + \delta_{ij}\right] x_{[1,n]}^{\ p-2} \quad (4.4)$$

(cf. (3.6)). Summing up each term of (4.4) at $j=i$ over all values of '$i$' from 1 to $n$ yields (4.1), by (4.2).•

**Theorem 4.1.** Given $n \in \omega_1$, for each $\overline{x}_{[1,n]} \in \overline{E}_n - \{\overline{0}_{[1,n]}\}$:

$$\Delta_{[1,n]} x_{[1,n]}^{\ 2-n} = 0 \text{ if } n \in \omega_1 - \{2\}, \quad (4.5)$$

$$\Delta_{[1,2]} \ln x_{[1,2]} = 0 \text{ if } n = 2, \tag{4.6}$$

subject to (3.11).

**Proof:** Under the hypothesis of Lemma 4.1, it immediately follows from (4.1) that for each $\bar{x}_{[1,n]} \in \bar{E}_n - \{\bar{0}_{[1,n]}\}$:

$$\Delta_{[1,n]} x_{[1,n]}^{\ p} = 0 \text{ if and only if } p=0 \text{ or } p=2-n. \tag{4.5$_1$}$$

Hence, relation (4.5$_1$) yields the nontrivial equation (4.5) if $p=2-n\neq 0$ and it yields the trivial equation $\Delta_{[1,n]} 1 = 0$ if $p=0$ or if particularly $p=2-n=0$. At $n=2$, the nontrivial equation (4.6) is proved by the following straightforward computations which are analogous to (4.3) and (4.4). For each $i \in \omega_{1,2}$: for each $j \in \omega_{1,2}$:

$$\nabla_i \ln x_{[1,2]} = x_{[1,2]}^{\ -2} x_i, \tag{4.7}$$

$$\nabla_i \nabla_j \ln x_{[1,2]} = \nabla_j x_{[1,2]}^{\ -1} x_i = x_{[1,2]}^{\ -2}\left(x_{[1,2]}^{\ 2}\delta_{ij} - 2x_i x_j\right). \tag{4.8}$$

Summing up each term of (4.8) at $j = i$ over the values 1 and 2 of '$i$' yields (4.6).•

**Definition 4.2.** 1) Given $n \in \omega_1$, for each $\bar{x}_{[1,n]} \in \bar{E}_n - \{\bar{x}_{*[1,n]}\}$:

a) If $n \geq 3$ then

$$\eta^{\langle n \rangle}\!\left(x'_{[1,n]}\right) \stackrel{\scriptscriptstyle\equiv}{=} a_n x'^{\ 2-n}_{[1,n]}, \tag{4.9}$$

where

$$a_n \stackrel{\scriptscriptstyle\equiv}{=} \frac{1}{(n-2)S_n(1)} = \frac{\Gamma(n/2)}{2(n-2)\pi^{n/2}}. \tag{4.10}$$

b) If $n=2$ then

$$\eta^{\langle n \rangle}\!\left(x'_{[1,2]}\right) \stackrel{\scriptscriptstyle\equiv}{=} -a_2 \ln x'_{[1,2]}, \tag{4.11}$$

where

$$a_2 \stackrel{\scriptscriptstyle\equiv}{=} \frac{1}{S_2(1)} = \frac{1}{2\pi}. \tag{4.12}$$

c) If $n=1$ then

$$\eta^{\langle 1 \rangle}\!\left(x'_{[1,1]}\right) \stackrel{\scriptscriptstyle\equiv}{=} a_1 x'_{[1,1]} = a_1|x'_1| = a_1|x_1 - x_{*1}|, \tag{4.13}$$

where

$$a_1 \stackrel{\scriptscriptstyle\equiv}{=} -\frac{1}{2}. \tag{4.14}$$

2) It is understood that $S_n(1)$, defined as

$$S_n(1) = nV_n(1) = \frac{2\pi^{n/2}}{\Gamma(n/2)} \text{ for each } n \in \omega_2. \quad (4.15)$$

is the area of an ($n$–1)-dimensional spherical surface in $\overline{E}_n$ of radius 1 (see Appendix A), i.e. the area of $\overline{B}_n^b(1, x_{*[1,n]})$, while $V_n(1)$ is the volume of $\overline{B}_n^c(1, x_{*[1,n]})$ and $\Gamma$ is the gamma-function. It is also understood that '$x'_{[1,n]}$' is defined by (2.44). Accordingly, (4.13) is the instance of (4.9) at $n=1$. At the same time, $\overline{B}_1^b(1, x_{*[1,n]}) = \{-1,1\}$, so that '$S_{\{1\}}(1)$' is not defined. It is, however, convenient to utilize (4.10) at $n=1$ with allowance for (1.13) and for $\Gamma(1/2) = \pi^{n/2}$ in order to formally define '$S_1(1)$' as:

$$S_1(1) \equiv 2. \quad (4.16)$$

This definition can be understood as expressing the fact that the boundary $\overline{B}_1^b(1, x_{*[1,n]})$ consists of *two points*: –1 and 1.

**Theorem 4.2.** For each $n \in \omega_1$, independent of $a_n$,

$$\Delta_{[1,n]} \eta^{\langle n \rangle}(x'_{[1,n]}) \equiv \sum_{i=1}^{n} \frac{\partial^2}{\partial x_i^2} \eta^{\langle n \rangle}(x'_{[1,n]}) = 0 \text{ for each } \overline{x}_{[1,n]} \in \overline{E}_{\{n\}} - \{\overline{x}_{*[1,n]}\}, \quad (4.17)$$

whereas owing to (4.10), (4.12), and (4.14),

$$\Delta_{[1,n]} \eta^{\langle n \rangle}(x'_{[1,n]}) = \delta^{\langle n \rangle}(\overline{x}'_{[1,n]}) \text{ for each } \overline{x}_{[1,n]} \in \overline{E}_n, \quad (4.18)$$

where $\delta^{\langle n \rangle}$ is the *n*-ary Dirac $\delta$-function. Accordingly, the *associated function* $\eta^{\langle n \rangle}$ of the functional form $\eta^{\langle n \rangle}(x'_{[1,n]})$ is a *spherically symmetrical harmonic* if $\eta^{\langle n \rangle}(x'_{[1,n]})$ applies in an *n*-dimensional domain of $\overline{E}_{\{n\}}$ that *does not contain* its *singular point* $\overline{x}_{*[1,n]}$, and $\eta^{\langle n \rangle}$ becomes the *Green function of an n-dimensional Poisson (inhomogeneous Laplace) equation* if $\eta^{\langle n \rangle}(x'_{[1,n]})$ applies in an *n*-dimensional domain of $\overline{E}_n$ that *contains* its singular point $\overline{x}_{*[1,n]}$.

**Proof:** (4.17) follows from Theorem 4.1 by Definition 4.2. (4.18) is proved in Appendix B.•

### 4.2. Specification of the integrands occurring in the basic binary integro-differential tautologies

**Definition 4.3.** Given $n \in \omega_1$, given $\overline{x}_{*[1,n]} \in \overline{E}_n$, for each $x_0 \in R$, for each $\overline{x}_{[1,n]} \in \overline{E}_n - \{\overline{x}_{*[1,n]}\}$, henceforth

$$\phi^{\langle 1,n \rangle}(x_0, \overline{x}_{[1,n]}) \equiv \eta^{\langle n \rangle}(x'_{[1,n]}) = \eta^{\langle n \rangle}(|\overline{x}_{[1,n]} - \overline{x}_{*[1,n]}|) \quad (4.19)$$

subject to (2.9)–(2.11), (2.43), (2.44), (4.10), (4.12), and (4.14).

**Corollary 4.1.** By definition (4.19), '$\phi^{\langle 1,n \rangle}(x_0, \bar{x}_{[1,n]})$' is a real-valued functional form, which along with all its partial derivatives with respect to '$x_0$', '$x_1$', '$x_n$' is defined and continuous for each $\langle x_0, \bar{x}_{[1,n]} \rangle \in R \times [\bar{E}_{\{n\}} - \{\bar{x}_{*[1,n]}\}]$, so that symbolically

$$\phi^{\langle 1,n \rangle} \in \mathcal{C}_\infty(R \times [\bar{E}_{\{n\}} - \{\bar{x}_{*[1,n]}\}]), \qquad (4.20)$$

where $\phi^{\langle 1,n \rangle}$ is *the full associated function of the functional form* '$\phi^{\langle 1,n \rangle}(x_0, \bar{x}_{[1,n]})$'.•

**Definition 4.4.** Given $n \in \omega_1$, given a real-valued function $\psi^{\langle 1,n \rangle}$ satisfying (3.1), for each $x_0 \in R$: for each $\bar{x}_{[1,n]} \in \bar{X}_n^c$:

$$f^{\langle 1,n \rangle}(x_0, \bar{x}_{[1,n]}) \triangleq -D_{[0,n]} \psi^{\langle 1,n \rangle}(x_0, \bar{x}_{[1,n]}), \qquad (4.21)$$

subject to (3.12).•

**Convention 4.1.** Besides Conventions 3.2 and 3.3, which remain in force, we shall use the following abbreviations analogous to those introduced in Definition 3.2 and Convention 3.3(1):

$$f^{\langle 1,n \rangle} \triangleq f^{\langle 1,n \rangle}(x_{*0} + \lambda x'_{[1,n]}, \bar{x}_{[1,n]}) \triangleq -\left[D_{[0,n]} \psi^{\langle 1,n \rangle}(x_0, \bar{x}_{[1,n]})\right]_{x_0 = x_{*0} + \lambda x'_{[1,n]}}, \qquad (4.22)$$

$$\eta^{\langle n \rangle} \triangleq \eta^{\langle n \rangle}(x'_{[1,n]}), \eta^{\langle n \rangle}_{x_i} \triangleq \eta^{\langle n \rangle}_{x_i}(x'_{[1,n]}) \triangleq \frac{\partial \eta^{\langle n \rangle}(x'_{[1,n]})}{\partial x_i} \text{ for each } i \in \omega_{1,n}. \qquad (4.23)$$

**Lemma 4.2.**

$$K^{(0)}(\eta^{\langle n \rangle}, \psi^{\langle 1,n \rangle}, \tau_n, \bar{x}_{[1,n]}) = -\eta^{\langle n \rangle} f^{\langle 1,n \rangle}; \qquad (4.24)$$

$$K^{(1)}(\eta^{\langle n \rangle}, \psi^{\langle 1,n \rangle}, \tau_n, \bar{x}_{[1,n]}) \triangleq -\frac{\lambda}{x'_{[1,n]}} \left[(n-1)\eta^{\langle n \rangle} + 2\sum_{i=1}^n x'_i \eta^{\langle n \rangle}_{x_i}\right] \psi^{\langle 1,n \rangle}_{x_0}$$

$$= \begin{cases} \dfrac{\lambda(n-3)a_n}{x'^{n-1}_{[1,n]}} \psi^{\langle 1,n \rangle}_{x_0} & \text{if } n \neq 2 \text{ (a)} \\ \dfrac{\lambda a_2 (\ln x'_{[1,2]} + 2)}{x'_{[1,2]}} \psi^{\langle 1,2 \rangle}_{x_0} & \text{if } n = 2 \text{ (b)} \end{cases}; \qquad (4.25)$$

$$K^{(3)}(\eta^{\langle n \rangle}, \psi^{\langle 1,n \rangle}, \tau_n, \bar{x}_{[1,n]}) \triangleq -\sum_{i=1}^n \frac{\partial}{\partial x_i} L_i(\eta^{\langle n \rangle}, \psi^{\langle 1,n \rangle}, \tau_n, \bar{x}_{[1,n]}), \qquad (4.26)$$

subject to

$$L_i(\eta^{\langle n \rangle}, \psi^{\langle 1,n \rangle}, \tau_n, \bar{x}_{[1,n]}) \triangleq \eta^{\langle n \rangle} \left[\psi^{\langle 1,n \rangle}_{x_i} - \frac{\lambda x'_i}{x'_{[1,n]}} \psi^{\langle 1,n \rangle}_{x_0}\right] \text{ for each } i \in \omega_{1,n}; \qquad (4.27)$$

$$K^{(0)}\left(\psi^{\langle 1,n\rangle}, \tau_n, \eta^{\langle n\rangle}, \overline{x}_{[1,n]}\right) = K^{(1)}\left(\psi^{\langle 1,n\rangle}, \tau_n, \eta^{\langle n\rangle}, \overline{x}_{[1,n]}\right) = 0; \tag{4.28}$$

$$K^{(3)}\left(\psi^{\langle 1,n\rangle}, \eta^{\langle n\rangle}, \tau_n, \overline{x}_{[1,n]}\right) \equiv -\sum_{i=1}^{n} \frac{\partial}{\partial x_i} L_i\left(\psi^{\langle 1,n\rangle}, \eta^{\langle n\rangle}, \tau_n, \overline{x}_{[1,n]}\right), \tag{4.29}$$

subject to

$$L_i\left(\psi^{\langle 1,n\rangle}, \eta^{\langle n\rangle}, \tau_n, \overline{x}_{[1,n]}\right) \equiv \psi^{\langle 1,n\rangle} \eta^{\langle n\rangle}_{x_i} = \begin{cases} \dfrac{(2-n)a_n x_i' \psi^{\langle 1,n\rangle}}{x_{[1,n]}'^{\,n}} & \text{if } n \neq 2 \text{ (a)} \\ -\dfrac{a_2 x_i' \psi^{\langle 1,2\rangle}}{x_{[1,2]}'^{\,2}} & \text{if } n = 2 \text{ (b)} \end{cases} \tag{4.30}$$

for each $i \in \omega_{1,n}$;

$$K^{(2)}\left(\eta^{\langle n\rangle}, \psi^{\langle 1,n\rangle}, \tau, \overline{x}_{[1,n]}\right) = K^{(2)}\left(\psi^{\langle 1,n\rangle}, \eta^{\langle n\rangle}, \tau, \overline{x}_{[1,n]}\right) = \sum_{i=1}^{n} \eta^{\langle n\rangle}_{x_i} \psi^{\langle 1,n\rangle}_{x_i}$$

$$= \begin{cases} \dfrac{(2-n)a_n}{x_{[1,n]}'^{\,n}} \sum_{i=1}^{n} x_i' \psi^{\langle 1,n\rangle}_{x_i} & \text{if } n \neq 2 \text{ (a)} \\ -\dfrac{a_2}{x_{[1,2]}'^{\,2}} \sum_{i=1}^{2} x_i' \psi^{\langle 1,2\rangle}_{x_i} & \text{if } n = 2 \text{ (b)} \end{cases}. \tag{4.31}$$

**Proof:** 1) Substitution of (4.19) into (3.18), (3.21) and (3.22) yields (4.24) (subject to (4.22)), (4.26), and (4.27) respectively, whereas the same substitution into (3.19) yields the first, general expression for $K^{(1)}\left(\eta^{\langle n\rangle}, \psi^{\langle 1,n\rangle}, \tau_n, \overline{x}_{[1,n]}\right)$ in (4.25). At the same time, by (2.43)-(2.46), (4.3), and (4.7), it follows from (4.9), (4.11), and (4.13) that for each $i \in \omega_{1,n}$:

$$\eta^{\langle n\rangle}_{x_i} = \frac{\partial \eta^{\langle n\rangle}\left(x_{[1,n]}'\right)}{\partial x_i} = \begin{cases} a_n \dfrac{\partial x_{[1,n]}'^{\,2-n}}{\partial x_i} = \dfrac{(2-n)a_n x_i'}{x_{[1,n]}'^{\,n}} & \text{if } n \neq 2 \text{ (a)} \\ -a_2 \dfrac{\partial \ln x_{[1,2]}'}{\partial x_i} = -\dfrac{a_2 x_i'}{x_{[1,2]}'^{\,2}} & \text{if } n = 2 \text{ (b)} \end{cases} \tag{4.32}$$

and hence

$$\sum_{i=1}^{n} x_i' \eta^{\langle n\rangle}_{x_i} = \begin{cases} \dfrac{(2-n)a_n}{x_{[1,n]}'^{\,n}} \sum_{i=1}^{n} x_i' x_i' = \dfrac{(2-n)a_n}{x_{[1,n]}'^{\,n-2}} & \text{if } n \neq 2 \text{ (a)} \\ -\dfrac{a_2}{x_{[1,2]}'^{\,2}} \sum_{i=1}^{2} x_i' x_i' = -a_2 & \text{if } n = 2 \text{ (b)} \end{cases} \tag{4.33}$$

At $n \neq 2$, the expression in the square brackets in (4.25) reduces with the help of (4.9) and (4.33a) thus:

$$(n-1)\eta^{\langle n\rangle} + 2\sum_{i=1}^{n} x_i' \eta^{\langle n\rangle}_{x_i} = \frac{[(n-1)+2(2-n)]a_n}{x_{[1,n]}'^{\,n-2}} = -\frac{(n-3)a_n}{x_{[1,n]}'^{\,n-2}}. \tag{4.34a}$$

At $n = 2$, that expression reduces with the help of (4.11) and (4.33b) thus:

$$\eta^{\langle 2\rangle} + 2\sum_{i=1}^{2} x'_i \eta^{\langle 2\rangle}_{x_i} = -a_2\left[\ln x'_{[1,2]} + 2\right]. \tag{4.34b}$$

The final result in (4.25) is immediately obtained with the help of (4.34a) and (4.34b).

2) By Definitions 3.3 and 4.2 and by Theorem 4.2, it follows that for each $\bar{x}_{[1,n]} \in \bar{E}_{\{n\}} - \{\bar{x}_{*[1,n]}\}$:

$$D_{[0,n]}\eta^{\langle n\rangle}(x'_{[1,n]}) = \Delta_{[1,n]}\eta^{\langle n\rangle}(x'_{[1,n]}) = 0, \tag{4.35}$$

$$\frac{\partial \eta^{\langle n\rangle}(x'_{[1,n]})}{\partial x_0} = 0. \tag{4.36}$$

Equations (4.28)–(4.31) follow from the variants of (3.18)–(3.22) with '$\phi$' and '$\psi$' exchanged subject to (4.19) by (4.35) and (4.36); the final result in (4.31) is deduced immediately with the help of (4.32).•

**Comment 4.1.** According to (4.25a), the identity (tautology)

$$K^{(1)}\left(\eta^{\langle n\rangle}, \psi^{\langle 1,n\rangle}, \tau_3, \bar{x}_{[1,n]}\right) = 0 \text{ for each } \bar{x}_{[1,3]} \in \bar{E}_3 \tag{4.25!}$$

is an exclusive fundamental property of $\bar{E}_3$, which implies that $\psi^{\langle 1,n\rangle}(x_{*0} + \lambda x'_{[1,n]}, \bar{x}_{[1,n]})$ can satisfy the non-dispersive causality principle and which also implies that special theory of relativity is possible if and only if $n=3$, i.e. only for a three-dimensional affine Euclidean space. These implications will be explicated in due course later on.•

**Comment 4.2.** Proceeding from (3.24) instead of (3.22), one arrives at (4.24) subject to

$$L_i\left(\phi^{\langle 1,n\rangle}, \psi^{\langle 1,n\rangle}, \tau_n, \bar{x}_{[1,n]}\right) \equiv \eta^{\langle n\rangle}\left[\frac{\partial \psi^{\langle 1,n\rangle}}{\partial x_i} - \frac{2\lambda x'_i}{x'_{[1,n]}}\psi^{\langle 1,n\rangle}_{x_0}\right] \text{ for each } i \in \omega_{1,n}, \tag{4.37}$$

instead of (4.25), whereas the variant of (3.24) with '$\phi$' and '$\psi$' exchanged subject to (4.19) reduces to (4.31) by (4.36). Equation (4.25) is simpler and more convenient in further applications than the equivalent equation (4.37).•

### 4.3. Two special basic integro-differential tautologies for a single biased functional form

**Theorem 4.3**: *The second special integro-differential tautology for a single biased functional form.* Given $n \in \omega_1$,

$$\begin{aligned}&\Omega^{(3)}\left(\psi^{\langle 1,n\rangle}, \eta^{\langle n\rangle}, \tau_n, \bar{C}_n^{cc}(r_1, r_2, \bar{x}_{*[1,n]})\right) \\ &= \sum_{\alpha \in \{0,1,3\}} \Omega^{(\alpha)}\left(\eta^{\langle n\rangle}, \psi^{\langle 1,n\rangle}, \tau_n, \bar{C}_n^{cc}(r_1, r_2, \bar{x}_{*[1,n]})\right)\end{aligned} \tag{4.38}$$

subject to the following definitions.

1) From the instance of (3.31) subject to (3.32) (see particularly (3.36) subject to (3.37) and (3.38)) with $\psi^{\langle 1,n \rangle}$ and $\eta^{\langle n \rangle}$ in place of $\phi^{\langle 1,n \rangle}$ and $\psi^{\langle 1,n \rangle}$ respectively, it follows that for each $n \in \omega_2$:

$$\Omega^{(3)}\left(\psi^{\langle 1,n \rangle}, \eta^{\langle n \rangle}, \tau_n, \overline{C}_n^{cc}(r_1, r_2, \overline{x}_{*[1,n]})\right) = -\Sigma\left(\psi^{\langle 1,n \rangle}, \eta^{\langle n \rangle}, \tau_n, \overline{C}_n^{b}(r_1, r_2, \overline{x}_{*[1,n]})\right) \quad (4.39)$$

subject to

$$\Sigma\left(\psi^{\langle 1,n \rangle}, \eta^{\langle n \rangle}, \tau_n, \overline{C}_n^{b}(r_1, r_2, \overline{x}_{*[1,n]})\right)$$
$$= \Sigma\left(\psi^{\langle 1,n \rangle}, \eta^{\langle n \rangle}, \tau_n, \overline{B}_n^{b-}(r_1, \overline{x}_{*[1,n]})\right) + \Sigma\left(\psi^{\langle 1,n \rangle}, \eta^{\langle n \rangle}, \tau_n, \overline{B}_n^{b+}(r_2, \overline{x}_{*[1,n]})\right), \quad (4.40)$$

where

$$\Sigma\left(\psi^{\langle 1,n \rangle}, \eta^{\langle n \rangle}, \tau_n, \overline{B}_n^{b-}(r_1, \overline{x}_{*[1,n]})\right)$$
$$= -\frac{1}{r_1} \int_{\overline{B}_n^b(r_1, \overline{x}_{*[1,n]})} \sum_{i=1}^n L_i\left(\psi^{\langle 1,n \rangle}, \eta^{\langle n \rangle}, \tau_n, \overline{x}_{[1,n]}\right) x_i'\, ds_{\{n\}}, \quad (4.41)$$

$$\Sigma\left(\psi^{\langle 1,n \rangle}, \eta^{\langle n \rangle}, \tau_n, \overline{B}_n^{b+}(r_2, \overline{x}_{*[1,n]})\right)$$
$$= \frac{1}{r_2} \int_{\overline{B}_n^b(r_2, \overline{x}_{*[1,n]})} \sum_{i=1}^n L_i\left(\psi^{\langle 1,n \rangle}, \eta^{\langle n \rangle}, \tau_n, \overline{x}_{[1,n]}\right) x_i'\, ds_{\{n\}}, \quad (4.42)$$

which are the pertinent instances of (3.37) and (3.38) respectively. By (4.30), it follows that

$$\sum_{i=1}^n L_i\left(\psi^{\langle 1,n \rangle}, \eta^{\langle n \rangle}, \tau_n, \overline{x}_{[1,n]}\right) x_i'$$
$$= \begin{cases} \dfrac{a_n(2-n)}{x_{[1,n]}'^{\,n-2}} \psi^{\langle 1,n \rangle}\left(x_{*0} + \lambda x_{[1,n]}', \overline{x}_{[1,n]}\right) & \text{if } n > 2 \text{ (a)} \\ -a_2 \psi^{\langle 1,2 \rangle}\left(x_{*0} + \lambda x_{[1,2]}', \overline{x}_{[1,2]}\right) & \text{if } n = 2 \text{ (b)} \end{cases} \quad (4.43)$$

Hence, by (4.10), (4.12), and (4.43), equations (4.41) and (4.42) become

$$\Sigma\left(\psi^{\langle 1,n \rangle}, \eta^{\langle n \rangle}, \tau_n, \overline{B}_n^{b-}(r_1, \overline{x}_{*[1,n]})\right) = \frac{1}{S_n(1) r_1^{n-1}} \int_{\overline{B}_n^b(r_1, \overline{x}_{*[1,n]})} \psi^{\langle 1,n \rangle}\left(x_{*0} + \lambda r_1, \overline{x}_{[1,n]}\right) ds_{\{n\}}, \quad (4.44)$$

$$\Sigma\left(\psi^{\langle 1,n \rangle}, \eta^{\langle n \rangle}, \tau_n, \overline{B}_n^{b+}(r_2, \overline{x}_{*[1,n]})\right) = -\frac{1}{S_n(1) r_2^{n-1}} \int_{\overline{B}_n^b(r_2, \overline{x}_{*[1,n]})} \psi^{\langle 1,n \rangle}\left(x_{*0} + \lambda r_2, \overline{x}_{[1,n]}\right) ds_{\{n\}} \quad (4.45)$$

for each $n \in \omega_2$. Ultimately, combination of (4.39) (4.40), (4.44), and (4.45) yields:

$$\Omega^{(3)}\left(\psi^{\langle 1,n \rangle}, \eta^{\langle n \rangle}, \tau_n, \overline{C}_n^{cc}(r_1, r_2, \overline{x}_{*[1,n]})\right)$$
$$= -\frac{1}{S_n(1) r_1^{n-1}} \int_{\overline{B}_n^b(r_1, \overline{x}_{*[1,n]})} \psi^{\langle 1,n \rangle}\left(x_{*0} + \lambda r_1, \overline{x}_{[1,n]}\right) ds_{\{n\}} \quad (4.46)$$
$$+ \frac{1}{S_n(1) r_2^{n-1}} \int_{\overline{B}_n^b(r_2, \overline{x}_{*[1,n]})} \psi^{\langle 1,n \rangle}\left(x_{*0} + \lambda r_2, \overline{x}_{[1,n]}\right) ds_{\{n\}}.$$

2) From the instance of (3.39) subject to (3.40) and (3.41) with $\psi^{\langle 1,1 \rangle}$ and $\eta^{\langle 1 \rangle}$ in place of $\phi^{\langle 1,1 \rangle}$ and $\psi^{\langle 1,1 \rangle}$ respectively, it follows that

$$\Omega^{(3)}\left(\psi^{\langle 1,1 \rangle},\eta^{\langle 1 \rangle},\tau_n,\overline{C}_1^{cc}(r_1,r_2,x_{*1})\right)$$
$$=\Omega^{(3)}\left(\psi^{\langle 1,1 \rangle},\eta^{\langle 1 \rangle},\tau_1,[x_{*1}-r_2,x_{*1}+r_2]\right)-\Omega^{(3)}\left(\psi^{\langle 1,1 \rangle},\eta^{\langle 1 \rangle},\tau_1,[x_{*1}+r_1,x_{*1}+r_1]\right) \quad (4.47)$$

subject to

$$\Omega^{(3)}\left(\psi^{\langle 1,1 \rangle},\eta^{\langle 1 \rangle},\tau_1,[x_{*1}+r_1,x_{*1}+r_1]\right)\equiv -\int_{x_{*1}-r_1}^{x_{*1}+r_1}\frac{\partial}{\partial x_1}L_1\left(\psi^{\langle 1,1 \rangle},\eta^{\langle 1 \rangle},\tau_1,x_1\right)dx_1$$
$$=L_1\left(\psi^{\langle 1,1 \rangle},\eta^{\langle 1 \rangle},\tau_1,x_{*1}-r_1\right)-L_1\left(\psi^{\langle 1,1 \rangle},\eta^{\langle 1 \rangle},\tau_1,x_{*1}+r_1\right), \quad (4.48)$$

$$\Omega^{(3)}\left(\psi^{\langle 1,1 \rangle},\eta^{\langle 1 \rangle},\tau_1,[x_{*1}-r_2,x_{*1}+r_2]\right)\equiv -\int_{x_{*1}-r_2}^{x_{*1}+r_2}\frac{\partial}{\partial x_1}L_1\left(\psi^{\langle 1,1 \rangle},\eta^{\langle 1 \rangle},\tau_1,x_1\right)dx_1$$
$$=L_1\left(\psi^{\langle 1,1 \rangle},\eta^{\langle 1 \rangle},\tau_1,x_{*1}-r_2\right)-L_1\left(\psi^{\langle 1,1 \rangle},\eta^{\langle 1 \rangle},\tau_1,x_{*1}+r_2\right). \quad (4.49)$$

At the same time, (4.30a) for $n=1$ reduces to

$$L_1\left(\psi^{\langle 1,1 \rangle},\eta^{\langle 1 \rangle},\tau_1,x_1\right)=\frac{a_1 x_1' \psi^{\langle 1,1 \rangle}\left(x_{*0}+\lambda|x_1'|,x_1\right)}{|x_1'|}=-\frac{x_1' \psi^{\langle 1,1 \rangle}\left(x_{*0}+\lambda|x_1'|,x_1\right)}{2|x_1'|} \quad (4.50)$$

subject to (4.14). Hence, (4.48) and (4.49) become

$$\Omega^{(3)}\left(\psi^{\langle 1,1 \rangle},\eta^{\langle 1 \rangle},\tau_1,[x_{*1}+r_1,x_{*1}+r_1]\right)$$
$$=-\frac{(-r_1)\psi^{\langle 1,1 \rangle}\left(x_{*0}+\lambda r_1,x_{*1}-r_1\right)}{2r_1}+\frac{r_1\psi^{\langle 1,1 \rangle}\left(x_{*0}+\lambda r_1,x_{*1}+r_1\right)}{2r_1} \quad (4.51)$$
$$=\frac{1}{2}\left[\psi^{\langle 1,1 \rangle}\left(x_{*0}+\lambda r_1,x_{*1}-r_1\right)+\psi^{\langle 1,1 \rangle}\left(x_{*0}+\lambda r_1,x_{*1}+r_1\right)\right],$$

$$\Omega^{(3)}\left(\psi^{\langle 1,1 \rangle},\eta^{\langle 1 \rangle},\tau_1,[x_{*1}+r_2,x_{*1}+r_2]\right)$$
$$=-\frac{(-r_2)\psi^{\langle 1,1 \rangle}\left(x_{*0}+\lambda r_2,x_{*1}-r_2\right)}{2r_2}+\frac{r_2\psi^{\langle 1,1 \rangle}\left(x_{*0}+\lambda r_2,x_{*1}+r_2\right)}{2r_2} \quad (4.52)$$
$$=\frac{1}{2}\left[\psi^{\langle 1,1 \rangle}\left(x_{*0}+\lambda r_2,x_{*1}-r_2\right)+\psi^{\langle 1,1 \rangle}\left(x_{*0}+\lambda r_2,x_{*1}+r_2\right)\right],$$

Consequently, (4.47) becomes

$$\Omega^{(3)}\left(\psi^{\langle 1,1 \rangle},\eta^{\langle 1 \rangle},\tau_1,\overline{C}_1^{cc}(r_1,r_2,x_{*1})\right)$$
$$=-\frac{1}{2}\left[\psi^{\langle 1,1 \rangle}\left(x_{*0}+\lambda r_1,x_{*1}-r_1\right)+\psi^{\langle 1,1 \rangle}\left(x_{*0}+\lambda r_1,x_{*1}+r_1\right)\right] \quad (4.53)$$
$$+\frac{1}{2}\left[\psi^{\langle 1,1 \rangle}\left(x_{*0}+\lambda r_2,x_{*1}-r_2\right)+\psi^{\langle 1,1 \rangle}\left(x_{*0}+\lambda r_2,x_{*1}+r_2\right)\right],$$

which is analogous to (4.46).

3) From the instance of (3.31) subject to (3.32) (see particularly (3.36) subject to (3.37) and (3.38)) with $\eta^{\langle n\rangle}$ in place of $\phi^{\langle 1,n\rangle}$, it follows that for each $n \in \omega_2$:

$$\Omega^{(3)}\left(\eta^{\langle n\rangle},\psi^{\langle 1,n\rangle},\tau_n,\overline{C}_n^{cc}\left(r_1,r_2,\overline{x}_{*[1,n]}\right)\right) = -\Sigma\left(\eta^{\langle n\rangle},\psi^{\langle 1,n\rangle},\tau_n,\overline{C}_n^b\left(r_1,r_2,\overline{x}_{*[1,n]}\right)\right) \quad (4.54)$$

subject to

$$\Sigma\left(\eta^{\langle n\rangle},\psi^{\langle 1,n\rangle},\tau_n,\overline{C}_n^b\left(r_1,r_2,\overline{x}_{*[1,n]}\right)\right)$$
$$= \Sigma\left(\eta^{\langle n\rangle},\psi^{\langle 1,n\rangle},\tau_n,\overline{B}_n^{b-}\left(r_1,\overline{x}_{*[1,n]}\right)\right) + \Sigma\left(\eta^{\langle n\rangle},\psi^{\langle 1,n\rangle},\tau_n,\overline{B}_n^{b+}\left(r_2,\overline{x}_{*[1,n]}\right)\right), \quad (4.55)$$

where

$$\Sigma\left(\eta^{\langle n\rangle},\psi^{\langle 1,n\rangle},\tau_n,\overline{B}_n^{b-}\left(r_1,\overline{x}_{*[1,n]}\right)\right)$$
$$= -\frac{1}{r_1} \int_{\overline{B}_n^{b-}(r_1,\overline{x}_{*[1,n]})} \sum_{i=1}^{n} L_i\left(\eta^{\langle n\rangle},\psi^{\langle 1,n\rangle},\tau_n,\overline{x}_{[1,n]}\right) x_i' \, ds_{\{n\}}, \quad (4.56)$$

$$\Sigma\left(\eta^{\langle n\rangle},\psi^{\langle 1,n\rangle},\tau_n,\overline{B}_n^{b+}\left(r_2,\overline{x}_{*[1,n]}\right)\right)$$
$$= \frac{1}{r_2} \int_{\overline{B}_n^b(r_2,\overline{x}_{*[1,n]})} \sum_{i=1}^{n} L_i\left(\eta^{\langle n\rangle},\psi^{\langle 1,n\rangle},\tau_n,\overline{x}_{[1,n]}\right) x_i' \, ds_{\{n\}}, \quad (4.57)$$

which are the pertinent instances of (3.37) and (3.38) respectively. By (4.27), it follows that

$$\sum_{i=1}^{n} L_i\left(\eta^{\langle n\rangle},\psi^{\langle 1,n\rangle},\tau_n,\overline{x}_{[1,n]}\right) x_i' = x_{[1,n]}' \eta^{\langle n\rangle}\left(x_{[1,n]}'\right) \Lambda\left(\psi^{\langle 1,n\rangle},\tau_n,\overline{x}_{[1,n]}\right), \quad (4.58)$$

where

$$\Lambda\left(\psi^{\langle 1,n\rangle},\tau_n,\overline{x}_{[1,n]}\right)$$
$$\equiv \sum_{i=1}^{n} \frac{x_i'}{x_{[1,n]}'} \psi_{x_i}^{\langle 1,n\rangle}\left(x_{*0} + \lambda x_{[1,n]}',\overline{x}_{[1,n]}\right) - \lambda \psi_{x_0}^{\langle 1,n\rangle}\left(x_{*0} + \lambda x_{[1,n]}',\overline{x}_{[1,n]}\right) \quad (4.59)$$

and, by (4.9) and (4.11),

$$\eta^{\langle n\rangle}\left(x_{[1,n]}'\right) = \begin{cases} a_n x_{[1,n]}'^{2-n} & \text{if } n \neq 2 \quad (a) \\ -a_2 \ln x_{[1,2]}' & \text{if } n = 2 \quad (b) \end{cases} \quad (4.60)$$

subject to (4.10) and (4.12). Hence, equations (4.56) and (4.57) become

$$\Sigma\left(\eta^{\langle n\rangle},\psi^{\langle 1,n\rangle},\tau_n,\overline{B}_n^{b-}\left(r_1,\overline{x}_{*[1,n]}\right)\right) = -\eta^{\langle n\rangle}(r_1) \int_{\overline{B}_n^b(r_1,\overline{x}_{*[1,n]})} \Lambda\left(\psi^{\langle 1,n\rangle},\tau_n,\overline{x}_{[1,n]}\right) ds_{\{n\}}, \quad (4.61)$$

$$\Sigma\left(\eta^{\langle n\rangle},\psi^{\langle 1,n\rangle},\tau_n,\overline{B}_n^{b+}\left(r_2,\overline{x}_{*[1,n]}\right)\right) = \eta^{\langle n\rangle}(r_2) \int_{\overline{B}_n^b(r_2,\overline{x}_{*[1,n]})} \Lambda\left(\psi^{\langle 1,n\rangle},\tau_n,\overline{x}_{[1,n]}\right) ds_{\{n\}} \quad (4.62)$$

for each $n \in \omega_2$. Ultimately, combination of (4.54), (4.55), (4.61), and (4.62) yields:

$$\Omega^{(3)}\left(\eta^{\langle n \rangle}, \psi^{\langle 1,n \rangle}, \tau_n, \overline{C}_n^{cc}\left(r_1, r_2, \overline{x}_{*[1,n]}\right)\right)$$
$$= \eta^{\langle n \rangle}(r_1) \int_{\overline{B}_n^b(r_1, \overline{x}_{*[1,n]})} \Lambda\left(\psi^{\langle 1,n \rangle}, \tau_n, \overline{x}_{[1,n]}\right) ds_{\{n\}} - \eta^{\langle n \rangle}(r_2) \int_{\overline{B}_n^b(r_2, \overline{x}_{*[1,n]})} \Lambda\left(\psi^{\langle 1,n \rangle}, \tau_n, \overline{x}_{[1,n]}\right) ds_{\{n\}} \tag{4.63}$$

for each $n \in \omega_2$.

4) From the instance of (3.39) subject to (3.40) and (3.41), with $\eta^{\langle 1 \rangle}$ in place of $\phi^{\langle 1,1 \rangle}$, it follows that:

$$\Omega^{(3)}\left(\eta^{\langle 1 \rangle}, \psi^{\langle 1,1 \rangle}, \tau_1, \overline{C}_1^{cc}(r_1, r_2, x_{*1})\right)$$
$$= \Omega^{(3)}\left(\eta^{\langle 1 \rangle}, \psi^{\langle 1,1 \rangle}, \tau_1, [x_{*1} - r_2, x_{*1} + r_2]\right) \tag{4.64}$$
$$- \Omega^{(3)}\left(\eta^{\langle 1 \rangle}, \psi^{\langle 1,1 \rangle}, \tau_1, [x_{*1} + r_1, x_{*1} + r_1]\right),$$

where

$$\Omega^{(3)}\left(\eta^{\langle 1 \rangle}, \psi^{\langle 1,1 \rangle}, \tau_1, [x_{*1} - r_1, x_{*1} + r_1]\right)$$
$$= L_1\left(\eta^{\langle 1 \rangle}, \psi^{\langle 1,1 \rangle}, \tau_1, x_{*1} - r_1\right) - L_1\left(\eta^{\langle 1 \rangle}, \psi^{\langle 1,1 \rangle}, \tau_1, x_{*1} + r_1\right), \tag{4.65}$$

$$\Omega^{(3)}\left(\eta^{\langle 1 \rangle}, \psi^{\langle 1,1 \rangle}, \tau_1, [x_{*1} - r_2, x_{*1} + r_2]\right)$$
$$= L_1\left(\eta^{\langle 1 \rangle}, \psi^{\langle 1,1 \rangle}, \tau_1, x_{*1} - r_2\right) - L_1\left(\eta^{\langle 1 \rangle}, \psi^{\langle 1,1 \rangle}, \tau_1, x_{*1} + r_2\right). \tag{4.66}$$

At the same time, by (4.13), the equation (4.27) at $n=1$ reduces to

$$L_1\left(\eta^{\langle 1 \rangle}, \psi^{\langle 1,1 \rangle}, \tau_1, x_1\right) = -\frac{1}{2}|x_1'| \Psi_1\left(x_{*0} + \lambda|x_1'|, x_{1,}\right), \tag{4.67}$$

where

$$\Psi_1\left(x_{*0} + \lambda|x_1'|, x_{1,}\right) \equiv \psi_{x_1}^{\langle 1,1 \rangle}\left(x_{*0} + \lambda|x_1'|, x_{1,}\right) - \frac{\lambda x_1'}{|x_1'|} \psi_{x_0}^{\langle 1,1 \rangle}\left(x_{*0} + \lambda|x_1'|, x_{1,}\right). \tag{4.68}$$

Hence, (4.65) and (4.66) become

$$\Omega^{(3)}\left(\eta^{\langle 1 \rangle}, \psi^{\langle 1,1 \rangle}, \tau_1, [x_{*1} - r_1, x_{*1} + r_1]\right)$$
$$= -\frac{r_1}{2}\left[\Psi_1(x_{*0} + \lambda r_1, x_{*1} - r_1) - \Psi_1(x_{*0} + \lambda r_1, x_{*1} + r_1)\right], \tag{4.69}$$

$$\Omega^{(3)}\left(\eta^{\langle 1 \rangle}, \psi^{\langle 1,1 \rangle}, \tau_1, [x_{*1} - r_2, x_{*1} + r_2]\right)$$
$$= -\frac{r_2}{2}\left[\Psi_1(x_{*0} + \lambda r_2, x_{*1} - r_2) - \Psi_1(x_{*0} + \lambda r_2, x_{*1} + r_2)\right]. \tag{4.70}$$

Consequently, (4.64) becomes

$$\Omega^{(3)}\left(\eta^{\langle 1 \rangle}, \psi^{\langle 1,1 \rangle}, \tau_1, \overline{C}_{\{1\}}^{cc}(r_1, r_2, x_{*1})\right)$$
$$= \frac{r_1}{2}\left[\Psi_1(x_{*0} + \lambda r_1, x_{*1} - r_1) - \Psi_1(x_{*0} + \lambda r_1, x_{*1} + r_1)\right] \tag{4.71}$$
$$- \frac{r_2}{2}\left[\Psi_1(x_{*0} + \lambda r_2, x_{*1} - r_2) - \Psi_1(x_{*0} + \lambda r_2, x_{*1} + r_2)\right],$$

which is analogous to (4.63).

5) For each $n \in \omega_1$ and each $\alpha \in \{0,1\}$, by (4.24) and (4.25), the pertinent instances of (3.29) with $\eta^{\langle n \rangle}$ in place of $\phi^{\langle 1,n \rangle}$ become:

$$\Omega^{(0)}\left(\eta^{\langle n \rangle}, \psi^{\langle 1,n \rangle}, \tau_n, \overline{C}_n^{cc}(r_1, r_2, \overline{x}_{*[1,n]})\right)$$
$$= -\int_{\overline{C}_n^{cc}(r_1, r_2, \overline{x}_{*[1,n]})} \eta^{\langle n \rangle}(x'_{[1,n]}) f^{\langle 1,n \rangle}(x_{*0} + \lambda x'_{[1,n]}, \overline{x}_{[1,n]}) dv_{\{n\}} \tag{4.72}$$

subject to (4.60), and

$$\Omega^{(1)}\left(\eta^{\langle n \rangle}, \psi^{\langle 1,n \rangle}, \tau_n, \overline{C}_n^{cc}(r_1, r_2, \overline{x}_{*[1,n]})\right)$$
$$= \lambda(n-3) \int_{\overline{C}_n^{cc}(r_1, r_2, \overline{x}_{*[1,n]})} \zeta^{\langle n \rangle}(x'_{[1,n]}) \left[\frac{\partial \psi^{\langle 1,n \rangle}(x_0, \overline{x}_{[1,n]})}{\partial x_0}\right]_{x_0 = x_{*0} + \lambda x'_{[1,n]}} dv_{\{n\}} \tag{4.73}$$

subject to

$$\zeta^{\langle n \rangle}(x'_{[1,n]}) \stackrel{\sim}{=} \begin{cases} a_n x'^{1-n}_{[1,n]} & \text{if } n \neq 2 \text{ (a)} \\ -\dfrac{a_2(\ln x'_{[1,2]} + 2)}{x'_{[1,2]}} & \text{if } n = 2 \text{ (b)} \end{cases}. \tag{4.74}$$

Hence, specifically, for convenience in further references and uses,

$$\Omega^{(0)}\left(\eta^{\langle n \rangle}, \psi^{\langle 1,n \rangle}, \tau_n, \overline{C}_n^{cc}(r_1, r_2, \overline{x}_{*[1,n]})\right)$$
$$= -a_n \int_{\overline{C}_n^{cc}(r_1, r_2, \overline{x}_{*[1,n]})} \frac{1}{x'^{n-2}_{[1,n]}} f^{\langle 1,n \rangle}(x_{*0} + \lambda x'_{[1,n]}, \overline{x}_{[1,n]}) dv_{\{n\}} \text{ for } n \neq 2, \tag{4.72a}$$

$$\Omega^{(0)}\left(\eta^{\langle 2 \rangle}, \psi^{\langle 1,2 \rangle}, \tau_2, \overline{C}_2^{cc}(r_1, r_2, \overline{x}_{*[1,2]})\right)$$
$$= a_2 \int_{\overline{C}_2^{cc}(r_1, r_2, \overline{x}_{*[1,2]})} (\ln x'_{[1,2]}) f^{\langle 1,2 \rangle}(x_{*0} + \lambda x'_{[1,2]}, \overline{x}_{[1,2]}) dv_{\{2\}} \text{ for } n = 2, \tag{4.72b}$$

$$\Omega^{(1)}\left(\eta^{\langle n \rangle}, \psi^{\langle 1,n \rangle}, \tau_n, \overline{C}_n^{cc}(r_1, r_2, \overline{x}_{*[1,n]})\right)$$
$$= \lambda(n-3)a_n \int_{\overline{C}_n^{cc}(r_1, r_2, \overline{x}_{*[1,n]})} \frac{1}{x'^{n-1}_{[1,n]}} \psi^{\langle 1,n \rangle}_{x_0}(x_{*0} + \lambda x'_{[1,n]}, \overline{x}_{[1,n]}) dv_{\{n\}} \text{ for } n \neq 2, \tag{4.73a}$$

$$\Omega^{(1)}\left(\eta^{\langle 2 \rangle}, \psi^{\langle 1,2 \rangle}, \tau_2, \overline{C}_2^{cc}(r_1, r_2, \overline{x}_{*[1,2]})\right)$$
$$= \lambda a_2 \int_{\overline{C}_2^{cc}(r_1, r_2, \overline{x}_{*[1,2]})} \frac{\ln x'_{[1,2]} + 2}{x'_{[1,2]}} \psi^{\langle 1,2 \rangle}_{x_0}(x_{*0} + \lambda x'_{[1,2]}, \overline{x}_{[1,2]}) dv_{\{2\}}. \tag{4.73b}$$

**Proof:** The proof of the theorem is included in its statement.•

**Comment 4.3.** After establishing the integro-differential tautology (4.37), the definition (4.21) of the functional form $f^{\langle 1,n \rangle}(x_0, \overline{x}_{[1,n]})$, which occurs in the integrand of (4.72), can be turned into the *inhomogeneous differential equation*

$$D_{[0,n]}\psi^{\langle 1,n\rangle}\left(x_0,\overline{x}_{[1,n]}\right)=-f^{\langle 1,n\rangle}\left(x_0,\overline{x}_{[1,n]}\right), \tag{4.75}$$

in which $\psi^{\langle 1,n\rangle}\left(x_0,\overline{x}_{[1,n]}\right)$ is an unknown functional form, while $f^{\langle 1,n\rangle}\left(x_0,\overline{x}_{[1,n]}\right)$ is a known (given) functional form that has the appropriate properties so as not to alter the *syntactic form* of tautology (4.37). However, from the standpoint of *semantic analysis,* In this case, tautology (4.37) subject to (4.75) should be regarded as an *inhomogeneous integro-differential equation for* $\psi^{\langle 1,n\rangle}\left(x_0,\overline{x}_{[1,n]}\right)$ in the domain $\overline{C}_{\{n\}}^{cc}\left(r_1,r_2,\overline{x}_{*[1,n]}\right)$, which is implied by the inhomogeneous differential equation (4.75).•

**Comment 4.4.** At *n*=3, relation (4.73a) turns into the identity:

$$\Omega^{(1)}\left(\eta^{\langle 3\rangle},\psi^{\langle 1,3\rangle},\tau_3,\overline{C}_3^{cc}\left(r_1,r_2,\overline{x}_{*[1,3]}\right)\right)=0, \tag{4.73!}$$

which immediately follows from identity (4.25!), given in Comment 4.1, and which is, like the latter identity, the *exclusive fundamental property of* $\overline{E}_{\{3\}}$, underlying both *the principle of non-dispersed causality* and *the principle of special relativity – the principles that, owing to (4.73!), prove to be possible only in a three-dimensional affine Euclidean space* (to be explicated in due course).•

**Theorem 4.4:** *The first special integro-differential tautology for a single biased functional form.* Given $n\in\omega_1$,

$$\sum_{\alpha=0}^{3}\Omega^{(\alpha)}\left(\eta^{\langle n\rangle},\psi^{\langle 1,n\rangle},\tau_n,\overline{C}_n^{cc}\left(r_1,r_2,\overline{x}_{*[1,n]}\right)\right)=0, \tag{4.76}$$

where the summands with $\alpha\in\{0,1,3\}$ have been specified in Theorem 4.3, and

$$\begin{aligned}&\Omega^{(2)}\left(\eta^{\langle n\rangle},\psi^{\langle 1,n\rangle},\tau_n,\overline{C}_n^{cc}\left(r_1,r_2,\overline{x}_{*[1,n]}\right)\right)\\&=\int_{\overline{C}_n^{cc}\left(r_1,r_2,\overline{x}_{*[1,n]}\right)}\sum_{i=1}^{n}\frac{\partial\eta^{\langle n\rangle}\left(x'_{[1,n]}\right)}{\partial x_i}\frac{\partial\psi^{\langle 1,n\rangle}\left(x_{*0}+\lambda x'_{[1,n]},\overline{x}_{[1,n]}\right)}{\partial x_i}dv_{\{n\}}\\&=-\frac{1}{S_n(1)}\int_{\overline{C}_n^{cc}\left(r_1,r_2,\overline{x}_{*[1,n]}\right)}\frac{1}{x'^n_{[1,n]}}\sum_{i=1}^{n}x'_i\frac{\partial\psi^{\langle 1,n\rangle}\left(x_{*0}+\lambda x'_{[1,n]},\overline{x}_{[1,n]}\right)}{\partial x_i}dv_{\{n\}},\end{aligned} \tag{4.77}$$

while the summand with $\alpha=2$ is given by the pertinent instance of (3.29) with $\eta^{\langle n\rangle}$ in place of $\phi^{\langle 1,n\rangle}$ subject to (4.31), so that

$$\Omega^{(2)}\left(\eta^{\langle n \rangle}, \psi^{\langle 1,n \rangle}, \tau_n, \overline{C}_n^{cc}\left(r_1, r_2, \overline{x}_{*[1,n]}\right)\right)$$

$$= \int_{\overline{C}_n^{cc}\left(r_1, r_2, \overline{x}_{*[1,n]}\right)} \sum_{i=1}^{n} \frac{\partial \eta^{\langle n \rangle}\left(x'_{[1,n]}\right)}{\partial x_i} \cdot \frac{\partial \psi^{\langle 1,n \rangle}\left(x_{*0} + \lambda x'_{[1,n]}, \overline{x}_{[1,n]}\right)}{\partial x_i} dv_{\{n\}} \quad (4.78)$$

$$= -\frac{1}{S_n(1)} \int_{\overline{C}_n^{cc}\left(r_1, r_2, \overline{x}_{*[1,n]}\right)} \frac{1}{x'^{n}_{[1,n]}} \sum_{i=1}^{n} x'_i \frac{\partial \psi^{\langle 1,n \rangle}\left(x_{*0} + \lambda x'_{[1,n]}, \overline{x}_{[1,n]}\right)}{\partial x_i} dv_{\{n\}},$$

by (4.32) subject t0 (4.10) and (4.12).•

## 5. General properties of improper *n*-fold volume integrals

**Definition 5.1.** 1) For each $n \in \omega_1$, the integrals involved into either of the special integro-differential formulas (4.38) and (4.75) are of two kinds, namely *n*-fold volume integrals over the closed *n*-dimensional spherical layer $\overline{C}_n^{cc}(r_1, r_2, \overline{x}_{*[1,n]})$ and $(n-1)$-fold surface ones over the *n*-dimensional spherical surfaces $\overline{B}_n^b(r_1, \overline{x}_{*[1,n]})$ and $\overline{B}_n^b(r_2, \overline{x}_{*[1,n]})$ forming the disconnected boundary $\overline{C}_n^b(r_1, r_2, \overline{x}_{*[1,n]})$ of $\overline{C}_n^{cc}(r_1, r_2, \overline{x}_{*[1,n]})$. At $n=1$ $\overline{B}_1^b(r_1, x_{*1}) = \langle x_{*1} - r_1, x_{*1} + r_1 \rangle$ and $\overline{B}_1^b(r_2, x_{*1}) = \langle x_{*1} - r_2, x_{*1} + r_2 \rangle$, so that by a 0-fold surface integral over $\overline{B}_1^b(r_1, x_{*1})$, e.g., we understand the difference of values of the pertinent antiderivative functional form at the boundary points $x_{*1} + r_1$ and $x_{*1} - r_1$ in this order.

2) If '$\Phi^{\langle n \rangle}(\overline{x}_{[1,n]})$' is a placeholder for the integrand of an *n*-fold volume integral over $\overline{C}_n^{cc}(r_1, r_2, \overline{x}_{*[1,n]})$, i.e. if $\Phi^{\langle n \rangle}(\overline{x}_{[1,n]})$ *is* such an integrand, then

$$\int_{\overline{B}_n^c(r_2, \overline{x}_{*[1,n]})} \Phi^{\langle n \rangle}(\overline{x}_{[1,n]}) dv_{\{n\}} \stackrel{\triangle}{=} \lim_{r_1 \to 0} \int_{\overline{C}_n^{cc}(r_1, r_2, \overline{x}_{*[1,n]})} \Phi^{\langle n \rangle}(\overline{x}_{[1,n]}) dv_{\{n\}}, \tag{5.1}$$

$$\int_{\overline{E}_n - \overline{B}_n^o(r_1, \overline{x}_{*[1,n]})} \Phi^{\langle n \rangle}(\overline{x}_{[1,n]}) dv_{\{n\}} \stackrel{\triangle}{=} \lim_{r_2 \to \infty} \int_{\overline{C}_n^{cc}(r_1, r_2, \overline{x}_{*[1,n]})} \Phi^{\langle n \rangle}(\overline{x}_{[1,n]}) dv_{\{n\}}, \tag{5.2}$$

$$\int_{\overline{E}_n} \Phi^{\langle n \rangle}(\overline{x}_{[1,n]}) dv_{\{n\}} \stackrel{\triangle}{=} \lim_{r_2 \to \infty} \int_{\overline{B}_n^c(r_2, \overline{x}_{*[1,n]})} \Phi^{\langle n \rangle}(\overline{x}_{[1,n]}) dv_{\{n\}}$$
$$= \lim_{r_1 \to 0} \int_{\overline{E}_\{ - \overline{B}_n^o(r_1, \overline{x}_{*[1,n]})} \Phi^{\langle n \rangle}(\overline{x}_{[1,n]}) dv_{\{n\}}. \tag{5.3}$$

The definiendum of definition (5.1), (5.2), or (5.3) is called an *improper n-fold volume integral of first, second,* or *third kind* respectively provided that *its definiens has a denotatum*. If an improper integral [of any kind] *exists*, i.e. *has a denotatum*, then it is said to *converge* and vice versa.

3) An improper integral is said to converge absolutely or to be an absolutely converging if and only if converges the improper integral that is obtained by replacing the integrand of the former with its absolute value. Thus, any one of the improper integrals defined by (5.1)–(5.3) converges absolutely if the corresponding one of the following limits exists:

$$\int_{\overline{B}_n^c(r_2, \overline{x}_{*[1,n]})} \left| \Phi^{\langle n \rangle}(\overline{x}_{[1,n]}) \right| dv_{\{n\}} \stackrel{\triangle}{=} \lim_{r_1 \to 0} \int_{\overline{C}_n^{cc}(r_1, r_2, \overline{x}_{*[1,n]})} \left| \Phi^{\langle n \rangle}(\overline{x}_{[1,n]}) \right| dv_{\{n\}}, \tag{5.4}$$

$$\int_{\overline{E}_n - \overline{B}_n^o(r_1, \overline{x}_{*[1,n]})} \left| \Phi^{\langle n \rangle}(\overline{x}_{[1,n]}) \right| dv_{\{n\}} \stackrel{\triangle}{=} \lim_{r_2 \to \infty} \int_{\overline{C}_n^{cc}(r_1, r_2, \overline{x}_{*[1,n]})} \left| \Phi^{\langle n \rangle}(\overline{x}_{[1,n]}) \right| dv_{\{n\}}, \tag{5.5}$$

$$\int_{\overline{E}_n} \left|\Phi^{\langle n\rangle}(\overline{x}_{[1,n]})\right| dv_{\{n\}} \equiv \lim_{r_2 \to \infty} \int_{\overline{B}_n^c(r_2,\overline{x}_{*[1,n]})} \left|\Phi^{\langle n\rangle}(\overline{x}_{[1,n]})\right| dv_{\{n\}}$$
$$= \lim_{r \to 0} \int_{\overline{E}_n - \overline{B}_n^o(r_1,\overline{x}_{*[1,n]})} \left|\Phi^{\langle n\rangle}(\overline{x}_{[1,n]})\right| dv_{\{n\}}. \quad (5.6)$$

**Theorem 5.1:** *A sufficient condition for convergence of an improper integral.* An improper integral of any one of the three kinds converges if it converges absolutely.

**Proof**: Let

$$\int_{\overline{X}_n^c} \Phi^{\langle n\rangle}(\overline{x}_{[1,n]}) dv_{\{n\}}$$

be an improper integral over any closed domain $\overline{X}_n^c \subseteq \overline{E}_n$ (as $\overline{B}_n^c(r_1,\overline{x}_{*[1,n]})$, $\overline{E}_n - \overline{B}_n^o(r_1,\overline{x}_{*[1,n]})$, or $\overline{E}_n$). By the hypothesis (antecedent) of the theorem, the corresponding improper integral

$$\int_{\overline{X}_n^c} \left|\Phi^{\langle n\rangle}(\overline{x}_{[1,n]})\right| dv_{\{n\}}$$

exists (converges). Since

$$\Phi^{\langle n\rangle}(\overline{x}_{[1,n]}) = \left|\Phi^{\langle n\rangle}(\overline{x}_{[1,n]})\right| - \left[\left|\Phi^{\langle n\rangle}(\overline{x}_{[1,n]})\right| - \Phi^{\langle n\rangle}(\overline{x}_{[1,n]})\right],$$

therefore

$$\int_{\overline{X}_{\{n\}}^c} \Phi^{\langle n\rangle}(\overline{x}_{[1,n]}) dv_{\{n\}} = \int_{\overline{X}_{\{n\}}^c} \left|\Phi^{\langle n\rangle}(\overline{x}_{[1,n]})\right| dv_{\{n\}}$$
$$- \int_{\overline{X}_{\{n\}}^c} \left[\left|\Phi^{\langle n\rangle}(\overline{x}_{[1,n]})\right| - \Phi^{\langle n\rangle}(\overline{x}_{[1,n]})\right] dv_{\{n\}}. \quad (5.7)$$

In this case,

$$0 \leq \int_{\overline{X}_{\{n\}}^c} \left[\left|\Phi^{\langle n\rangle}(\overline{x}_{[1,n]})\right| - \Phi^{\langle n\rangle}(\overline{x}_{[1,n]})\right] dv_{\{n\}} \leq 2 \int_{\overline{X}_{\{n\}}^c} \left|\Phi^{\langle n\rangle}(\overline{x}_{[1,n]})\right| dv_{\{n\}},$$

because

$$\left|\Phi^{\langle n\rangle}(\overline{x}_{[1,n]})\right| - \Phi^{\langle n\rangle}(\overline{x}_{[1,n]}) = \begin{cases} 0 & \text{if } \Phi^{\langle n\rangle}(\overline{x}_{[1,n]}) > 0 \\ 2\left|\Phi^{\langle n\rangle}(\overline{x}_{[1,n]})\right| & \text{if } \Phi^{\langle n\rangle}(\overline{x}_{[1,n]}) \leq 0 \end{cases}.$$

Thus, both integrals on the right-hand side of the tautology (5.14) converge and hence the integral on the left-hand side of the tautology (5.14) converges as well. QED.•

**Theorem 5.2:** *A sufficient condition for convergence of an improper integral of first kind.* Given $n \in \omega_1$, a volume improper integral of first kind, defined by (5.1), converges if there exist three real numbers: $r_0 \in (0,r_2)$, $p_n < n$, and $k_n > 0$ such that for each $\overline{x}_{[1,n]} \in \overline{B}_n^c(r_0,\overline{x}_{*[1,n]})$:

$$\left|\Phi^{\langle n\rangle}\!\left(\overline{x}_{[1,n]}\right)\right| \le \frac{k_n}{\left|\overline{x}_{[1,n]} - \overline{x}_{*[1,n]}\right|^{p_n}} \quad \text{for each } \overline{x}_{[1,n]} \in \overline{B}_n^c\!\left(r_0, \overline{x}_{*[1,n]}\right). \tag{5.8}$$

**Proof:** 1) The definiendum of (5.1) can be represented as:

$$\int\limits_{\overline{B}_n^c(r_2,\overline{x}_{*[1,n]})} \Phi^{\langle n\rangle}\!\left(\overline{x}_{[1,n]}\right) dv_{\{n\}} = \int\limits_{\overline{B}_n^c(r_0,\overline{x}_{*[1,n]})} \Phi^{\langle n\rangle}\!\left(\overline{x}_{[1,n]}\right) dv_{\{n\}} + \int\limits_{\overline{C}_n^{cc}(r_0,r_2,\overline{x}_{*[1,n]})} \Phi^{\langle n\rangle}\!\left(\overline{x}_{[1,n]}\right) dv_{\{n\}}, \tag{5.9}$$

where the second integral on the right-hand side is an ordinary Riemann integral, while the first integral is an improper integral of first kind, whose value can be estimated with the help of (5.8) thus:

$$\begin{aligned}
\left|\int\limits_{\overline{B}_n^c(r_0,\overline{x}_{*[1,n]})} \Phi^{\langle n\rangle}\!\left(\overline{x}_{[1,n]}\right) dv_{\{n\}}\right| &\le \int\limits_{\overline{B}_n^c(r_0,\overline{x}_{*[1,n]})} \left|\Phi^{\langle n\rangle}\!\left(\overline{x}_{[1,n]}\right)\right| dv_{\{n\}} \\
&\le \int\limits_0^{r_0} \frac{k_n dV_n(r)}{r^{p_n}} = k_n S_n(1) \int\limits_0^{r_0} r^{n-p_n-1} dr = \frac{k_n S_n(1)}{n-p_n} r_0^{n-p_n},
\end{aligned} \tag{5.10}$$

because

$$dV_n(r) = S_n(r)\,dr = r^{n-1} S_n(1)\,dr. \tag{5.11}\bullet$$

**Theorem 5.3:** *A sufficient condition for convergence of an improper integral of second kind.* Given $n \in \omega_1$, the volume improper integral of second kind, defined by (5.2), converges if there exist three real numbers: $r_* \in (r_1, \infty)$, $q_n > n$, and $K_n > 0$ such that

$$\left|\Phi_{\langle n\rangle}\!\left(\overline{x}_{[1,n]}\right)\right| \le \frac{K_n}{\left|\overline{x}_{[1,n]} - \overline{x}_{*[1,n]}\right|^{q_n}} \quad \text{for each } \overline{x}_{[1,n]} \in \overline{E}_n - \overline{B}_n^o\!\left(r_*, \overline{x}_{*[1,n]}\right). \tag{5.12}$$

**Proof:** The definiendum of (5.2) can be represented as:

$$\begin{aligned}
&\int\limits_{\overline{E}_n - \overline{B}_n^o(r_1,\overline{x}_{*[1,n]})} \Phi^{\langle n\rangle}\!\left(\overline{x}_{[1,n]}\right) dv_{\{n\}} \\
&= \int\limits_{\overline{C}_n^{cc}(r_1,r_*,\overline{x}_{*[1,n]})} \Phi^{\langle n\rangle}\!\left(\overline{x}_{[1,n]}\right) dv_{\{n\}} + \int\limits_{\overline{E}_n - \overline{B}_n^o(r_*,\overline{x}_{*[1,n]})} \Phi^{\langle n\rangle}\!\left(\overline{x}_{[1,n]}\right) dv_{\{n\}},
\end{aligned} \tag{5.13}$$

where the first integral on the right-hand side is an ordinary Riemann integral, while the second integral is an improper integral of second kind, i.e.

$$\int\limits_{\overline{E}_n - \overline{B}_n^o(r_*,\overline{x}_{*[1,n]})} \Phi^{\langle n\rangle}\!\left(\overline{x}_{[1,n]}\right) dv_{\{n\}} = \lim_{r_2 \to \infty} \int\limits_{\overline{C}_n^{cc}(r_*,r_2,\overline{x}_{*[1,n]})} \Phi^{\langle n\rangle}\!\left(\overline{x}_{[1,n]}\right) dv_{\{n\}}. \tag{5.14}$$

whose value can be estimated with the help of (5.12) thus:

$$\left| \int_{\overline{C}_n^{cc}(r_*,r_2,\overline{x}_{*[1,n]})} \Phi^{\langle n\rangle}(\overline{x}_{[1,n]}) dv_{\{n\}} \right| \leq \int_{\overline{C}_n^{cc}(r_*,r_2,\overline{x}_{*[1,n]})} \left|\Phi^{\langle n\rangle}(\overline{x}_{[1,n]})\right| dv_{\{n\}} \leq \int_{r_*}^{R} \frac{K_n dV_n(r)}{r^{q_n}}$$

$$= K_n S_n(1) \int_{r_*}^{r_2} r^{n-q_n-1} dr = \frac{K_n S_n(1)}{q_n - n} \left( \frac{1}{r_*^{q_n-n}} - \frac{1}{r_2^{q_n-n}} \right) \leq \frac{K_n S_n(1)}{(q_n - n) r_*^{q_n-n}}, \quad (5.15)$$

whence

$$\left| \lim_{r_2 \to \infty} \int_{\overline{C}_n^{cc}(r_*,r_2,\overline{x}_{*[1,n]})} \Phi^{\langle n\rangle}(\overline{x}_{[1,n]}) dv_{\{n\}} \right| \leq \frac{K_n S_n(1)}{(q_n - n) r_*^{q_n-n}}. \quad (5.16)\bullet$$

**Theorem 5.4:** *A sufficient condition for convergence of an improper integral of third kind.* Given $n \in \omega_1$, the volume improper integral of third kind, defined by (5.3) converges if hypotheses of Theorems 5.2 and 5.3 hold simultaneously subject to $r_* > r_0$.

**Proof:** The difeniendum of (5.3) can be represented as:

$$\int_{\overline{E}_n} \Phi^{\langle n\rangle}(\overline{x}_{[1,n]}) dv_{\{n\}} = \int_{\overline{B}_n^c(r_0,\overline{x}_{*[1,n]})} \Phi^{\langle n\rangle}(\overline{x}_{[1,n]}) dv_{\{n\}} + \int_{\overline{C}_n^{cc}(r_0,r_*,\overline{x}_{*[1,n]})} \Phi^{\langle n\rangle}(\overline{x}_{[1,n]}) dv_{\{n\}}$$
$$+ \int_{\overline{E}_n - \overline{B}_n^o(r_*,\overline{x}_{*[1,n]})} \Phi^{\langle n\rangle}(\overline{x}_{[1,n]}) dv_{\{n\}}, \quad (5.17)$$

where the first integral on the right-hand side is a converging improper intrgral of first kind (by Theorem 5.2), the second integral is an ordinary Riemann integral, and the third integral is a converging improper integral of second kind (by Theorem 5.3).$\bullet$

# 6. The special second integro-differential formula for a biased functional form in a closed *n*-dimensional sphere

**Theorem 6.1.** 1) Given $\lambda \in \{-1,1\}$, given $n \in \omega_1$, as $r_1 \to 0$ formula (4.38) reduces to:

$$\psi^{\langle 1,n\rangle}(x_{*0},\overline{x}_{*[1,n]}) - \frac{1}{S_n(1) r_2^{n-1}} \int_{\overline{B}_n^b(r_2,\overline{x}_{*[1,n]})} \psi^{\langle 1,n\rangle}(x_{*0} + \lambda r_2, \overline{x}_{[1,n]}) ds_{\{n\}}$$
$$= - \sum_{\alpha \in \{0,1,3\}} \Omega^{(\alpha)}\left(\eta^{\langle n\rangle}, \psi^{\langle 1,n\rangle}, \tau_n, \overline{B}_n^c(r_2,\overline{x}_{*[1,n]})\right), \quad (6.1)$$

where

$$\Omega^{(0)}\left(\eta^{\langle n\rangle}, \psi^{\langle 1,n\rangle}, \tau_n, \overline{B}_n^c(r_2,\overline{x}_{*[1,n]})\right)$$
$$= -a_n \int_{\overline{B}_n^c(r_2,\overline{x}_{*[1,n]})} \frac{1}{x'^{n-2}_{[1,n]}} f^{\langle 1,n\rangle}(x_{*0} + \lambda x'_{[1,n]}, \overline{x}_{[1,n]}) dv_{\{n\}} \text{ if } n \neq 2, \quad (6.2a)$$

$$\Omega^{(0)}\left(\eta^{\langle 2\rangle},\psi^{\langle 1,2\rangle},\tau_2,\overline{B}_2^c\left(r_2,\overline{x}_{*[1,2]}\right)\right)$$
$$= a_2 \int_{\overline{B}_2^c(r_2,\overline{x}_{*[1,2]})}\left(\ln x'_{[1,2]}\right)f^{\langle 1,2\rangle}\left(x_{*0}+\lambda x'_{[1,2]},\overline{x}_{[1,2]}\right)dv_{\{2\}} \text{ if } n=2, \tag{6.2b}$$

$$\Omega^{(1)}\left(\eta^{\langle n\rangle},\psi^{\langle 1,n\rangle},\tau_n,\overline{B}_n^c\left(r_2,\overline{x}_{*[1,n]}\right)\right)$$
$$=\lambda(n-3)a_n \int_{\overline{B}_n^c(r_2,\overline{x}_{*[1,n]})}\frac{1}{x'^{n-1}_{[1,n]}}\psi_{x_0}^{\langle 1,n\rangle}\left(x_{*0}+\lambda x'_{[1,n]},\overline{x}_{[1,n]}\right)dv_{\{n\}} \text{ if } n\ne 2, \tag{6.3a}$$

$$\Omega^{(1)}\left(\eta^{\langle 2\rangle},\psi^{\langle 1,2\rangle},\tau_2,\overline{B}_{\{2}^c\left(r_2,\overline{x}_{*[1,2]}\right)\right)$$
$$=\lambda a_2 \int_{\overline{B}_2^c(r_2,\overline{x}_{*[1,2]})}\frac{\ln x'_{[1,2]}+2}{x'_{[1,2]}}\psi_{x_0}^{\langle 1,2\rangle}\left(x_{*0}+\lambda x'_{[1,2]},\overline{x}_{[1,2]}\right)dv_{\{2\}}; \tag{6.3b}$$

and also

$$\Omega^{(3)}\left(\eta^{\langle n\rangle},\psi^{\langle 1,n\rangle},\tau_n,\overline{B}_n^c\left(r_2,\overline{x}_{*[1,n]}\right)\right)=-\eta^{\langle n\rangle}(r_2)\int_{\overline{B}_n^b(r_2,\overline{x}_{*[1,n]})}\Lambda\left(\psi^{\langle 1,n\rangle},\tau_n,\overline{x}_{[1,n]}\right)ds_{\{n\}} \tag{6.4}$$

subject to (4.59) with $x'_{[1,n]}=r_2$, i.e. subject to

$$\Lambda\left(\psi^{\langle 1,n\rangle},\tau_n,\overline{x}_{[1,n]}\right)$$
$$\equiv \sum_{i=1}^n \frac{x'_i}{r_2}\psi_{x_i}^{\langle 1,n\rangle}\left(x_{*0}+\lambda r_2,\overline{x}_{[1,n]}\right)-\lambda\psi_{x_0}^{\langle 1,n\rangle}\left(x_{*0}+\lambda r_2,\overline{x}_{[1,n]}\right). \tag{6.5}$$

$\eta^{\langle n\rangle}\left(x'_{[1,n]}\right)$ is given by (4.60) subject to (4.10) and (4.12). The $n$-fold integrals occurring (6.2a)–(6.3b) are converging improper integrals of first kind.

2) For $n=1$, as $r_1$ is contracted to 0, the formula (4.38) reduces to:

$$\psi^{\langle 1,1\rangle}(x_{*0},x_{*1})-\frac{1}{2}\left[\psi^{\langle 1,1\rangle}(x_{*0}+\lambda r_2,x_{*1}-r_2)+\psi^{\langle 1,1\rangle}(x_{*0}+\lambda r_2,x_{*1}+r_2)\right]$$
$$=-\sum_{\alpha\in\{0,1,3\}}\Omega^{(\alpha)}\left(\eta^{\langle 1\rangle},\psi^{\langle 1,1\rangle},\tau_1,[x_{*1}+r_2,x_{*1}+r_2]\right), \tag{6.6}$$

where

$$\Omega^{(0)}\left(\eta^{\langle 1\rangle},\psi^{\langle 1,1\rangle},\tau_1,[x_{*1}+r_2,x_{*1}+r_2]\right)=-\frac{1}{2}\int_{x_{*1}-r_2}^{x_{*1}+r_2}|x'_1|f^{\langle 1,1\rangle}(x_{*0}+\lambda|x'_1|,x_1)dx_1, \tag{6.7}$$

$$\Omega^{(1)}\left(\eta^{\langle 1\rangle},\psi^{\langle 1,1\rangle},\tau_1,[x_{*1}+r_2,x_{*1}+r_2]\right)=\lambda\int_{x_{*1}-r_2}^{x_{*1}+r_2}\psi_{x_0}^{\langle 1,n\rangle}(x_{*0}+\lambda|x'_1|,x_1)dx_1, \tag{6.8}$$

$$\Omega^{(3)}\left(\psi^{\langle 1,1\rangle},\eta^{\langle 1\rangle},\tau_1,[x_{*1}+r_2,x_{*1}+r_2]\right)$$
$$=-\frac{r_2}{2}\left[\Psi_1(x_{*0}+\lambda r_2,x_{*1}-r_2)-\Psi_1(x_{*0}+\lambda r_2,x_{*1}+r_2)\right] \tag{6.9}$$

subject to (4.68), i.e. subject to

$$\Psi_1(x_{*0} + \lambda r_2, x_{1,}) \equiv \psi_{x_1}^{\langle 1,1 \rangle}(x_{*0} + \lambda r_2, x_{1,}) - \frac{\lambda x_1'}{r_2} \psi_{x_0}^{\langle 1,1 \rangle}(x_{*0} + \lambda r_2, x_{1,}). \tag{6.10}$$

It is understood that (6.6)–(6.10) are instances of (6.1)–(6.5) at $n=2$.

**Proof:** 1) According to Hypothesis 3.1, for each $n \in \omega_1$ the composite-plain functional form $\psi^{\langle 1,n \rangle}(x_{*0} + \lambda x'_{[1,n]}, \bar{x}_{[1,n]})$ is continuous and hence bounded in the closure of a neighborhood of the point $\bar{x}_{*[1,n]}$ – such a neighborhood, e.g., as $\bar{B}_n^o(r_1, \bar{x}_{*[1,n]})$ or $\bar{B}_n^o(r_2, \bar{x}_{*[1,n]})$.

i) Hence, it follows from (4.44) that for each $n \in \omega_2$:

$$\begin{aligned}
&\lim_{r_1 \to 0} \Sigma\left(\psi^{\langle 1,n \rangle}, \eta^{\langle n \rangle}, \tau_n, \bar{B}_n^{b-}(r_1, \bar{x}_{*[1,n]})\right) \\
&= \psi^{\langle 1,n \rangle}(x_{*0}, \bar{x}_{[1,n]}) \frac{1}{S_n(1)} \lim_{r_1 \to 0} \frac{1}{r_1^{n-1}} \int_{\bar{B}_n^b(r_1, \bar{x}_{*[1,n]})} ds_{\{n\}} \\
&= \psi^{\langle 1,n \rangle}(x_{*0}, \bar{x}_{[1,n]}) \frac{1}{S_n(1)} \lim_{r_1 \to 0} \frac{1}{r_1^{n-1}} r_1^{n-1} S_n(1) = \psi^{\langle 1,n \rangle}(x_{*0}, \bar{x}_{[1,n]}),
\end{aligned} \tag{6.11}$$

because by definition

$$\int_{\bar{B}_n^b(r_1, \bar{x}_{*[1,n]})} ds_{\{n\}} = S_n(r_1) = r_1^{n-1} S_n(1). \tag{6.12}$$

Consequently, the entire expression on the left-hand side of equation (6.1) along with the preceding sign minus, which has been transferred into the right-hand side of that equation, follows from (4.46) as $r_1$ is contracted to 0.

ii) Likewise, for $n=1$, the entire expression on the left-hand side of equation (6.6) along with the preceding sign minus, which has been transferred into the right-hand side of that equation, follows from (4.53) as $r_1$ is contracted to 0.

2. i) For $n \geq 2$, it follows from (4.59) that

$$\begin{aligned}
&\left|\Lambda\left(\psi^{\langle 1,n \rangle}, \tau_n, \bar{x}_{[1,n]}\right)\right| \\
&\leq \left|\sum_{i=1}^n \frac{x_i'}{x'_{[1,n]}} \psi_{x_i}^{\langle 1,n \rangle}(x_{*0} + \lambda x'_{[1,n]}, \bar{x}_{[1,n]})\right| + \left|\lambda \psi_{x_0}^{\langle 1,n \rangle}(x_{*0} + \lambda x'_{[1,n]}, \bar{x}_{[1,n]})\right| \\
&\leq \sum_{i=1}^n \left|\frac{x_i'}{x'_{[1,n]}}\right| \left|\psi_{x_i}^{\langle 1,n \rangle}(x_{*0} + \lambda x'_{[1,n]}, \bar{x}_{[1,n]})\right| + |\lambda| \left|\psi_{x_0}^{\langle 1,n \rangle}(x_{*0} + \lambda x'_{[1,n]}, \bar{x}_{[1,n]})\right| \\
&\leq \sum_{i=1}^n \left|\psi_{x_i}^{\langle 1,n \rangle}(x_{*0} + \lambda x'_{[1,n]}, \bar{x}_{[1,n]})\right| + \left|\psi_{x_0}^{\langle 1,n \rangle}(x_{*0} + \lambda x'_{[1,n]}, \bar{x}_{[1,n]})\right|,
\end{aligned} \tag{6.13}$$

because

$$0 \leq \left|\frac{x_i'}{x'_{[1,n]}}\right| \leq 1 \text{ for each } i \in \omega_{1,n} \tag{6.14}$$

and also because $|\lambda|=1$. Owing to Definitions 2.12 and 3.2 and owing to Hypothesis 3.1, for each $n \in \omega_1$, all composite-plain functional forms occurring in the final expression of (6.13) are bounded for each $\bar{x}_{[1,n]} \in B_n^c(r_2, \bar{x}_{*[1,n]})$, and hence so is $\Lambda(\psi^{\langle 1,n \rangle}, \tau_n, \bar{x}_{[1,n]})$ for each $n \in \omega_2$, in accordance with (6.13). Consequently, given $r_1 \in (0, r_2]$, let $M(\psi^{\langle 1,n \rangle}, \tau_n, r_1)$ be the maximum of $|\Lambda(\psi^{\langle 1,n \rangle}, \tau_n, \bar{x}_{[1,n]})|$ for all $\bar{x}_{[1,n]} \in \bar{B}_n^b(r_1, \bar{x}_{*[1,n]})$. Then it follows from (4.61) subject to (4.59) and (4.60) that for each $n \in \omega_2$:

$$\lim_{r_1 \to 0} \left| \Sigma\left(\eta^{\langle n \rangle}, \psi^{\langle 1,n \rangle}, \tau_n, \bar{B}_n^{b-}(r_1, \bar{x}_{*[1,n]})\right) \right|$$
$$\leq \lim_{r_1 \to 0} \left| \eta^{\langle n \rangle}(r_1) \right| \int_{\bar{B}_n^b(r_1, \bar{x}_{*[1,n]})} \left| \Lambda(\psi^{\langle 1,n \rangle}, \tau_n, \bar{x}_{[1,n]}) \right| ds_{\{n\}} \quad (6.15)$$
$$\leq \lim_{r_1 \to 0} \left| \eta^{\langle n \rangle}(r_1) \right| M(\psi^{\langle 1,n \rangle}, \tau_n, r_1) \int_{\bar{B}_n^b(r_1, \bar{x}_{*[1,n]})} ds_{\{n\}} = 0,$$

because $M(\psi^{\langle 1,n \rangle}, \tau_n, r_1)$ is a finite real number, while

$$\lim_{r_1 \to 0} \left| \eta^{\langle n \rangle}(r_1) \right| \int_{\bar{B}_n^b(r_1, \bar{x}_{*[1,n]})} ds_{\{n\}} = S_n(1) \lim_{r_1 \to 0} \left| \eta^{\langle n \rangle}(r_1) \right| r_1^{n-1}$$
$$= |a_n| \cdot \begin{cases} \lim_{r_1 \to 0} r_1^{2-n+n-1} = \lim_{r_1 \to 0} r_1 = 0 \text{ if } n \geq 3 \text{ (a)} \\ \lim_{r_1 \to 0} r_1 \ln r_1 = 0 \qquad \qquad \text{if } n = 2 \text{ (b)} \end{cases} \quad (6.16)$$

Hence, it follows from (4.63) that

$$\lim_{r_1 \to 0} \Omega^{(3)}\left(\eta^{\langle n \rangle}, \psi^{\langle 1,n \rangle}, \tau_n, \bar{C}_n^{cc}(r_1, r_2, \bar{x}_{*[1,n]})\right)$$
$$= \Omega^{(3)}\left(\eta^{\langle n \rangle}, \psi^{\langle 1,n \rangle}, \tau_n, \bar{B}_n^b(r_2, \bar{x}_{*[1,n]})\right) \quad (6.17)$$
$$= -\eta^{\langle n \rangle}(r_2) \int_{\bar{B}_n^b(r_2, \bar{x}_{*[1,n]})} \Lambda(\psi^{\langle 1,n \rangle}, \tau_n, \bar{x}_{[1,n]}) ds_{\{n\}},$$

which proves (6.4) subject to (6.5).

ii) For $n=1$, it follows from (4.69) subject to (4.68) that

$$\lim_{\tau_1 \to 0} \Omega^{(3)}\left(\eta^{\langle 1 \rangle}, \psi^{\langle 1,1 \rangle}, \tau_1, [x_{*1} - r_1, x_{*1} + r_1]\right) = 0. \quad (6.18)$$

Therefore, as $r_1$ is contracted to 0, (4.71) yields:

$$\lim_{\tau_1 \to 0} \Omega^{(3)}\left(\eta^{\langle 1 \rangle}, \psi^{\langle 1,1 \rangle}, \tau_1, \bar{C}_1^{cc}(r_1, r_2, x_{*1})\right)$$
$$= \Omega^{(3)}\left(\eta^{\langle 1 \rangle}, \psi^{\langle 1,1 \rangle}, \tau_1, [x_{*1} - r_2, x_{*1} + r_2]\right) \quad (6.19)$$
$$= -\frac{r_2}{2} \left[ \Psi_1(x_{*0} + \lambda r_2, x_{*1} - r_2) - \Psi_1(x_{*0} + \lambda r_2, x_{*1} + r_2) \right]$$

which proves (6.9) subject to (6.10).

3) Integrals (6.2a)–(6.3b) are improper $n$-fold integrals of first kind subject to the general definition (5.1), whereas the integrals (6.7) and (6.8) are the particular cases of (6.2a) and (6.2b) corresponding to $n=1$. That is, formally, for each $n \in \omega_1$ and for each $\alpha \in \{0,1\}$:

$$\Omega^{(\alpha)}\left(\eta^{\langle n \rangle}, \psi^{\langle 1,n \rangle}, \tau_n, \overline{B}_n^c\left(r_2, \overline{x}_{*[1,n]}\right)\right) \\ \doteq \lim_{r_1 \to 0} \Omega^{(\alpha)}\left(\eta^{\langle n \rangle}, \psi^{\langle 1,n \rangle}, \tau_n, \overline{C}_n^{cc}\left(r_1, r_2, \overline{x}_{*[1,n]}\right)\right) \quad (6.20)$$

subject to (4.72a)–(4.73b). In what follows, we shall prove that all those integrals converge. To this end, let for each $n \in \omega_1$:

$$\Xi^{(0)}\left(\psi^{\langle 1,n \rangle}, \tau_n, \overline{x}_{[1,n]}\right) \doteq -f^{\langle 1,n \rangle}\left(x_{*0} + \lambda x'_{[1,n]}, \overline{x}_{[1,n]}\right) \\ = \sum_{i=1}^{n} \psi_{x_i x_i}^{\langle 1,n \rangle}\left(x_{*0} + \lambda x'_{[1,n]}, \overline{x}_{[1,n]}\right) - \psi_{x_0 x_0}^{\langle 1,n \rangle}\left(x_{*0} + \lambda x'_{[1,n]}, \overline{x}_{[1,n]}\right), \quad (6.21)$$

$$\Xi^{(1)}\left(\psi^{\langle 1,n \rangle}, \tau_n, \overline{x}_{[1,n]}\right) \doteq \psi_{x_0}^{\langle 1,n \rangle}\left(x_{*0} + \lambda x'_{[1,n]}, \overline{x}_{[1,n]}\right). \quad (6.22)$$

Then $K^{(0)}\left(\eta^{\langle n \rangle}, \psi^{\langle 1,n \rangle}, \tau_n, \overline{x}_{[1,n]}\right)$ and $K^{(1)}\left(\eta^{\langle n \rangle}, \psi^{\langle 1,n \rangle}, \tau_n, \overline{x}_{[1,n]}\right)$, defined by (4.24) and (4.25), become

$$K^{(0)}\left(\eta^{\langle n \rangle}, \psi^{\langle 1,n \rangle}, \tau_n, \overline{x}_{[1,n]}\right) \\ = \begin{cases} \dfrac{a_n}{x'^{n-2}_{[1,n]}} \Xi^{(0)}\left(\psi^{\langle 1,n \rangle}, \tau_n, \overline{x}_{[1,n]}\right) & \text{if } n \neq 2 \text{ (a)} \\ -a_2 \left(\ln x'_{[1,2]}\right) \Xi^{(0)}\left(\psi^{\langle 1,2 \rangle}, \tau_2, \overline{x}_{[1,2]}\right) & \text{if } n = 2 \text{ (b)} \end{cases}, \quad (6.23)$$

$$K^{(1)}\left(\eta^{\langle n \rangle}, \psi^{\langle 1,n \rangle}, \tau_n, \overline{x}_{[1,n]}\right) \\ = \lambda(n-3) \cdot \begin{cases} \dfrac{a_n}{x'^{n-1}_{[1,n]}} \Xi^{(1)}\left(\psi^{\langle 1,n \rangle}, \tau_n, \overline{x}_{[1,n]}\right) & \text{if } n \neq 2 \text{ (a)} \\ -\dfrac{a_2 \left(\ln x'_{[1,2]} + 2\right)}{x'_{[1,2]}} \Xi^{(1)}\left(\psi^{\langle 1,2 \rangle}, \tau_2, \overline{x}_{[1,2]}\right) & \text{if } n = 2 \text{ (b)} \end{cases}. \quad (6.24)$$

Owing to Definitions 2.12 and 3.2 and Hypothesis 3.1, for each $n \in \omega_1$, the composite-plain functional forms $\psi_{x_i x_i}^{\langle 1,n \rangle}\left(x_{*0} + \lambda x'_{[1,n]}, \overline{x}_{[1,n]}\right)$ and $\psi_{x_0 x_0}^{\langle 1,n \rangle}\left(x_{*0} + \lambda x'_{[1,n]}, \overline{x}_{[1,n]}\right)$ occurring in $\Xi^{(0)}\left(\psi^{\langle 1,n \rangle}, \tau_n, \overline{x}_{[1,n]}\right)$, defined by (6.21), and the composite-plain functional form $\psi_{x_0}^{\langle 1,n \rangle}\left(x_{*0} + \lambda x'_{[1,n]}, \overline{x}_{[1,n]}\right)$ occurring in $\Xi^{(1)}\left(\psi^{\langle 1,n \rangle}, \tau_n, \overline{x}_{[1,n]}\right)$, defined by (6.22), are bounded for each $\overline{x}_{[1,n]} \in B_n^c\left(r_2, \overline{x}_{*[1,n]}\right)$. we shall therefore take for granted the following property of boundedness of the functional forms $\Xi^{(0)}\left(\psi^{\langle 1,n \rangle}, \tau_n, \overline{x}_{[1,n]}\right)$ and $\Xi^{(1)}\left(\psi^{\langle 1,n \rangle}, \tau_n, \overline{x}_{[1,n]}\right)$.

**Axiom 6.1.** Given $n \in \omega_1$, given $\alpha \in \{0,1\}$, there exist two real numbers $r_0^{(\alpha)} \in (0,1)$ and $c_n^{(\alpha)} > 0$ such that

$$\left| \Xi^{(\alpha)}\left(\psi^{\langle 1,n \rangle}, \tau_n, \bar{x}_{[1,n]}\right) \right| \leq c_n^{(\alpha)} \text{ for every } \bar{x}_{[1,n]} \in \bar{B}_n^c\left(r_0^{(\alpha)}, \bar{x}_{*[1,n]}\right) \tag{6.25}$$

In this case, without loss of generality, it is assumed that $0 < r_1 < r_0 < r_2 < \infty$, where $r_0$ is the *smallest* one of the two numbers $r_0^{(0)}$ and $r_0^{(1)}$.

i) By Axiom 6.1, equations (6.23a) and (6.24a) imply the following corollary.

**Corollary 6.1.** Given $n \in \omega_1 - \{2\}$, given $\alpha \in \{0,1\}$, there exist three real numbers $r_0 \in (0,1)$, $p_n^{(\alpha)} < n$, and $k_n^{(\alpha)} > 0$ such that

$$\left| K^{(\alpha)}\left(\eta^{\langle n \rangle}, \psi^{\langle 1,n \rangle}, \tau_n, \bar{x}_{[1,n]}\right) \right| \leq \frac{k_n^{(\alpha)}}{\left| \bar{x}_{[1,n]} - \bar{x}_{*[1,n]} \right|^{p_n^{(\alpha)}}} \tag{6.26}$$

for every $\bar{x}_{[1,n]} \in \bar{B}_n^c\left(r_0, \bar{x}_{*[1,n]}\right)$,

the understanding being that

$$k_n^{(0)} \doteq a_n c_n^{(0)}, \; k_n^{(1)} \doteq |n-3| a_n c_n^{(1)}, \; p_n^{(0)} = n-2, \; p_n^{(1)} = n-1. \tag{6.27}$$

The above corollary is an instance of the hypothesis of Theorem 5.2. Therefore, the following proofs of convergence of integrals (6.2a) and (6.3a) are the pertinent instances of the proof of Theorem 5.2.

The integrals (4.72a) and (4.73a) can be represented as:

$$\Omega^{(\alpha)}\left(\eta^{\langle n \rangle}, \psi^{\langle 1,n \rangle}, \tau_n, \bar{C}_n^{cc}\left(r_1, r_2, \bar{x}_{*[1,n]}\right)\right) = \Omega^{(\alpha)}\left(\eta^{\langle n \rangle}, \psi^{\langle 1,n \rangle}, \tau_n, \bar{C}_n^{cc}\left(r_1, r_0, \bar{x}_{*[1,n]}\right)\right) \\ + \Omega^{(\alpha)}\left(\eta^{\langle n \rangle}, \psi^{\langle 1,n \rangle}, \tau_n, \bar{C}_n^{cc}\left(r_0, r_2, \bar{x}_{*[1,n]}\right)\right), \tag{6.28}$$

for $\alpha=0$ and $\alpha=1$ respectively. The first integral on the right-hand side of this identity can be majorated (estimated by absolute value from above) with the help of (6.26) thus:

$$\left| \Omega^{(\alpha)}\left(\eta^{\langle n \rangle}, \psi^{\langle 1,n \rangle}, \tau_n, \bar{C}_n^{cc}\left(r_1, r_0, \bar{x}_{*[1,n]}\right)\right) \right|$$

$$\leq \int_{\bar{C}_n^{cc}(r_1, r_0, \bar{x}_{*[1,n]})} \left| K^{(\alpha)}\left(\eta^{\langle n \rangle}, \psi^{\langle 1,n \rangle}, \tau_n, \bar{x}_{[1,n]}\right) \right| dv_{\{n\}} \leq \int_{r_1}^{r_0} \frac{k_n^{(\alpha)} dV_n(r)}{r} \tag{6.29}$$

$$= k_n^{(\alpha)} S_n(1) \int_{r_1}^{r_0} r^{n-p_n^{(\alpha)}-1} dr = \frac{k_n^{(\alpha)} S_n(1)}{n - p_n^{(\alpha)}} \left( r_0^{n-p_n^{(\alpha)}} - r_1^{n-p_n^{(\alpha)}} \right) \leq \frac{k_n^{(\alpha)} S_n(1)}{n - p_n^{(\alpha)}} r_0^{n-p_n^{(\alpha)}}$$

(cf. (5.11) and (5.10)). Consequently, by (6.29), it follows that

$$\left| \Omega^{(\alpha)}\left(\eta^{\langle n \rangle}, \psi^{\langle 1,n \rangle}, \tau_n, \bar{B}_n^c\left(r_0, \bar{x}_{*[1,n]}\right)\right) \right|$$

$$= \lim_{r_1 \to 0} \left| \Omega^{(\alpha)}\left(\eta^{\langle n \rangle}, \psi^{\langle 1,n \rangle}, \tau_n, \bar{C}_n^{cc}\left(r_1, r_0, \bar{x}_{*[1,n]}\right)\right) \right| \leq \frac{k_n^{(\alpha)} S_n(1)}{n - p_n^{(\alpha)}} r_0^{n-p_n^{(\alpha)}}, \tag{6.30}$$

which is the pertinent instance of (5.10). By (6.23), at $\alpha=0$ or $\alpha=0$, (6.30) yields:

$$\left|\Omega^{(0)}\left(\eta^{\langle n\rangle},\psi^{\langle 1,n\rangle},\tau_n,\overline{B}_n^c\left(r_0,\overline{x}_{*[1,n]}\right)\right)\right| \leq \frac{1}{2}k_n^{(0)}S_n(1)r_0^2, \tag{6.31}$$

$$\left|\Omega^{(1)}\left(\eta^{\langle n\rangle},\psi^{\langle 1,n\rangle},\tau_n,\overline{B}_n^c\left(r_0,\overline{x}_{*[1,n]}\right)\right)\right| \leq k_n^{(1)}S_n(1)r_0. \tag{6.32}$$

The second integral on the right-hand side of this identity (6.28) is an ordinary Riemann integral. Therefore, for each $\alpha \in \{0,1\}$, (6.20) subject to (6.28) becomes

$$\begin{aligned}&\Omega^{(\alpha)}\left(\eta^{\langle n\rangle},\psi^{\langle 1,n\rangle},\tau_n,\overline{B}_n^c\left(r_2,\overline{x}_{*[1,n]}\right)\right)\\ &= \Omega^{(\alpha)}\left(\eta^{\langle n\rangle},\psi^{\langle 1,n\rangle},\tau_n,\overline{B}_n^c\left(r_0,\overline{x}_{*[1,n]}\right)\right)\\ &+ \Omega^{(\alpha)}\left(\eta^{\langle n\rangle},\psi^{\langle 1,n\rangle},\tau_n,\overline{C}_n^{cc}\left(r_0,r_2,\overline{x}_{*[1,n]}\right)\right),\end{aligned} \tag{6.33}$$

subject (6.30)–(6.32), which proves convergence of improper integrals (6.2a) and (6.3a).

ii) At $n=2$, definitions (6.21) and (6.22) become

$$\begin{aligned}\Xi^{(0)}\left(\psi^{\langle 1,2\rangle},\tau_2,\overline{x}_{[1,2]}\right) &\equiv -f^{\langle 1,2\rangle}\left(x_{*0}+\lambda x'_{[1,2]},\overline{x}_{[1,2]}\right)\\ &= \sum_{i=1}^{2}\psi^{\langle 1,2\rangle}_{x_ix_i}\left(x_{*0}+\lambda x'_{[1,2]},\overline{x}_{[1,2]}\right) - \psi^{\langle 1,2\rangle}_{x_0x_0}\left(x_{*0}+\lambda x'_{[1,2]},\overline{x}_{[1,2]}\right),\end{aligned} \tag{6.34}$$

$$\Xi^{(1)}\left(\psi^{\langle 1,2\rangle},\tau_2,\overline{x}_{[1,2]}\right) \equiv \psi^{\langle 1,2\rangle}_{x_0}\left(x_{*0}+\lambda x'_{[1,2]},\overline{x}_{[1,2]}\right). \tag{6.35}$$

At the same time, by (6.23b) and (6.24b), the integrals (4.72b) and (4.73b) subject to (4.12) can be rewritten as:

$$\begin{aligned}&\Omega^{(0)}\left(\eta^{\langle 2\rangle},\psi^{\langle 1,2\rangle},\tau_2,\overline{C}_2^{cc}\left(r_1,r_2,\overline{x}_{*[1,2]}\right)\right)\\ &= \int_{\overline{C}_2^{cc}(r_1,r_2,\overline{x}_{*[1,2]})} K^{(0)}\left(\eta^{\langle 2\rangle},\psi^{\langle 1,2\rangle},\tau_n,\overline{x}_{[1,2]}\right)dv_{\{2\}}\\ &= -\frac{1}{2\pi}\int_{\overline{C}_2^{cc}(r_1,r_2,\overline{x}_{*[1,2]})}\left(\ln x'_{[1,2]}\right)\Xi^{(0)}\left(\psi^{\langle 1,2\rangle},\tau_2,\overline{x}_{[1,2]}\right)dv_{\{2\}},\end{aligned} \tag{6.36}$$

$$\begin{aligned}&\Omega^{(1)}\left(\eta^{\langle 2\rangle},\psi^{\langle 1,2\rangle},\tau_2,\overline{C}_2^{cc}\left(r_1,r_2,\overline{x}_{*[1,2]}\right)\right)\\ &= \int_{\overline{C}_2^{cc}(r_1,r_2,\overline{x}_{*[1,2]})} K^{(1)}\left(\eta^{\langle 2\rangle},\psi^{\langle 1,2\rangle},\tau_n,\overline{x}_{[1,2]}\right)dv_{\{2\}}\\ &= \frac{\lambda}{2\pi}\int_{\overline{C}_2^{cc}(r_1,r_2,\overline{x}_{*[1,2]})} x'^{-1}_{[1,2]}\left(\ln x'_{[1,2]}+2\right)\Xi^{(1)}\left(\psi^{\langle 1,2\rangle},\tau_2,\overline{x}_{[1,2]}\right)dv_{\{2\}}.\end{aligned} \tag{6.37}$$

In a polar coordinate system $(r,\varphi)$ such that

$$\overline{x}_{[1,2]} = \overline{x}_{*[1,2]} + r\overline{n}_{[1,2]}(\varphi),\ x'_{[1,2]} = r,\ \overline{n}_{[1,2]}(\varphi) = \langle\cos\varphi,\sin\varphi\rangle, \tag{6.38}$$

the integral (6.36) becomes:

$$\Omega^{(0)}\left(\eta^{\langle 2\rangle},\psi^{\langle 1,2\rangle},\tau_2,\overline{C}_2^{cc}\left(r_1,r_2,\overline{x}_{*[1,2]}\right)\right) = \int_{r_1}^{r_2}\Psi^{(0)}\left(\psi^{\langle 1,2\rangle},\tau_2,r\right)r\ln r\,dr, \tag{6.39}$$

where

$$\Psi^{(0)}(\psi^{\langle 1,2\rangle},\tau_2,r) \doteq -\frac{1}{2\pi}\int_0^{2\pi}\Xi^{(0)}(\psi^{\langle 1,2\rangle},\tau_2,\bar{x}_{*[1,2]}+r\bar{n}_{[1,2]}(\varphi))d\varphi$$
$$= \frac{1}{2\pi}\int_0^{2\pi} f^{\langle 1,2\rangle}(x_{*0}+\lambda x'_{[1,2]},\bar{x}_{*[1,2]}+r\bar{n}_{[1,2]}(\varphi))d\varphi, \qquad (6.40)$$

and integral (6.37) can conveniently be written as:

$$\Omega^{(1)}(\eta^{\langle 2\rangle},\psi^{\langle 1,2\rangle},\tau_2,\overline{C}_2^{cc}(r_1,r_2,\bar{x}_{*[1,2]}))$$
$$= \Omega_1^{(1)}(\eta^{\langle 2\rangle},\psi^{\langle 1,2\rangle},\tau_2,\overline{C}_2^{cc}(r_1,r_2,\bar{x}_{*[1,2]})) \qquad (6.41)$$
$$+ \Omega_2^{(1)}(\eta^{\langle 2\rangle},\psi^{\langle 1,2\rangle},\tau_2,\overline{C}_2^{cc}(r_1,r_2,\bar{x}_{*[1,2]}))$$

subject to

$$\Omega_1^{(1)}(\eta^{\langle 2\rangle},\psi^{\langle 1,2\rangle},\tau_2,\overline{C}_2^{cc}(r_1,r_2,\bar{x}_{*[1,2]})) = \int_{r_1}^{r_2}\Psi^{(1)}(\psi^{\langle 1,2\rangle},\tau_2,r)\ln r\, dr, \qquad (6.42)$$

$$\Omega_1^{(2)}(\eta^{\langle 2\rangle},\psi^{\langle 1,2\rangle},\tau_2,\overline{C}_2^{cc}(r_1,r_2,\bar{x}_{*[1,2]})) = 2\int_{r_1}^{r_2}\Psi^{(1)}(\psi^{\langle 1,2\rangle},\tau_2,r)dr, \qquad (6.43)$$

where

$$\Psi^{(1)}(\psi^{\langle 1,2\rangle},\tau_2,r) \doteq \frac{\lambda}{2\pi}\int_0^{2\pi}\Xi^{(1)}(\psi^{\langle 1,2\rangle},\tau_2,\bar{x}_{*[1,2]}+r\bar{n}_{[1,2]}(\varphi))d\varphi$$
$$= \frac{\lambda}{2\pi}\int_0^{2\pi}\psi_{x_0}^{\langle 1,2\rangle}(x_{*0}+\lambda x'_{[1,2]},\bar{x}_{*[1,2]}+r\bar{n}_{[1,2]}(\varphi))d\varphi. \qquad (6.44)$$

By Axiom 6.1, it follows from (6.40) and (6.44) that for each $\alpha \in \{0,1\}$:

$$\left|\Psi^{(\alpha)}(\psi^{\langle 1,2\rangle},\tau_2,r)\right| \le \frac{1}{2\pi}\int_0^{2\pi}\left|\Xi^{(\alpha)}(\psi^{\langle 1,2\rangle},\tau_2,\bar{x}_{*[1,2]}+r\bar{n}_{[1,2]}(\varphi))\right|d\varphi$$
$$= \frac{1}{2\pi}\int_0^{2\pi} c_2^{(\alpha)}d\varphi = c_2^{(\alpha)} \text{ for each } r \in (r_1,r_0), \qquad (6.45)$$

At the same time, each one of the integrals (6.39), (6.42), and (6.43) can be represented in the form analogous to (6.28), e.g.

$$\Omega^{(0)}(\eta^{\langle 2\rangle},\psi^{\langle 1,2\rangle},\tau_2,\overline{C}_2^{cc}(r_1,r_2,\bar{x}_{*[1,2]}))$$
$$= \Omega^{(0)}(\eta^{\langle 2\rangle},\psi^{\langle 1,2\rangle},\tau_2,\overline{C}_2^{cc}(r_1,r_0,\bar{x}_{*[1,2]})) \qquad (6.46)$$
$$+ \Omega^{(0)}(\eta^{\langle 2\rangle},\psi^{\langle 1,2\rangle},\tau_2,\overline{C}_2^{cc}(r_0,r_1,\bar{x}_{*[1,2]}))$$

and similarly with '$\Omega_1^{(1)}$' or '$\Omega_1^{(1)}$' in place of '$\Omega^{(0)}$'. In this case, by (6.45), it follows from the variants of (6.39), (6.42), and (6.43) with $r_0$ in place of $r_2$ that

$$\left|\Omega^{(0)}\left(\eta^{\langle 2\rangle},\psi^{\langle 1,2\rangle},\tau_2,\overline{C}_2^{cc}(r_1,r_0,\overline{x}_{*[1,2]})\right)\right| = \left|\int_{r_1}^{r_0}\Psi^{(0)}\left(\psi^{\langle 1,2\rangle},\tau_2,r\right)r\ln r\,dr\right|$$

$$\leq -\int_{r_1}^{r_0}\left|\Psi^{(0)}\left(\psi^{\langle 1,2\rangle},\tau_2,r\right)\right|r\ln r\,dr \leq -c_2^{(0)}\int_{r_1}^{r_0}r\ln r\,dr \quad (6.47)$$

$$= \frac{1}{2}c_2^{(0)}\left[r_1^2\left(\ln r_1 - \frac{1}{2}\right) - r_0^2\left(\ln r_0 - \frac{1}{2}\right)\right],$$

$$\left|\Omega_1^{(1)}\left(\eta^{\langle 2\rangle},\psi^{\langle 1,2\rangle},\tau_2,\overline{C}_2^{cc}(r_1,r_0,\overline{x}_{*[1,2]})\right)\right| = \left|\int_{r_1}^{r_0}\Psi^{(1)}\left(\psi^{\langle 1,2\rangle},\tau_2,r\right)r\ln r\,dr\right|$$

$$\leq -\int_{r_1}^{r_0}\left|\Psi^{(1)}\left(\psi^{\langle 1,2\rangle},\tau_2,r\right)\right|\ln r\,dr \leq -c_2^{(1)}\int_{r_1}^{r_0}\ln r\,dr \quad (6.48)$$

$$= c_2^{(1)}\left[r_1(\ln r_1 - 1) - r_0(\ln r_0 - 1)\right],$$

$$\left|\Omega_2^{(1)}\left(\eta^{\langle 2\rangle},\psi^{\langle 1,2\rangle},\tau_2,\overline{C}_2^{cc}(r_1,r_0,\overline{x}_{*[1,2]})\right)\right| = 2\left|\int_{r_1}^{r_0}\Psi^{(1)}\left(\psi^{\langle 1,2\rangle},\tau_2,r\right)dr\right|$$

$$\leq c_2^{(1)}\int_{r_1}^{r_0}dr = c_2^{(1)}(r_0 - r_1). \quad (6.49)$$

By Axiom 6.1, $r_0 \in (0,1)$ and hence

$$\ln r < 0 \text{ for each } r \in (0, r_0). \quad (6.50)$$

Therefore, the sign '−', occurring in certain terms of the trains (6.47) and (6.48), guarantees that these terms are *strictly positive*. The following two elementary integrals are calculated by integration by parts:

$$\int_{r_1}^{r_2} r\ln r\,dr = \frac{1}{2}\left(r^2\ln r\right)_{r_1}^{r_2} - \frac{1}{2}\int_{r_1}^{r_2}r\,dr = \frac{1}{2}\left(r_2^2\ln r_2 - r_1^2\ln r_1\right)$$

$$-\frac{1}{4}\left(r_2^2 - r_1^2\right) = \frac{1}{2}\left[r_2^2\left(\ln r_2 - \frac{1}{2}\right) - r_1^2\left(\ln r_1 - \frac{1}{2}\right)\right], \quad (6.51)$$

$$\int_{r_1}^{r_2}\ln r\,dr = \left(r\ln r\right)_{r_1}^{r_2} - \int_{r_1}^{r_2}dr = r_2(\ln r_2 - 1) - r_1(\ln r_1 - 1). \quad (6.52)$$

From (6.47)–(6.49), it follows that

$$\lim_{r_1 \to 0}\left|\Omega^{(0)}\left(\eta^{\langle 2\rangle},\psi^{\langle 1,2\rangle},\tau_2,\overline{C}_2^{cc}(r_1,r_0,\overline{x}_{*[1,2]})\right)\right| \leq -\frac{1}{2}c_2^{(0)}r_0^2\left(\ln r_0 - \frac{1}{2}\right) > 0, \quad (6.53)$$

$$\lim_{r_1 \to 0}\left|\Omega_1^{(1)}\left(\eta^{\langle 2\rangle},\psi^{\langle 1,2\rangle},\tau_2,\overline{C}_2^{cc}(r_1,r_0,\overline{x}_{*[1,2]})\right)\right| \leq -c_2^{(1)}r_0(\ln r_0 - 1) > 0, \quad (6.54)$$

$$\lim_{r_1 \to 0}\left|\Omega_2^{(1)}\left(\eta^{\langle 2\rangle},\psi^{\langle 1,2\rangle},\tau_2,\overline{C}_2^{cc}(r_1,r_0,\overline{x}_{*[1,2]})\right)\right| \leq c_2^{(1)}r_0 > 0, \quad (6.55)$$

because, by l'Hospitale's rule,

$$\lim_{r_1 \to 0}(r_1^m \ln r_1) = \lim_{r_1 \to 0}\left[\frac{d \ln r_1}{dr_1}\left(\frac{dr_1^{-m}}{dr_1}\right)^{-1}\right] = \lim_{r_1 \to 0}\left[r_1^{-1}\left(-mr_1^{-(m+1)}\right)^{-1}\right] \qquad (6.56)$$
$$= -\frac{1}{m}\lim_{r_1 \to 0} r_1^m = 0 \text{ for each } m \in \omega_1.$$

Thus, the integrals (6.39), (6.42), and (6.43) converge as $r_1 \to \infty$. QED.•

**Comment 6.1.** If the functional form '$\psi^{\langle 1,n \rangle}(x_0, \bar{x}_{[1,n]})$' is regarded as a solution of the inhomogeneous differential equation (4.75) and if, hence, the functional form '$f^{\langle 1,n \rangle}(x_0, \bar{x}_{[1,n]})$' is regarded as a *given* one, and not as the definiendum of definition (4.21) then Theorem 6.1 has the following important implication. The *special second integro-differential formula* (*SIDF*) (6.1) subject to (6.2a)–(6.5) uniquely determines the value of the functional form '$\psi^{\langle 1,n \rangle}(x_{*0}, \bar{x}_{*[1,n]})$' at any given point $\langle x_{*0}, \bar{x}_{*[1,n]} \rangle \in R \times \overline{B}^o_{\{n\}}(r_2, \bar{x}_{*[1,n]})$ and hence it uniquely determines the function $\psi$ on the open sphere $\overline{B}^o_{\{n\}}(r_2, \bar{x}_{*[1,n]})$ if and only if the following conditions are satisfied in that order:

a) Exactly one of the two values -1 and 1 of the parameter '$\lambda$' is fixed.

b) $n=3$.

c) The values of the functional forms '$\psi^{\langle 1,3 \rangle}(x_{*0} + \lambda r_2, \bar{x}_{[1,3]})$' and '$\psi^{\langle 1,3 \rangle}_{x_i}(x_{*0} + \lambda r_2, \bar{x}_{[1,3]})$' at each $i \in \omega_{0,3}$ are given for each $\bar{x}_{[1,3]} \in \overline{B}^b_3(r_2, \bar{x}_{*[1,3]})$, i.e. on the 2-dimensional boundary spherical surface. Hence, the 2-fold surface integral, occurring on the left-hand side of the identity (6.1) at $n=3$ and the 2-fold surface integral, occurring on the right-hand side of the identity (6.4) subject to (6.5) at $n=3$ are given.

In principle, the order of conditions b) and c) and hence the very condition c) can be changed as follows:

b′) The values of the functional forms '$\psi^{\langle 1,n \rangle}(x_{*0} + \lambda r_2, \bar{x}_{[1,n]})$' and '$\psi^{\langle 1,n \rangle}_{x_i}(x_{*0} + \lambda x'_{[1,n]}, \bar{x}_{[1,n]})$' for each $i \in \omega_{0,n}$ are given for each $\bar{x}_{[1,n]} \in \overline{B}^b_n(r_2, \bar{x}_{*[1,n]})$, i.e. on the ($n$–1)-dimensional boundary spherical surface. Hence, the ($n$–1)-fold surface integral, occurring on the left-hand side of the identity (6.1) and the ($n$–1)-fold surface integral, occurring on the right-hand side of the identity (6.4) subject to (6.5), are given.

c′) $n=3$.

In this case, condition c) is the instance of condition b′) subject to n=3, i.e. subject to condition c′).•

# 7. The special second integro-differential formula for a biased functional form in an infinite *n*-dimensional Euclidean space

## 7.1. The orthochronous axiom

**Preliminary Remark 7.1.** After condition a) of Comment 6.1 is satisfied, the radius $r_2$ of the *n*-dimensional sphere $B_n^c(r_2, \bar{x}_{*[1,n]})$ can be made as large as desired. As $r_2 \to \infty$, condition (3.1) that determines differential properties of the functional forms '$\psi^{\langle 1,n \rangle}(x_0, \bar{x}_{[1,n]})$' and '$\psi^{\langle 1,n \rangle}(x_{*0} + \lambda r_2, \bar{x}_{[1,n]})$' becomes:

$$\psi^{\langle 1,n \rangle} \in \mathcal{C}_2(R \times \bar{E}_n), \quad (7.1)$$

while the (*n*–1)-dimensional boundary spherical surface $B_n^b(r_2, \bar{x}_{*[1,n]})$ is said to be an *infinitely remote* one. In order to allow the limiting transition $r_2 \to \infty$ in the SIDF (6.1) subject to (6.2a)–(6.5), we shall adopt the following axiom.

**Axiom 7.1:** *The orthochronous axiom.* Given $\lambda \in \{-1,1\}$, given $n \in \omega_1$, there is a real-valued function $\psi^{\langle 1,n \rangle}$ subject to (7.1), which satisfies Theorem 6.1 and which has the following additional property. With $\dot{r}_0$ being the *largest* one of the two numbers $r_0^{(0)}$ and $r_0^{(1)}$ of Axiom 6.1, if $r_2$ becomes large enough ($r_2 \to \infty$) then there is a strictly positive real number $r_* \in (\dot{r}_0, r_2)$ such that for each $m \in \omega_{0,2}$ there are two strictly positive real numbers $p_{*n,m} > m$, and $C_{n,m} > 0$ such that

$$\left. \frac{\left|\partial^m \psi^{\langle 1,n \rangle}(x_0, \bar{x}_{[1,n]})\right|}{\partial x_i^m} \right|_{x_0 = x_{*0} + \lambda x'_{[1,n]}} \leq \frac{C_{n,m}}{\left|\bar{x}_{[1,n]} - \bar{x}_{*[1,n]}\right|^{p_{*n,m}}} \quad (7.2)$$

for each $i \in \omega_{0,n}$ and every $\left|\bar{x}_{[1,n]} - \bar{x}_{*[1,n]}\right| > r_*$;

the inequality $\left|\bar{x}_{[1,n]} - \bar{x}_{*[1,n]}\right| > r_*$ is equivalent to the relation $\bar{x}_{[1,n]} \in \bar{E}_n - \bar{B}_n^c(r_*, \bar{x}_{*[1,n]})$.•

**Lemma 7.1.** Under definition (4.21) or under inhomogeneous differential equation (4.75), given $\lambda \in \{-1,1\}$, given $n \in \omega_1$,

$$\left|f^{\langle 1,n \rangle}(x_{*0} + \lambda x'_{[1,n]}, \bar{x}_{[1,n]})\right| \leq \frac{(n+1)C_{n,2}}{\left|\bar{x}_{[1,n]} - \bar{x}_{*[1,n]}\right|^{p_{*n,2}}} \text{ for each } \left|\bar{x}_{[1,n]} - \bar{x}_{*[1,n]}\right| > r_*. \quad (7.3)$$

**Proof:** By (7.2) at $m=2$, it follows either from definition (4.21) or from equation (4.75), subject to (3.12) that

$$\left| f^{\langle 1,n \rangle}\left(x_{*0} + \lambda x'_{[1,n]}, \overline{x}_{[1,n]}\right) \right| \leq \left| \sum_{i=1}^{n} \psi^{\langle 1,n \rangle}_{x_i x_i}\left(x_0, \overline{x}_{[1,n]}\right) - \psi^{\langle 1,n \rangle}_{x_0 x_0}\left(x_0, \overline{x}_{[1,n]}\right) \right|_{x_0 = x_{*0} + \lambda x'_{[1,n]}}$$

$$\leq \left[ \sum_{i=1}^{n} \left| \psi^{\langle 1,n \rangle}_{x_i x_i}\left(x_0, \overline{x}_{[1,n]}\right) \right| + \left| \psi^{\langle 1,n \rangle}_{x_0 x_0}\left(x_0, \overline{x}_{[1,n]}\right) \right| \right]_{x_0 = x_{*0} + \lambda x'_{[1,n]}} \leq \frac{(n+1)C_{n,2}}{\left| \overline{x}_{[1,n]} - \overline{x}_{*[1,n]} \right|^{p_{*n,2}}} \quad (7.3_1)\bullet$$

for each $\left| \overline{x}_{[1,n]} - \overline{x}_{*[1,n]} \right| > r_*$.

**Comment 7.1.** 1) Since $x_{*0} - r_2 \to -\infty$ and $x_{*0} + r_2 \to +\infty$ as $r_2 \to \infty$, it is therefore understood that if a certain functional form $\psi^{\langle 1,n \rangle}\left(x_{*0} - r_2, \overline{x}_{[1,n]}\right)$, having the respective *associative function* $\psi^{\langle 1,n \rangle}$, and its pertinent partial derivatives satisfy Axiom 7.1 with $\lambda=-1$ than that axiom is not supposed to apply with $\lambda=1$ to the functional form $\psi^{\langle 1,n \rangle}\left(x_{*0} + r_2, \overline{x}_{[1,n]}\right)$, having the *same* associative function. That is to say, according to Axiom 7.1, all functions of the given class $\mathcal{C}_2(R \times \overline{E}_n)$ are divided into two subclasses: those satisfying condition (7.2) at $\lambda=-1$ and those satisfying condition (7.2) at $\lambda=1$, the understanding being that one and the same function cannot satisfy condition (7.2) at both values of '$\lambda$'. Consequently, for avoidance confusion, instead of the functional variable '$\psi^{\langle 1,n \rangle}$', whose range comprises functions of both subclasses, we can introduce two different functional variables: '$\psi^{-1\langle 1,n \rangle}$', whose range comprises functions of the former subclass, and '$\psi^{1\langle 1,n \rangle}$', whose range comprises functions of the latter subclass. In this case, a function $\psi^{-1\langle 1,n \rangle}$ is said to be *retarded* in the sense that it can be associated only with a retarded functional form and a function $\psi^{1\langle 1,n \rangle}$ is said to be *advanced* in the sense that it can be associated only with an advanced functional form. Accordingly, before passing to the limit $r_2 \to \infty$ in the SIDF (6.1), the latter is supposed to be specified, explicitly or mentally, for composite-plain functional forms of each one of the two kinds separately, namely for $\psi^{-1\langle 1,n \rangle}\left(x_{*0} - x'_{[1,n]}, \overline{x}_{[1,n]}\right)$, which satisfies condition (7.2) at $\lambda=-1$, and for $\psi^{1\langle 1,n \rangle}\left(x_{*0} + x'_{[1,n]}, \overline{x}_{[1,n]}\right)$, which satisfies condition (7.2) at $\lambda=1$. At the same time, in order to treat of both above functional forms simultaneously, we should use the placeholder '$\psi^{\lambda\langle 1,n \rangle}\left(x_{*0} + \lambda x'_{[1,n]}, \overline{x}_{[1,n]}\right)$', which is as ambiguous as the former one '$\psi^{\langle 1,n \rangle}\left(x_{*0} + \lambda r_2, \overline{x}_{[1,n]}\right)$'. Like remarks apply with '$f^{-1\langle 1,n \rangle}$', '$f^{1\langle 1,n \rangle}$', and '$f^{\lambda\langle 1,n \rangle}$' in place of '$\psi^{-1\langle 1,n \rangle}$', '$\psi^{1\langle 1,n \rangle}$', and '$\psi^{\lambda\langle 1,n\rangle}$' respectively. In the sequel, we shall use '$\psi^{\lambda\langle 1,n\rangle}$' and '$f^{\lambda\langle 1,n\rangle}$' interchangeably with '$\psi^{\langle 1,n\rangle}$' and '$f^{\langle 1,n\rangle}$' respectively and we shall specify '$\psi^{\langle 1,n\rangle}$' as '$\psi^{-1\langle 1,n\rangle}$' or '$\psi^{1\langle 1,n\rangle}$' and '$f^{\langle 1,n\rangle}$' as '$f^{-1\langle 1,n\rangle}$' or '$f^{1\langle 1,n\rangle}$' without any further comments.

2) Since Axiom 7.1 discriminates, by means of its condition (7.2), between retarded and advanced functions, it is qualified *orthochronous*, although this qualifier is, strictly speaking, justified only in the case of $n=3$; in all other cases, it would have been more correct to qualify that axiom *quasi-orthochronous*.

3) Condition (7.2) at $m=0$ guarantees that the second term on the left-hand side of identity (6.1) and that the second term on the left-hand side of identity (6.6) vanish as $r_2 \to \infty$. Condition (7.2) at $m=1$ guarantees that the expression on the right-hand side of identity (6.4) and that the expression on the right-hand side of identity (6.9) vanish as $r_2 \to \infty$. Also, as $r_2 \to \infty$, condition (7.2) at $m=1$ and $i=0$ guarantees convergence of the integrals occurring in (6.3a), (6.3b), and (6.8), whereas that condition at $m=2$ guarantees convergence of the integrals occurring in (6.2a), (6.2b), and (6.7). Thus, in the latter two cases, Axiom 7.1 adjusts the asymptotic behavior of the pertinent derivative biased functional forms to the hypotheses of Theorems 5.3 and 5.4. The above said is demonstrated by the proofs of the following Lemmas 7.2–7.4.•

## 7.2. Derivation of the special second integro-differential formula for a biased functional form in an infinite *n*-dimensional Euclidean space

**Lemma 7.2.** Given $\lambda \in \{-1, 1\}$,

$$\lim_{r_2 \to 0} \frac{1}{S_n(1) r_2^{n-1}} \int_{\overline{B}_n^b(r_2, \overline{x}_{*[1,n]})} \psi^{\langle 1,n\rangle}(x_{*0} + \lambda r_2, \overline{x}_{[1,n]}) ds_{\{n\}} = 0, \quad (7.4)$$

$$\lim_{r_2 \to 0} \Omega^{(3)}\left(\eta^{\langle n\rangle}, \psi^{\langle 1,n\rangle}, \tau_n, \overline{B}_n^c(r_2, \overline{x}_{*[1,n]})\right) = 0 \quad (7.5)$$

for each $n \in \omega_2$, and

$$\lim_{r_2 \to 0} \frac{1}{2}\left[\psi^{\langle 1,1\rangle}(x_{*0} + \lambda r_2, x_{*1} - r_2) + \psi^{\langle 1,1\rangle}(x_{*0} + \lambda r_2, x_{*1} + r_2)\right] = 0, \quad (7.6)$$

$$\lim_{r_2 \to 0} \Omega^{(3)}\left(\eta^{\langle 1\rangle}, \psi^{\langle 1,1\rangle}, \tau_1, [x_{*1} + r_2, x_{*1} + r_2]\right) = 0. \quad (7.7)$$

**Proof:** 1) By (7.2) with $n \geq 2$ and $m=0$, it follows that if $r_2 > r_*$ then

$$\left| \int_{\overline{B}_n^b(r_2,\overline{x}_{*[1,n]})} \psi^{\langle 1,n \rangle}(x_{*0}+\lambda r_2, \overline{x}_{[1,n]}) ds_{\{n\}} \right| \leq \int_{\overline{B}_n^b(r_2,\overline{x}_{*[1,n]})} \left| \psi^{\langle 1,n \rangle}(x_{*0}+\lambda r_2, \overline{x}_{[1,n]}) \right| ds_{\{n\}}$$

$$\leq \frac{C_{n,0}}{r_2^{p_{*n,0}}} \int_{\overline{B}_n^b(r_2,\overline{x}_{*[1,n]})} ds_{\{n\}} = \frac{C_{n,0} r_2^{n-1} S_n(1)}{r_2^{p_{*n,0}}},$$
(7.4₁)

because

$$\int_{\overline{B}_n^b(r_2,\overline{x}_{*[1,n]})} ds_{\{n\}} = S_n(r_2) = r_2^{n-1} S_n(1) \tag{7.4_2}$$

(cf. (5.11)). Hence, at $r_2 > r_*$ the second term on the left-hand side of identity (6.1) can be majorated (estimated by absolute value from above) thus:

$$\frac{1}{S_n(1) r_2^{n-1}} \left| \int_{\overline{B}_n^b(r_2,\overline{x}_{*[1,n]})} \psi^{\langle 1,n \rangle}(x_{*0}+\lambda r_2, \overline{x}_{[1,n]}) ds_{\{n\}} \right| \leq \frac{C_{n,0}}{r_2^{p_{*n,0}}}. \tag{7.4_3}$$

Since $p_{*n,0} > 0$ by Axiom 7.1, therefore (7.4₃) implies (7.4).

2) By (7.2) with $n \geq 2$ and $m=1$, it follows from (6.5) at $r_2 > r_*$ that

$$\left| \Lambda(\psi^{\langle 1,n \rangle}, \tau_n, \overline{x}_{[1,n]}) \right|$$
$$\leq \left[ \sum_{i=1}^n \left| \frac{x_i'}{r_2} \right| \left| \psi_{x_i}^{\langle 1,n \rangle}(x_0, \overline{x}_{[1,n]}) \right| + \left| \lambda \psi_{x_0}^{\langle 1,n \rangle}(x_0, \overline{x}_{[1,n]}) \right| \right]_{x_0 = x_{*0}+\lambda r_2}$$
$$\leq \left[ \sum_{i=1}^n \left| \psi_{x_i}^{\langle 1,n \rangle}(x_0, \overline{x}_{[1,n]}) \right| + \left| \psi_{x_0}^{\langle 1,n \rangle}(x_0, \overline{x}_{[1,n]}) \right| \right]_{x_0 = x_{*0}+\lambda r_2} \leq \frac{(n+1)C_{n,1}}{r_2^{p_{*n,1}}}.$$
(7.5₁)

Therefore,

$$\left| \int_{\overline{B}_n^b(r_2,\overline{x}_{*[1,n]})} \Lambda(\psi^{\langle 1,n \rangle}, \tau_n, \overline{x}_{[1,n]}) ds_{\{n\}} \right| \leq \int_{\overline{B}_n^b(r_2,\overline{x}_{*[1,n]})} \left| \Lambda(\psi^{\langle 1,n \rangle}, \tau_n, \overline{x}_{[1,n]}) \right| ds_{\{n\}}$$

$$\leq \frac{(n+1)C_{n,1}}{r_2^{p_{*n,1}}} \int_{\overline{B}_n^b(r_2,\overline{x}_{*[1,n]})} ds_{\{n\}} = \frac{(n+1)S_n(1)C_{n,1} r_2^{n-1}}{r_2^{p_{*n,1}}},$$
(7.5₂)

by (7.4₂). At the same time, by (4.9)–(4.12),

$$\left| \eta^{\langle n \rangle}(r_2) \right| = \frac{1}{(n-2)S_n(1) r_2^{n-2}} \quad \text{for each } n \in \omega_3, \tag{7.5_3a}$$

$$\left| \eta^{\langle 2 \rangle}(r_2) \right| = \frac{1}{S_2(1)} \ln r_2 \quad \text{for } n=2. \tag{7.5_3b}$$

Hence, the expression on the right-hand side of (6.4) at $r_2 > r_*$ is majorated thus:

$$\left|\Omega^{(3)}\left(\eta^{\langle n\rangle},\psi^{\langle 1,n\rangle},\tau_n,\overline{B}_n^{\mathrm{b}}\left(r_2,\overline{x}_{*[1,n]}\right)\right)\right|\leq\frac{(n+1)C_{n,1}}{(n-2)r_2^{p_{*n,1}-1}}\text{ for each }n\in\omega_3,\qquad(7.5_4\mathrm{a})$$

$$\left|\Omega^{(3)}\left(\eta^{\langle 2\rangle},\psi^{\langle 1,2\rangle},\tau_2,\overline{B}_2^{\mathrm{c}}\left(r_2,\overline{x}_{*[1,2]}\right)\right)\right|\leq\frac{2C_{2,1}\ln r_2}{r_2^{p_{*n,1}-1}}\text{ for }n=2.\qquad(7.5_4\mathrm{b})$$

Since $p_{*2,1}>1$ by Axiom 7.1, therefore it follows by l'Hospitale's rule that

$$\lim_{r_2\to 0}\frac{\ln r_2}{r_2^{p_{*n,1}-1}}=\lim_{r_2\to 0}\left[\frac{d\ln r_2}{dr_2}\left(\frac{dr_2^{p_{*n,1}-1}}{dr_2}\right)^{-1}\right]=\frac{1}{p_{*n,1}-1}\lim_{r_2\to 0}\left[\frac{1}{r_2}\frac{1}{r_2^{p_{*n,1}-2}}\right]$$
$$=\frac{1}{p_{*n,1}-1}\lim_{r_2\to 0}\frac{1}{r_2^{p_{*n,1}-1}}=0.\qquad(7.5_5)$$

Hence, (7.5$_4$a) and (7.5$_4$b) imply (7.5).

3) By (7.2) with $n=1$ and $m=0$, it follows that if $r_2>r_*$ then

$$\frac{1}{2}\left|\psi^{\langle 1,1\rangle}\left(x_{*0}+\lambda r_2,x_{*1}-r_2\right)+\psi^{\langle 1,1\rangle}\left(x_{*0}+\lambda r_2,x_{*1}+r_2\right)\right|\leq\frac{C_{1,0}}{r_2^{p_{*1,0}}},\qquad(7.6_1)$$

which implies (7.6) for $p_{*1,0}>0$, by Axiom 7.1.

4) By (7.2) with $n=1$ and $m=1$, it follows from (6.10) at $r_2>r_*$ that

$$\left|\Psi_1\left(x_{*0}+\lambda r_2,x_{1,}\right)\right|=\left|\psi_{x_1}^{\langle 1,1\rangle}\left(x_{*0}+\lambda r_2,x_{1,}\right)-\frac{\lambda x_1'}{r_2}\psi_{x_0}^{\langle 1,1\rangle}\left(x_{*0}+\lambda r_2,x_{1,}\right)\right|$$
$$\leq\left|\psi_{x_1}^{\langle 1,1\rangle}\left(x_{*0}+\lambda r_2,x_{1,}\right)\right|+\left|\frac{\lambda x_1'}{r_2}\psi_{x_0}^{\langle 1,1\rangle}\left(x_{*0}+\lambda r_2,x_{1,}\right)\right|\qquad(7.7_1)$$
$$\leq\left|\psi_{x_1}^{\langle 1,1\rangle}\left(x_{*0}+\lambda r_2,x_{1,}\right)\right|+\left|\psi_{x_0}^{\langle 1,1\rangle}\left(x_{*0}+\lambda r_2,x_{1,}\right)\right|\leq\frac{2C_{1,1}}{r_2^{p_{*1,1}}}.$$

Hence, the expression on the right-hand side of (6.9) is majorated thus:

$$\left|\Omega^{(3)}\left(\psi^{\langle 1,1\rangle},\eta^{\langle 1\rangle},\tau_1,[x_{*1}+r_2,x_{*1}+r_2]\right)\right|$$
$$\leq\frac{r_2}{2}\left[\left|\Psi_1\left(x_{*0}+\lambda r_2,x_{*1}-r_2\right)\right|+\left|\Psi_1\left(x_{*0}+\lambda r_2,x_{*1}+r_2\right)\right|\right]\leq\frac{2C_{1,1}}{r_2^{p_{*1,1}-1}}.\qquad(7.7_2)$$

Since $p_{*1,1}>1$ by Axiom 7.1, therefore (7.7$_2$) implies (7.7).•

**Lemma 7.3.** Given $\lambda\in\{-1,1\}$, for each $n\in\omega_1$, for each $\alpha\in\{0,1\}$,

$$\Omega^{(\alpha)}\left(\eta^{\langle n\rangle},\psi^{\langle 1,n\rangle},\tau_n,\overline{E}_n\right)=\lim_{r_2\to 0}\Omega^{(\alpha)}\left(\eta^{\langle n\rangle},\psi^{\langle 1,n\rangle},\tau_n,\overline{B}_n^{\mathrm{c}}\left(r_2,\overline{x}_{*[1,n]}\right)\right),\qquad(7.8)$$

the understanding being that the improper integral of third kind converges.

**Proof:** Assuming that $r_2>r_*$, represent $\Omega^{(\alpha)}\left(\eta^{\langle n\rangle},\psi^{\langle 1,n\rangle},\tau_n,\overline{B}_n^{\mathrm{c}}\left(r_2,\overline{x}_{*[1,n]}\right)\right)$, which is an improper integral of first kind, as:

$$\Omega^{(\alpha)}\left(\eta^{\langle n\rangle},\psi^{\langle 1,n\rangle},\tau_n,\overline{B}_n^c(r_2,\overline{x}_{*[1,n]})\right) = \Omega^{(\alpha)}\left(\eta^{\langle n\rangle},\psi^{\langle 1,n\rangle},\tau_n,\overline{B}_n^c(r_*,\overline{x}_{*[1,n]})\right)$$
$$+ \Omega^{(\alpha)}\left(\eta^{\langle n\rangle},\psi^{\langle 1,n\rangle},\tau_n,\overline{C}_n^{cc}(r_*,r_2,\overline{x}_{*[1,n]})\right),\quad (7.9)$$

the understanding being that the first integral on the right-hand side of this identity is a well-defined improper integral of the first kind, whereas the second one is an ordinary Riemann integral, which for $\alpha=0$ or $\alpha=1$ is respectively the variant of (4.72) or (4.73) with $r_*$ in place of $r_1$. More specifically,

$$\Omega^{(\alpha)}\left(\eta^{\langle n\rangle},\psi^{\langle 1,n\rangle},\tau_n,\overline{C}_n^{cc}(r_*,r_2,\overline{x}_{*[1,n]})\right) \quad (7.10)$$

is the variant with $r_*$ in place of $r_1$ of (4.72a) if $\alpha=0$ and $n\neq 2$, (4.72b) if $\alpha=0$ and $n=2$, (4.73a) if $\alpha=1$ and $n\neq 2$, or of (4.73b) if $\alpha=1$ and $n=2$. These four Riemann integrals can be majorated as follows.

1) By condition (7.3), it follows from (4.72) with $r_*$ in place of $r_1$ that

$$\left|\Omega^{(0)}\left(\eta^{\langle n\rangle},\psi^{\langle 1,n\rangle},\tau_n,\overline{C}_n^{cc}(r_*,r_2,\overline{x}_{*[1,n]})\right)\right|$$
$$= \left|\int_{\overline{C}_n^{cc}(r_*,r_2,\overline{x}_{*[1,n]})} \eta^{\langle n\rangle}(x'_{[1,n]}) f^{\langle 1,n\rangle}(x_{*0}+\lambda x'_{[1,n]},\overline{x}_{[1,n]})\,dv_{\{n\}}\right|$$
$$\leq \int_{\overline{C}_n^{cc}(r_*,r_2,\overline{x}_{*[1,n]})} \left|\eta^{\langle n\rangle}(x'_{[1,n]})\right|\left|f^{\langle 1,n\rangle}(x_{*0}+\lambda x'_{[1,n]},\overline{x}_{[1,n]})\right|dv_{\{n\}} \quad (7.11)$$
$$\leq (n+1)C_{n,2}\int_{\overline{C}_n^{cc}(r_*,r_2,\overline{x}_{*[1,n]})} {x'_{[1,n]}}^{-p_{*n,2}}\left|\eta^{\langle n\rangle}(x'_{[1,n]})\right|dv_{\{n\}}.$$

In this case,

$$\int_{\overline{C}_n^{cc}(r_*,r_2,\overline{x}_{*[1,n]})} {x'_{[1,n]}}^{-p_{*n,2}}\left|\eta^{\langle n\rangle}(x'_{[1,n]})\right|dv_{\{n\}} = \int_{r_*}^{r_2} dr\, r^{-p_{*n,2}}\left|\eta^{\langle n\rangle}(r)\right|\int_{\overline{B}_{\{n\}}^b(r,\overline{x}_{*[1,n]})} ds_{\{n\}}$$
$$= S_n(1)\int_{r_*}^{r_2} r^{n-1-p_{*n,2}}\left|\eta^{\langle n\rangle}(r)\right|dr, \quad (7.11_1)$$

because

$$\int_{\overline{B}_n^b(r,\overline{x}_{*[1,n]})} ds_{\{n\}} = S_n(r) = r^{n-1}S_n(1). \quad (7.11_2)$$

(cf. (5.10) and (5.11)). Hence, the train of inequalities (7.11) reduces to

$$\left|\Omega^{(0)}\left(\eta^{\langle n\rangle},\psi^{\langle 1,n\rangle},\tau_n,\overline{C}_n^{cc}(r_*,r_2,\overline{x}_{*[1,n]})\right)\right|$$
$$\leq (n+1)C_{n,2}S_n(1)\int_{r_*}^{r_2} r^{n-1-p_{*n,2}}\left|\eta^{\langle n\rangle}(r)\right|dr. \quad (7.12)$$

2) Similarly, by condition (7.2) with $m=1$ and $i=0$, , it follows from (4.73) with $r_*$ in place of $r_1$ that

$$\left|\Omega^{(1)}\left(\eta^{\langle n\rangle},\psi^{\langle 1,n\rangle},\tau_n,\overline{C}_n^{cc}\left(r_*,r_2,\overline{x}_{*[1,n]}\right)\right)\right|$$
$$=\left|\lambda(n-3)\int_{\overline{C}_n^{cc}\left(r_*,r_2,\overline{x}_{*[1,n]}\right)}\zeta^{\langle n\rangle}\left(x'_{[1,n]}\right)\psi_{x_0}^{\langle 1,n\rangle}\left(x_{*0}+\lambda x'_{[1,n]},\overline{x}_{[1,n]}\right)dv_{\{n\}}\right|$$
$$\leq|n-3|\int_{\overline{C}_n^{cc}\left(r_*,r_2,\overline{x}_{*[1,n]}\right)}\left|\zeta^{\langle n\rangle}\left(x'_{[1,n]}\right)\right|\left|\psi_{x_0}^{\langle 1,n\rangle}\left(x_{*0}+\lambda x'_{[1,n]},\overline{x}_{[1,n]}\right)\right|dv_{\{n\}}$$
$$\leq|n-3|C_{n,1}\int_{\overline{C}_n^{cc}\left(r_*,r_2,\overline{x}_{*[1,n]}\right)}\left|x'_{[1,n]}\right|^{-p_{*n,1}}\left|\zeta^{\langle n\rangle}\left(x'_{[1,n]}\right)\right|dv_{\{n\}} \quad (7.13)$$

In this case, in analogy with (7.11),

$$\int_{\overline{C}_n^{cc}\left(r_*,r_2,\overline{x}_{*[1,n]}\right)}\left|x'_{[1,n]}\right|^{-p_{*n,1}}\left|\zeta^{\langle n\rangle}\left(x'_{[1,n]}\right)\right|dv_{\{n\}}=\int_{r_*}^{r_2}drr^{-p_{*n,1}}\left|\zeta^{\langle n\rangle}(r)\right|\int_{B_{\{n\}}^b\left(r,\overline{x}_{*[1,n]}\right)}ds_{\{n\}}$$
$$=S_n(1)\int_{r_*}^{r_2}r^{n-1-p_{*n,1}}\left|\zeta^{\langle n\rangle}(r)\right|dr, \quad (7.13_1)$$

so that (7.13) reduces to

$$\left|\Omega^{(1)}\left(\eta^{\langle n\rangle},\psi^{\langle 1,n\rangle},\tau_n,\overline{C}_n^{cc}\left(r_*,r_2,\overline{x}_{*[1,n]}\right)\right)\right|$$
$$\leq|n-3|C_{n,1}S_n(1)\int_{r_*}^{r_2}r^{n-1-p_{*n,1}}\left|\zeta^{\langle n\rangle}(r)\right|dr. \quad (7.14)$$

3) By (4.60a) (see also (4.9) subject to (4.10) and (4.13) subject to (4.14)) and by (4.74a), the integrals occurring on the right-hand sides of inequalities (7.12) and (7.14) can be calculated for each $n\in\omega_1-\{2\}$ thus:

$$\left|a_n\right|^{-1}\int_{r_*}^{r_2}r^{n-1-p_{*n,2}}\left|\eta^{\langle n\rangle}(r)\right|dr=\int_{r_*}^{r_2}r^{-(p_{*n,2}-1)}dr=-\frac{1}{p_{*n,2}-2}r^{-(p_{*n,2}-2)}\Big|_{r_*}^{r_2}$$
$$=\frac{1}{p_{*n,2}-2}\left(r_*^{-(p_{*n,2}-2)}-r_2^{-(p_{*n,2}-2)}\right)<\frac{1}{p_{*n,2}-2}r_*^{-(p_{*n,2}-2)}, \quad (7.15)$$

$$\left|a_n\right|^{-1}\int_{r_*}^{r_2}r^{n-1-p_{*n,1}}\left|\zeta^{\langle n\rangle}(r)\right|dr=\int_{r_*}^{r_2}r^{-p_{*n,1}}dr=-\frac{1}{p_{*n,2}-1}r^{-(p_{*n,1}-1)}\Big|_{r_*}^{r_2}$$
$$=\frac{1}{p_{*n,1}-1}\left(r_*^{-(p_{*n,1}-1)}-r_2^{-(p_{*n,1}-1)}\right)<\frac{1}{p_{*n,1}-1}r_*^{-(p_{*n,1}-1)}. \quad (7.16)$$

Hence, for each $n\in\omega_1-\{2\}$, inequalities (7.12) and (7.14) imply the following ones:

$$\left|\Omega^{(0)}\left(\eta^{\langle n\rangle},\psi^{\langle 1,n\rangle},\tau_n,\overline{C}_n^{cc}\left(r_*,r_2,\overline{x}_{*[1,n]}\right)\right)\right|$$

$$< \frac{(n+1)C_{n,2}S_n(1)|a_n|}{(p_{*n,2}-2)r_*^{p_{*n,2}-2}} = \frac{(n+1)C_{n,2}}{|n-2|(p_{*n,2}-2)r_*^{p_{*n,2}-2}}, \qquad (7.17)$$

$$\left|\Omega^{(1)}\left(\eta^{\langle n\rangle},\psi^{\langle 1,n\rangle},\tau_n,\overline{C}_n^{cc}\left(r_*,r_2,\overline{x}_{*[1,n]}\right)\right)\right|$$

$$< \frac{|n-3|C_{n,1}S_n(1)|a_n|}{(p_{*n,1}-1)r_*^{p_{*n,1}-1}} = \frac{|n-3|C_{n,1}}{|n-2|(p_{*n,1}-1)r_*^{p_{*n,1}-1}}. \qquad (7.18)$$

4) By (4.60b) (see also (4.11) subject to (4.12)) and by (4.74b), the integrals occurring on the right-hand sides of inequalities (7.12) and (7.14) can be calculated for $n=2$ by integration by part as follows:

$$2\pi \int_{r_*}^{r_2} r^{1-p_{*2,2}} \left|\eta^{\langle 2\rangle}(r)\right| dr = \int_{r_*}^{r_2} r^{-(p_{*2,2}-1)} \ln r\, dr$$

$$= -\frac{1}{p_{*2,2}-2}\left[ r^{-(p_{*2,2}-2)}\ln r \Big|_{r_*}^{r_2} - \int_{r_*}^{r_2} r^{-(p_{*2,2}-1)} dr \right]$$

$$= -\frac{1}{p_{*2,2}-2}\left[ r^{-(p_{*2,2}-2)}\ln r \Big|_{r_*}^{r_2} + \frac{1}{p_{*2,2}-2} r^{-(p_{*2,2}-2)}\Big|_{r_*}^{r_2} \right] \qquad (7.19)$$

$$= \frac{1}{p_{*2,2}-2}\left[ r_*^{-(p_{*2,2}-2)}\left(\ln r_* + \frac{1}{p_{*2,2}-2}\right) \right]$$

$$- \frac{1}{p_{*2,2}-2}\left[ r_2^{-(p_{*2,2}-2)}\left(\ln r_2 + \frac{1}{p_{*2,2}-2}\right) \right] > 0$$

$$2\pi \int_{r_*}^{r_2} r^{-(p_{*2,1}-1)} \left|\zeta^{\langle 2\rangle}(r)\right| dr = \int_{r_*}^{r_2} r^{-p_{*2,1}}(\ln r + 2) dr$$

$$= -\frac{1}{p_{*2,1}-1}\left[ r^{-(p_{*2,1}-1)}(\ln r + 2)\Big|_{r_*}^{r_2} - \int_{r_*}^{r_2} r^{-p_{*2,1}} dr \right]$$

$$= -\frac{1}{p_{*2,1}-1}\left[ r^{-(p_{*2,1}-1)}(\ln r + 2)\Big|_{r_*}^{r_2} + \frac{1}{p_{*2,1}-1} r^{-(p_{*2,1}-1)}\Big|_{r_*}^{r_2} \right] \qquad (7.20)$$

$$= \frac{1}{p_{*2,1}-1}\left[ r_*^{-(p_{*2,1}-1)}\left(\ln r_* + \frac{2p_{*2,1}-1}{p_{*2,1}-1}\right) \right]$$

$$- \frac{1}{p_{*2,1}-1}\left[ r_2^{-(p_{*2,1}-1)}\left(\ln r_2 + \frac{2p_{*2,1}-1}{p_{*2,1}-1}\right) \right] > 0.$$

Hence, for $n=2$, inequalities (7.12) and (7.14) imply the following ones:

$$\left|\Omega^{(0)}\left(\eta^{\langle 2\rangle},\psi^{\langle 1,2\rangle},\tau_2,\overline{C}_2^{cc}\left(r_*,r_2,\overline{x}_{*[1,2]}\right)\right)\right|$$

$$< \frac{2\,C_{,2}}{(p_{*2,2}-2)r_*^{p_{*2,2}-2}}\left(\ln r_* + \frac{1}{p_{*2,2}-2}\right), \tag{7.21}$$

$$\left|\Omega^{(1)}\left(\eta^{\langle 2\rangle},\psi^{\langle 1,2\rangle},\tau_2,\overline{C}_2^{cc}\left(r_*,r_2,\overline{x}_{*[1,2]}\right)\right)\right|$$

$$< \frac{C_{2,1}}{(p_{*2,1}-1)r_*^{p_{*2,1}-1}}\left(\ln r_* + \frac{2p_{*2,1}-1}{p_{*2,1}-1}\right). \tag{7.22}$$

Inequalities (7.17), (7.18), (7.21), and (7.12) mean that given $\lambda \in \{-1,1\}$, for each $n \in \omega_1$, for each $\alpha \in \{0,1\}$, integral (7.10), and hence (7.9) converges. QED.•

**Theorem 7.1.** Given $\lambda \in \{-1,1\}$, given $n \in \omega_1$, as $r_2 \to \infty$, the SIDF (6.1) reduces to:

$$\psi^{\langle 1,n\rangle}\left(x_{*0},\overline{x}_{*[1,n]}\right) = -\sum_{\alpha \in \{0,1\}}\Omega^{(\alpha)}\left(\eta^{\langle n\rangle},\psi^{\langle 1,n\rangle},\tau_n,\overline{E}_n\right) \tag{7.23}$$

subject to

$$\Omega^{(0)}\left(\eta^{\langle n\rangle},\psi^{\langle 1,n\rangle},\tau_n,\overline{E}_n\right) = -\int_{\overline{E}_n}\eta^{\langle n\rangle}\left(x'_{[1,n]}\right)f^{\langle 1,n\rangle}\left(x_{*0}+\lambda x'_{[1,n]},\overline{x}_{[1,n]}\right)dv_{\{n\}}, \tag{7.24}$$

$$\Omega^{(1)}\left(\eta^{\langle n\rangle},\psi^{\langle 1,n\rangle},\tau_n,\overline{E}_n\right)$$
$$= \lambda(n-3)\int_{\overline{E}_n}\zeta^{\langle n\rangle}\left(x'_{[1,n]}\right)\psi_{x_0}^{\langle 1,n\rangle}\left(x_{*0}+\lambda x'_{[1,n]},\overline{x}_{[1,n]}\right)dv_{\{n\}}, \tag{7.25}$$

where $\tau_n$ is the list of parameters $\overline{x}_{*0},\lambda,\overline{x}_{*[1,n]}$

$$x'_{[1,n]} \equiv \left|\overline{x}_{[1,n]} - \overline{x}_{*[1,n]}\right| = \sqrt{\sum_{i=1}^{n}(x_i - x_{*i})^2} \geq 0, \tag{7.26}$$

$$\eta^{\langle n\rangle}\left(x'_{[1,n]}\right) \equiv \begin{cases} a_n x'^{2-n}_{[1,n]} & \text{if } n \neq 2 \text{ (a)} \\ -a_2 \ln x'_{[1,2]} & \text{if } n = 2 \text{ (b)}, \\ a_1|x'_1| & \text{if } n = 1 \text{ (c)} \end{cases} \tag{7.27}$$

$$\zeta^{\langle n\rangle}\left(x'_{[1,n]}\right) \equiv \begin{cases} a_n x'^{1-n}_{[1,n]} & \text{if } n \neq 2 \text{ (a)} \\ -\dfrac{a_2(\ln x'_{[1,2]}+2)}{x'_{[1,2]}} & \text{if } n = 2 \text{ (b)}, \\ a_1 & \text{if } n \neq 2 \text{ (c)} \end{cases} \tag{7.28}$$

$$a_n \equiv \begin{cases} \dfrac{1}{(n-2)S_n(1)} = \dfrac{\Gamma(n/2)}{2(n-2)\pi^{n/2}} & \text{if } n \geq 3 \text{ (a)} \\ \dfrac{1}{S_2(1)} = \dfrac{1}{2\pi} & \text{if } n = 2 \text{ (b)}, \\ -\dfrac{1}{2} & \text{if } n = 1 \text{ (c)} \end{cases} \tag{7.29}$$

so that particularly,

$$a_{\{3\}} = \frac{1}{4\pi}, \ a_{\{4\}} = \frac{1}{4\pi^2}, \ a_{\{5\}} = \frac{1}{8\pi^2}, \text{ etc.} \quad (7.29a)$$

**Proof:** The theorem follows from Theorem 6.1 by Lemmas 7.2 and 7.3. •

**Corollary 7.1.** For $n=3$, given $\lambda \in \{-1,1\}$, given $\bar{x}_{*0} \in \bar{E}_{\{1\}}$, given $\bar{x}_{*[1,3]} \in \bar{E}_{\{3\}}$:

$$\psi^{\langle 1,3 \rangle}(x_{*0}, \bar{x}_{*[1,3]}) = \frac{1}{4\pi} \int_{\bar{E}_3} \frac{f^{\langle 1,3 \rangle}(x_{*0} + \lambda|\bar{x}_{[1,3]} - \bar{x}_{*[1,3]}|, \bar{x}_{[1,3]})}{|\bar{x}_{[1,3]} - \bar{x}_{*[1,3]}|} dx_1 dx_2 dx_3 \quad (7.30)$$

subject to

$$\bar{x}_{[1,3]} \triangleq \langle x_1, x_2, x_3 \rangle \in \bar{E}_3, \ \bar{x}_{*[1,3]} \triangleq \langle x_{*1}, x_{*2}, x_{*3} \rangle \in \bar{E}_3, \quad (7.31)$$

$$|\bar{x}_{[1,3]} - \bar{x}_{*[1,3]}| \triangleq \sqrt{\sum_{i=1}^{3}(x_i - x_{*i})^2} \geq 0. \quad (7.32)$$

**Proof:** The corollary is a detailed instance of Theorem 7.1 for $n=3$.•

### 7.3. Concluding comments

**1.** Given $\lambda \in \{-1,1\}$, given $n \in \omega_1$, if the function $\psi^{\langle 1,n \rangle}$, i.e. actually $\psi^{\lambda\langle 1,n \rangle}$, is supposed to be *determined* by the source function $f^{\langle 1,n \rangle}$, i.e. actually $f^{\lambda\langle 1,n \rangle}$, with the help of the *inhomogeneous differential equation* (4.75) subject to (3.12) in the whole of the direct product $R \times \bar{E}_n$ then $\psi^{\langle 1,n \rangle}$ satisfies the integral formula (7.23) subject to (7.24)–(7.29). However, a *non-dispersive* retarded or advanced solution for $\psi^{\langle 1,n \rangle}$, which is given by the equation (7.30) subject to (7.31) and (7.32), exists only in the case of $n=3$ under certain asymptotic conditions at infinity as specified by the pertinent instance of Axiom 7.1. The solution is *retarded* if $\lambda=-1$ and *advanced* if $\lambda=1$. The advanced solution is impracticable in the sense that it has no interpretations by physical phenomena. By contrast, the retarded solution expresses or just is *the non-dispersive principle of causality*, which is one of the most fundamental principles of nature. Therefore, Axiom 7.1 should be restricted to the case of $n=3$ and be supplemented by the following axiom.

**Axiom 7.2:** ***The axiom of non-dispersive causality.*** In the case of $n=3$, the functions $\psi^{-1\langle 1,3 \rangle}$ and $f^{-1\langle 1,3 \rangle}$ retain, whereas the functions $\psi^{1\langle 1,3 \rangle}$ and $f^{1\langle 1,3 \rangle}$, are disregarded, i.e. *the ranges of the functional placeholders (functional place-holding variables)* '$\psi^{1\langle 1,3 \rangle}$' and '$f^{1\langle 1,3 \rangle}$' are empty.•

The differential properties of formula (7.23), which make its instance at *n*=3 exclusive, are diacussed in Appendix C.

**2.** The electromagnetic field of an electric 4-current of density $\langle c\rho, \mathbf{j} \rangle$ (see, e.g. Landau & Lifshitz [1989, p. 71]), as an interpretand of $f^{\langle 1,3 \rangle}$, at large distances from the latter, i.e. the electromagnetic waves radiated by the 4-current, is described by a relativistically invariant retarded 4-vector potential $\langle \phi, \mathbf{A} \rangle$ (*ibid.* pp. 45, 170) of the field, as an interpretand of $\psi^{\langle 1,3 \rangle}$, which is the retarded solution of the inhomogeneous d'Alambertian (wave) equation (*ibid.* p. 158) (as (4.75) at *n*=3). From the standpoint of algebraic analysis, the transformation properties of a pseudo-Euclidean vector space if index 1 are largely same independent on the dimension of the basic Euclidean space. However, the boundary-free SIDF as given by Theorem 7.1 is not a solution of the inhomogeneous d'Alambertian (wave) equation exists in an infinite affine Euclidean space of any dimension other than 3. Therefore, *a three-dimensional affine Euclidean space is the only possible receptacle of nature*. For a detailed discussion of this fact, see section 1.•

## Appendix A: The volume of an *n*-dimensional sphere and the area of its surface

Given $n \in \omega_2$, the *n*-dimensional volume [measure] $dV_n(r)$ of a thin *n*-dimensional spherical layer $\overline{C}_n^{cc}(r+dr, \overline{x}_{*[1,n]})$ is related to the (*n*–1)-dimensional surface area $S_n(r)$ [measure] of the (*n*–1)-dimensional boundary spherical surface $\overline{B}_n^b(r, \overline{x}_{*[1,n]})$ as $dV_n(r) = S_n(r)dr$, whence $S_n(r) = dV_n(r)/dr$. Therefore, we proceed from calculation of $V_n(r)$, which can be defined as

$$V_n(r) \stackrel{=}{=} \int\limits_{|\overline{x}_{[1,n]} - \overline{x}_{*[1,n]}| \leq r} dx_1 \cdots dx_n \tag{A.1}$$

subject to (2.44), where $dx_i \geq 0$ for each $i \in \omega_{1,n}$. With the help of the substitution

$$\overline{\xi}_{[1,n]} = \langle \xi_1, \ldots, \xi_n \rangle = r^{-1} \langle x_1 - x_{*1}, \ldots, x_n - x_{*n} \rangle = r^{-1} (\overline{x}_{[1,n]} - \overline{x}_{*[1,n]}), \tag{A.2}$$

(A.1) becomes

$$V_n(r) = r^n V_n(1), \tag{A.3}$$

where $V_n(1)$, defined as

$$V_n(1) \equiv \int_{|\bar{\xi}_{[1,n]}|\leq 1} d\xi_1 \cdots d\xi_n \tag{A.4}$$

subject to

$$\xi_{[1,n]} \equiv |\xi_{[1,n]}| = \sqrt{\sum_{i=1}^{n} \xi_i^2} \geq 0, \tag{A.5}$$

is the volume of the sphere $\bar{B}_n^c(1, \bar{x}_{*[1,n]})$ of unit radius. Consequently,

$$S_n(r) = \frac{dV_n(r)}{dr} = nr^{n-1}V_n(1) = r^{n-1}S_n(1), \tag{A.6}$$

where $S_n(1)$, defined as

$$S_n(1) \equiv nV_n(1), \tag{A.7}$$

is the area of the spherical surface $\bar{B}_n^b(1, \bar{x}_{*[1,n]})$ of unit radius.

In order to calculate $V_n(1)$ conveniently, we shall use the beta and gamma functions, which are conventionally defined as:

$$B(z_1, z_2) \equiv \int_0^1 u^{z_1-1}(1-u)^{z_2-1} \, du \tag{A.8}$$

$$\Gamma(z) \equiv \int_0^\infty e^{-t} t^{z-1} \, dt \tag{A.9}$$

for all complex numbers $z_1$, $z_2$, and $z$ such that $\operatorname{Re} z_1 > 0$, $\operatorname{Re} z_2 > 0$, and $\operatorname{Re} z > 0$. The integrals, being the definientia of definitions (A.8) and (A.9), are said to be the *Eulerian integral of the first kind* and the *Eulerian integral of the second kind*, while the functions B and $\Gamma$ thus defined are said to be the *beta-function* and the *gamma-function* respectively. It is known that

$$B(z_1, z_2) = \frac{\Gamma(z_1)\Gamma(z_2)}{\Gamma(z_1+z_2)} \tag{A.10}$$

(see, e.g. Spiegel [1964, p. 282)] and also that

$$\Gamma(z+1) = z\Gamma(z), \tag{A.11}$$

which immediately follows from (A.9) by integration by parts:

$$\Gamma(z) = \int_0^\infty e^{-t} t^{z-1} \, dt = \frac{1}{z} e^{-t} t^z \Big|_0^\infty + \frac{1}{z}\int_0^\infty e^{-t} t^z \, dt = \frac{1}{z}\Gamma(z+1); \tag{A.11_1}$$

A generalized version of (A.10) is proved in Comment A.1 below in this appendix.

**Lemma A.1.** For each $n \in \omega_2$:

$$V_n(1) = \frac{\sqrt{\pi}\,\Gamma\!\left(\dfrac{n+1}{2}\right)}{\Gamma\!\left(\dfrac{n+2}{2}\right)} V_{n-1}(1). \tag{A.12}$$

**Proof:** By (A.5), the inequalities $\left|\bar{\xi}_{[1,n]}\right| \leq 1$ and $\left|\bar{\xi}_{[1,n]}\right|^2 \leq 1$ are equivalent. The latter can be rewritten as

$$\left|\xi_{[1,n-1]}\right|^2 = 1 - \xi_n^{\,2}. \tag{A.13}$$

where

$$\left|\bar{\xi}_{[1,n-1]}\right|^2 = \sum_{i=1}^{n-1} \xi_i^{\,2} > 0. \tag{A.14}$$

Hence, the integral (A.4) can be developed thus:

$$V_n(1) = \int_{\left|\bar{\xi}_{[1,n]}\right|^2 \leq 1} d\xi_1 \cdots d\xi_n = \int_{-1}^{1} d\xi_n \int_{\left|\bar{\xi}_{[1,n-1]}\right|^2 \leq 1 - \xi_n^{\,2}} d\xi_1 \cdots d\xi_{n-1} = \int_{-1}^{1} V_{n-1}\!\left(\sqrt{1 - \xi_n^{\,2}}\right) d\xi_n, \tag{A.15}$$

where $V_{n-1}\!\left(\sqrt{1-\xi_n^{\,2}}\right)$ is the (*n*-1)-dimensional volume [measure] of an (*n*–1)-dimensional sphere of radius $\sqrt{1-\xi_n^{\,2}}$. Thus, by (A.3) with '*n* – 1' in place of '*n*' and with $r = \sqrt{1-\xi_n^{\,2}}$, equation (A.15) becomes

$$V_n(1) = c_{n-1} V_{n-1}(1), \tag{A.16}$$

where

$$c_{n-1} \equiv \int_{-1}^{1} \left(1 - \xi_n^{\,2}\right)^{\frac{n-1}{2}} d\xi_n = 2\int_{0}^{1} \left(1 - \xi_n^{\,2}\right)^{\frac{n-1}{2}} d\xi_n = \int_{0}^{1} u^{-\frac{1}{2}} (1-u)^{\frac{n-1}{2}} du$$
$$= B\!\left(\frac{1}{2}, \frac{n+1}{2}\right). \tag{A.17}$$

In developing (A.17), use of the substitution $u = \xi_n^{\,2}$ and also use of the pertinent instance of definition (A.8) have been made. By (A.10), equation (A.17) becomes

$$c_{n-1} = \frac{\Gamma\!\left(\dfrac{1}{2}\right)\Gamma\!\left(\dfrac{n+1}{2}\right)}{\Gamma\!\left(\dfrac{n+2}{2}\right)} = \frac{\sqrt{\pi}\,\Gamma\!\left(\dfrac{n+1}{2}\right)}{\Gamma\!\left(\dfrac{n+2}{2}\right)}, \tag{A.18}$$

because (see, e.g., Spiegel [1964, p. 281, Problem 14])

$$\Gamma(1/2) = \sqrt{\pi}. \tag{A.19}$$

By (A.18) and (A.19), equation (A.16) coincides with (A.12). QED.•

**Theorem A.1.** For each $n \in \omega_1$:

$$V_n(1) = \frac{\pi^{n/2}}{\Gamma\left(\frac{n}{2}+1\right)} = \frac{2\pi^{n/2}}{n\Gamma(n/2)}, \tag{A.20}$$

$$S_n(1) = nV_n(1) = \frac{2\pi^{n/2}}{\Gamma(n/2)}, \tag{A.21}$$

**Proof:** It follows from (A.4) at $n=1$ that

$$V_1(1) = \int_{|\bar{\xi}_{[1,1]}|\le 1} d\xi_1 = \int_{|\xi_1|\le 1} d\xi_1 = \int_{-1}^{1} d\xi_1 = 2. \tag{A.22}$$

By (A.11), (A.19), and (A.22), equation (A.12) yields

$$V_n(1) = \pi^{\frac{n}{2}-1} V_1(1) \prod_{i=2}^{n} \frac{\Gamma\left(\frac{i+1}{2}\right)}{\Gamma\left(\frac{i+2}{2}\right)} = \frac{2\pi^{\frac{n-1}{2}}\Gamma\left(\frac{3}{2}\right)}{\Gamma\left(\frac{n+2}{2}\right)} = \frac{\pi^{\frac{n}{2}}}{\Gamma\left(\frac{n+2}{2}\right)}, \tag{A.23}$$

which proves (A.20). Identity (A.21) immediately follows from (A.7) by (A.20).•

**Corollary A.1.** 1) By (A.11) and (A.19), it follows from (A.20) and (A.21) that for each $m \in \omega_1$:

$$S_{2m}(1) = 2mV_{2m}(1) = \frac{2\pi^m}{(m-1)!}, \; S_{2m+1}(1) = (2m+1)V_{2m+1}(1) = \frac{2(2\pi)^m}{(2m-1)!!}, \tag{A.24}$$

whence for the first three values of '$m$':

$$S_2(1) = 2V_2(1) = 2\pi, \; S_3(1) = 3V_3(1) = 4\pi, \; S_4(1) = 4V_4(1) = 2\pi^2,$$
$$S_5(1) = 5V_5(1) = \frac{8\pi^2}{3}. \tag{A.24a}$$

2) From (A.20) it follows that for each $n \in \omega_2$:

$$S_{n+1}(1) = \frac{2\pi^{(n+1)/2}}{\Gamma(n/2)} = \frac{(2\pi)\pi^{(n-1)/2}}{\Gamma(n/2)} = 2\pi V_{n-1}(1), \tag{A.25}$$

whence

$$S_{2m+1}(1) = 2\pi V_{2m-1}(1), \; S_{2m}(1) = 2\pi V_{2(m-1)}(1) \tag{A.26}$$

subject to (A.24).•

**Comment A.1.** Identity (A.10) is an instance at $n=2$ of the following identity that is valid for any $n \in \omega_2$ complex numbers $z_1$, $z_2$, …, $z_n$, such that $\operatorname{Re} z_i > 0$ for each $i \in \omega_{1,n}$:

$$\mathrm{B}_n(z_1,\ldots,z_n) = \frac{\Gamma(z_1)\cdots\Gamma(z_n)}{\Gamma(z_1+\ldots+z_n)} = \frac{\prod_{i=1}^{n}\Gamma(z_i)}{\Gamma\left(\sum_{i=1}^{n}z_i\right)} \tag{A.27}$$

subject to

$$\mathrm{B}_n(z_1,\ldots,z_n) \doteq 2^{n-1}\int_{\substack{\overline{\xi}_{[1,n-1]}\leq 1,\\ \xi_1\geq 0,\ldots,\xi_{n-1}\geq 0}} \xi_1^{2z_1-1}\cdots\xi_{n-1}^{2z_{n-1}-1}\left(1-\left|\overline{\xi}_{[1,n-1]}\right|^2\right)^{z_n-1} d\xi_1\cdots d\xi_{n-1}, \tag{A.28}$$

so that $\mathrm{B}_2=\mathrm{B}$. Indeed, with the help of the substitution $t=u^2$, definition (A.9) becomes

$$\Gamma(z) = 2\int_0^\infty e^{-u^2} u^{2z-1}\, du. \tag{A.29}$$

Hence,

$$\begin{aligned}\Gamma(z_1)\cdots\Gamma(z_n) &= 2^n \int_0^\infty\cdots\int_0^\infty e^{-(u_1^2+\ldots+u_n^2)} u_1^{2z_1-1}\cdots u_n^{2z_n-1}\, du_1\cdots du_n \\ &= 2^n \int_{u_1\geq 0}\cdots\int_{u_n\geq 0} e^{-(u_1^2+\ldots+u_n^2)} u_1^{2z_1-1}\cdots u_n^{2z_n-1}\, du_1\cdots du_n.\end{aligned} \tag{A.30}$$

Instead of '$u_n$', let us introduce a new variable of integration '$v$' such as:

$$u_n = \sqrt{v^2 - u_1^2 - \ldots - u_{n-1}^2} \quad \text{subject to } u_1^2+\ldots+u_{n-1}^2 \leq v^2. \tag{A.31}$$

Thus, if the values of '$u_1$', ..., '$u_{n-1}$' are held constant then

$$du_n = \frac{v\, dv}{\sqrt{v^2 - u_1^2 - \ldots - u_{n-1}^2}}. \tag{A.32}$$

Therefore, (A.30) can be written as

$$\Gamma(z_1)\cdots\Gamma(z_n) = 2^n \int_0^\infty dv\, e^{-v^2} v \int_{\substack{u_1^2+\ldots+u_{n-1}^2\leq v^2,\\ u_1\geq 0,\ldots,u_{n-1}\geq 0}} u_1^{2z_1-1}\cdots u_{n-1}^{2z_{n-1}-1}\left(v^2 - u_1^2 - \ldots - u_{n-1}^2\right)^{z_n-1} du_1\cdots du_{n-1} \tag{A.33}$$

With the help of the substitutions

$$u_i = v\xi_i \quad \text{for each } i \in \omega_{1,n-1}, \tag{A.34}$$

equation (A.33) becomes

$$\begin{aligned}\Gamma(z_1)\cdots\Gamma(z_n) &= \left[2\int_0^\infty e^{-v^2} v^{2(z_1+\ldots+z_n)-1}\, dv\right] \\ &\quad \cdot 2^{n-1}\int_{\substack{\left|\overline{\xi}_{n-1}\right|^2\leq 1,\\ \xi_1\geq 0,\ldots,\xi_{n-1}\geq 0}} \xi_1^{2z_1-1}\cdots\xi_{n-1}^{2z_{n-1}-1}\left(1-\left|\overline{\xi}_{n-1}\right|^2\right)^{z_n-1} d\xi_1\cdots d\xi_{n-1}.\end{aligned} \tag{A.35}$$

By the pertinent variant of (A.29),

$$2\int_0^\infty e^{-v^2} v^{2(z_1+\ldots+z_n)-1}\, dv = \Gamma(z_1+\ldots+z_n). \tag{A.36}$$

Hence, by (A.28) and (A.36), equation (A.35) reduces to (A.29).•

# Appendix B: Normalization of the function $\eta^{\langle n \rangle}$

**Theorem B.1.** Given $n \in \omega_1$, given $\bar{x}_{*[1,n]} \in \bar{E}_n$, for each $r \in (0,\infty)$:

$$\int_{\bar{B}^{b+}_{\{n\}}(r,\bar{x}_{*[1,n]})} \sum_{i=1}^n \frac{\partial \eta^{\langle n \rangle}(x'_{[1,n]})}{\partial x_i} n_i(\bar{x}_{[1,n]})\, ds_{\{n\}} = -1 \text{ if } n \in \omega_2, \tag{B.1}$$

$$\left. \frac{\partial \eta^{\langle 1 \rangle}(x'_{[1,1]})}{\partial x_1} \right|_{x_{1,}=x_{*1}-r}^{x_{1,}=x_{*1}+r} = -1 \text{ if } n=1, \tag{B.2}$$

subject to Definitions 2.12 and 4.2.

**Proof:** By Definitions 2.12 and 4.2, the expression on the left-hand side of equation (B.1), or (B.2), can be developed respectively thus:

$$\int_{\bar{B}^{b+}_{\{n\}}(r,\bar{x}_{*[1,n]})} \sum_{i=1}^n \frac{\partial \eta^{\langle n \rangle}(x'_{[1,n]})}{\partial x_i} n_i(\bar{x}_{[1,n]})\, ds_{\{n\}}$$

$$= a_n(2-n) \int_{\bar{B}^{c-}_{\{n\}}(r,\bar{x}_{*[1,n]})} \sum_{i=1}^n \frac{x_i - x_{*i}}{r^n} \frac{x_i - x_{*i}}{r}\, ds_{\{n\}} \tag{B.3}$$

$$= a_n(2-n)r^{1-n} \int_{\bar{B}^{c-}_{\{n\}}(r,\bar{x}_{*[1,n]})} ds_{\{n\}} = a_n(2-n)r^{1-n}r^{n-1}S_{\{n\}}(1) = -1 \text{ if } n \geq 3,$$

$$\int_{\bar{B}^{b+}_{\{2\}}(r,\bar{x}_{*[1,2]})} \sum_{i=1}^2 \frac{\partial \eta^{\langle 2 \rangle}(x'_{[1,2]})}{\partial x_i} n_i(\bar{x}_{[1,2]})\, ds_{\{2\}}$$

$$= -a_2 \int_{\bar{B}^{b+}_{\{2\}}(r,\bar{x}_{*[1,2]})} \sum_{i=1}^n \frac{x_i - x_{*i}}{r^2} \frac{x_i - x_{*i}}{r}\, ds_{\{2\}} \tag{B.4}$$

$$= -a_2 r^{-1} \int_{\bar{B}^{b+}_{\{2\}}(r,\bar{x}_{*[1,2]})} ds_{\{n\}} = -2\pi r^{-1} r S_{\{n\}}(1) = -1 \text{ if } n=2,$$

$$\left. \frac{\partial \eta_{\{1\}}(x'_{[1,1]})}{\partial x_1} \right|_{x_{1,}=x_{*1}-r}^{x_{1,}=x_{*1}+r} = \left. \frac{\partial \eta_{\{1\}}(|x_{1,}-x_{*1,}|)}{\partial x_1} \right|_{x_{1,}=x_{*1}-r}^{x_{1,}=x_{*1}+r}$$

$$= -a_{\{1\}} \left. \frac{x_{1,}-x_{*1,}}{|x_{1,}-x_{*1,}|} \right|_{x_{1,}=x_{*1}-r}^{x_{1,}=x_{*1}+r} = -\frac{1}{2r}[r-(-r)] = -1. \tag{B.5}$$

QED.•

**Theorem B.2.** Given $n \in \omega_1$, given $\bar{x}_{*[1,n]} \in \bar{E}_{\{}$, for each $\bar{x}_{[1,n]} \in \bar{E}_n$:

$$\sum_{i=1}^{n} \frac{\partial^2 \eta^{\langle n \rangle}\left(\left|\bar{x}_{[1,n]} - \bar{x}_{*[1,n]}\right|\right)}{\partial x_i^2} = -\delta\left(\bar{x}_{[1,n]} - \bar{x}_{*[1,n]}\right). \tag{B.6}$$

**Proof:** Owing to Theorem B.1, integration of the expressions on both sides of equation (B.6) with respect to $\bar{x}_{[1,n]}$ over any closed domain $\bar{X}_n^c \subseteq \bar{E}_n$, which is *regular in each coordinate direction* and which contains the point $\bar{x}_{*[1,n]}$ in its interior $\bar{X}_n^o$, i.e. such that $\bar{x}_{*[1,n]} \in \bar{X}_n^o$, yields:

$$\int_{\bar{X}_n^c} \sum_{i=1}^{n} \frac{\partial^2 \eta^{\langle n \rangle}\left(\left|\bar{x}_{[1,n]} - \bar{x}_{*[1,n]}\right|\right)}{\partial x_i^2} dv_{\{n\}} = -\int_{\bar{X}_n^c} \delta\left(\bar{x}_{[1,n]} - \bar{x}_{*[1,n]}\right) dv_{\{n\}} = -1, \tag{B.7}$$

by the Ostrogradsky-Gauss theorem (see, e.g., Budak and Fomin [1978, pp. 218–224]). QED.•

## Appendix C: An *a posteriori* verification of Theorem 7.1

**Preliminary Remark C.1.** Given $\lambda \in \{-1,1\}$, given $n \in \omega_1$, in agreement with (3.29), formulas (7.24) and (7.25) can be rewritten as:

$$\Omega^{(0)}\left(\eta^{\langle n \rangle}, \psi^{\langle 1,n \rangle}, \tau_n, \bar{E}_n\right) = \int_{\bar{E}_n} K^{(0)}\left(\eta^{\langle n \rangle}, \psi^{\langle 1,n \rangle}, \tau_n, \bar{x}_{[1,n]}\right) dv_{\{n\}}, \tag{C.1}$$

$$\Omega^{(1)}\left(\eta^{\langle n \rangle}, \psi^{\langle 1,n \rangle}, \tau_n, \bar{E}_n\right) = \int_{\bar{E}_n} K^{(1)}\left(\eta^{\langle n \rangle}, \psi^{\langle 1,n \rangle}, \tau_n, \bar{x}_{[1,n]}\right) dv_{\{n\}}, \tag{C.2}$$

where, in agreement with (4.24), (4.25), and (4.74) (or (7.29)),

$$K^{(0)}\left(\eta^{\langle n \rangle}, \psi^{\langle 1,n \rangle}, \tau_n, \bar{x}_{[1,n]}\right) \equiv -\eta^{\langle n \rangle}\left(x'_{[1,n]}\right) f^{\langle 1,n \rangle}\left(x_{*0} + \lambda x'_{[1,n]}, \bar{x}_{[1,n]}\right), \tag{C.3}$$

$$K^{(0)}\left(\eta^{\langle n \rangle}, \psi^{\langle 1,n \rangle}, \tau_n, \bar{x}_{[1,n]}\right) \equiv \lambda(n-3)\zeta^{\langle n \rangle}\left(x'_{[1,n]}\right) \psi_{x_0}^{\langle 1,n \rangle}\left(x_{*0} + \lambda x'_{[1,n]}, \bar{x}_{[1,n]}\right). \tag{C.4}$$

So far we have studied integral properties of the definientia of (C.3) and (C.4), subject to (7.26)–(7.29) with respect to $x_0$ and $\bar{x}_{[1,n]}$. The following lemma is designed to study differential functional properties of the definientia with respect to $x_{*0}$ and $\bar{x}_{*[1,n]}$ with the purpose to establish, which ones of the integrals (C.1) and (C.2) with various $n$ can be differentiated under the sign of integral, especially by applying the pertinent variant of the d'Alambertian operator, defined as:

$$D_{*[0,n]} \equiv \left(\sum_{i=1}^{n} \frac{\partial^2}{\partial x_{*i}^2} - \frac{\partial^2}{\partial x_{*0}^2}\right). \tag{C.5}$$

**Lemma C.1.** 1) Given $\lambda \in \{-1,1\}$, given $n \in \omega_1$, for each $\bar{x}_{*0} \in R$, for each $\bar{x}_{*[1,n]} \in \bar{E}_n - \{\bar{x}_{[1,n]}\}$:

$$
\begin{aligned}
& D_{*[0,n]} K^{(0)}\left(\eta^{\langle n \rangle}, \psi^{\langle 1,n \rangle}, \tau_n, \bar{x}_{[1,n]}\right) \\
& = -D_{*[0,n]}\left[\eta^{\langle n \rangle}(x'_{[1,n]}) f^{\langle 1,n \rangle}(x_{*0} + \lambda x'_{[1,n]}, \bar{x}_{[1,n]})\right] \\
& = f^{\langle 1,n \rangle}(x_{*0}, \bar{x}_{*[1,n]}) \delta(\bar{x}_{[1,n]} - \bar{x}_{*[1,n]}) \\
& \quad - \lambda(n-3) \zeta^{\langle n \rangle}(x'_{[1,n]}) \frac{\partial f^{\langle 1,n \rangle}(x_{*0} + \lambda x'_{[1,n]}, \bar{x}_{[1,n]})}{\partial x_{*0}},
\end{aligned} \tag{C.6}
$$

$$
\begin{aligned}
& D_{*[0,n]} K^{(1)}\left(\eta^{\langle n \rangle}, \psi^{\langle 1,n \rangle}, \tau_n, \bar{x}_{[1,n]}\right) \\
& = \lambda(n-3) D_{*[0,n]}\left[\zeta^{\langle n \rangle}(x'_{[1,n]}) \psi^{\langle 1,n \rangle}_{x_0}(x_{*0} + \lambda x'_{[1,n]}, \bar{x}_{[1,n]})\right] \\
& = \lambda(n-3)(n-1) a_n \left[x'_{[1,n]}{}^{-(n+1)} \psi^{\langle 1,n \rangle}_{x_0}(x_{*0} + \lambda x'_{[1,n]}, \bar{x}_{[1,n]})\right. \\
& \quad \left. + 3\lambda x'_{[1,n]}{}^{-(n-1)} \psi^{\langle 1,n \rangle}_{x_0 x_0}(x_{*0} + \lambda x'_{[1,n]}, \bar{x}_{[1,n]})\right] \text{ if } n \neq 2,
\end{aligned} \tag{C.7a}
$$

$$
\begin{aligned}
& D_{*[0,2]} K^{(1)}\left(\eta^{\langle 2 \rangle}, \psi^{\langle 1,2 \rangle}, \tau_2, \bar{x}_{[1,2]}\right) \\
& = -\lambda D_{*[0,2]}\left[\zeta^{\langle 2 \rangle}(x'_{[1,2]}) \psi^{\langle 1,2 \rangle}_{x_0}(x_{*0} + \lambda x'_{[1,2]}, \bar{x}_{[1,2]})\right] \\
& = \lambda a_2 \left[\frac{\ln x'_{[1,2]} - 2}{x'_{[1,n]}{}^3} \psi^{\langle 1,2 \rangle}_{x_0}(x_{*0} + \lambda x'_{[1,2]}, \bar{x}_{[1,2]})\right. \\
& \quad \left. + \frac{a_2}{x'_{[1,2]}{}^2}(3\ln x'_{[1,2]} - 1) \psi^{\langle 1,2 \rangle}_{x_0 x_0}(x_{*0} + \lambda x'_{[1,2]}, \bar{x}_{[1,2]})\right] \text{ if } n = 2,
\end{aligned} \tag{C.7b}
$$

subject to (7.26)–(7.29) and (C.5). The character of singularity at the deleted point $\bar{x}_{[1,n]}$ of every term occurring in the final expression of each one of the three identities (C.6), (C.7a), and (C.7b) is self-evident. Since $\Delta_{*[1,n]} \eta^{\langle n \rangle}(x'_{[1,n]}) = 0$ for each $\bar{x}_{*[1,n]} \in \bar{E}_n - \{\bar{x}_{[1,n]}\}$, therefore we have utilized equation (B.6) in order to indicate the character of singularity of $\Delta_{*[1,n]} \eta^{\langle n \rangle}(x'_{[1,n]})$ at $\bar{x}_{*[1,n]} = \bar{x}_{[1,n]}$ with the understanding that $\delta(\bar{x}_{[1,n]} - \bar{x}_{*[1,n]}) = 0$ if $\bar{x}_{*[1,n]} \neq \bar{x}_{[1,n]}$.

**Proof:** 1) In analogy with (3.6), it follows from (2.43) and (2.44) that for each $i \in \omega_{1,n}$ and each $j \in \omega_{1,n}$:

$$
\frac{\partial x'_{[1,n]}}{\partial x_{*i}} = -\frac{x'_i}{x'_{[1,n]}}, \quad \frac{\partial^2 x'_{[1,n]}}{\partial x_{*i} \partial x_{*j}} = \frac{x'_{[1,n]}{}^2 \delta_{ij} - x'_i x'_j}{x'_{[1,n]}{}^3}, \tag{C.8}
$$

whence

$$
\sum_{i=1}^n \left(\frac{\partial x'_{[1,n]}}{\partial x_{*i}}\right)^2 = 1, \quad \sum_{i=1}^n \frac{\partial^2 x'_{[1,n]}}{\partial x_{*i}{}^2} = \sum_{i=1}^n \frac{x'_{[1,n]}{}^2 \delta_{ii} - x'_i x'_i}{x'_{[1,n]}{}^3} = \frac{n-1}{x'_{[1,n]}}. \tag{C.9}
$$

Also, either

$$f^{\langle 1,n\rangle}\left(x_{*0},\bar{x}_{*[1,n]}\right) \doteq -D_{*[0,n]}\psi^{\langle 1,n\rangle}\left(x_{*0},\bar{x}_{*[1,n]}\right), \tag{C.10}$$

by (4.21), or

$$D_{*[0,n]}\psi^{\langle 1,n\rangle}\left(x_{*0},\bar{x}_{*[1,n]}\right) = -f^{\langle 1,n\rangle}\left(x_{*0},\bar{x}_{*[1,n]}\right), \tag{C.11}$$

by (4.75), and

$$\sum_{i=1}^{n}\frac{\partial^2 \eta^{\langle n\rangle}\left(\left|\bar{x}_{[1,n]}-\bar{x}_{*[1,n]}\right|\right)}{\partial x_{*i}^{2}} = -\delta\left(\bar{x}_{[1,n]}-\bar{x}_{*[1,n]}\right), \tag{C.12}$$

by (B.6).

2) By (2.43) and (2.44), it follows that

$$\frac{\partial}{\partial x_{*i}}\left[\eta^{\langle n\rangle}\left(x'_{[1,n]}\right)f^{\langle 1,n\rangle}\left(x_{*0}+\lambda x'_{[1,n]},\bar{x}_{[1,n]}\right)\right]$$

$$= \frac{\partial \eta^{\langle n\rangle}\left(x'_{[1,n]}\right)}{\partial x_{*i}}f^{\langle 1,n\rangle}\left(x_{*0}+\lambda x'_{[1,n]},\bar{x}_{[1,n]}\right) \tag{C.13}$$

$$+ \lambda\eta^{\langle n\rangle}\left(x'_{[1,n]}\right)\frac{\partial x'_{[1,n]}}{\partial x_{*i}}\frac{\partial f^{\langle 1,n\rangle}\left(x_{*0}+\lambda x'_{[1,n]},\bar{x}_{[1,n]}\right)}{\partial x_{*0}},$$

whence

$$\frac{\partial^2}{\partial x_{*i}^{2}}\left[\eta^{\langle n\rangle}\left(x'_{[1,n]}\right)f^{\langle 1,n\rangle}\left(x_{*0}+\lambda x'_{[1,n]},\bar{x}_{[1,n]}\right)\right]$$

$$= \frac{\partial^2 \eta^{\langle n\rangle}\left(x'_{[1,n]}\right)}{\partial x_{*i}^{2}}f^{\langle 1,n\rangle}\left(x_{*0}+\lambda x'_{[1,n]},\bar{x}_{[1,n]}\right)$$

$$+ \lambda\left[\eta^{\langle n\rangle}\left(x'_{[1,n]}\right)\frac{\partial^2 x'_{[1,n]}}{\partial x_{*i}^{2}}+2\frac{\partial \eta^{\langle n\rangle}\left(x'_{[1,n]}\right)}{\partial x_{*i}}\frac{\partial x'_{[1,n]}}{\partial x_{*i}}\right]\frac{\partial f^{\langle 1,n\rangle}\left(x_{*0}+\lambda x'_{[1,n]},\bar{x}_{[1,n]}\right)}{\partial x_{*0}} \tag{C.14}$$

$$+ \lambda^2\eta^{\langle n\rangle}\left(x'_{[1,n]}\right)\left(\frac{\partial x'_{[1,n]}}{\partial x_{*i}}\right)^2\frac{\partial^2 f^{\langle 1,n\rangle}\left(x_{*0}+\lambda x'_{[1,n]},\bar{x}_{[1,n]}\right)}{\partial x_{*0}^{2}};$$

Therefore,

$$\left(\sum_{i=1}^{n}\frac{\partial^2}{\partial x_{*i}^{2}}-\frac{\partial^2}{\partial x_{*0}^{2}}\right)\left[\eta^{\langle n\rangle}\left(x'_{[1,n]}\right)f^{\langle 1,n\rangle}\left(x_{*0}+\lambda x'_{[1,n]},\bar{x}_{[1,n]}\right)\right]$$

$$= -\delta\left(\bar{x}_{[1,n]}-\bar{x}_{*[1,n]}\right)f^{\langle 1,n\rangle}\left(x_{*0}+\lambda x'_{[1,n]},\bar{x}_{[1,n]}\right) \tag{C.15}$$

$$+ \frac{\lambda}{x'_{[1,n]}}\left[(n-1)\eta^{\langle n\rangle}\left(x'_{[1,n]}\right)+2\sum_{i=1}^{n}x'_{i}\frac{\partial \eta^{\langle n\rangle}\left(x'_{[1,n]}\right)}{\partial x_{i}}\right]\frac{\partial f^{\langle 1,n\rangle}\left(x_{*0}+\lambda x'_{[1,n]},\bar{x}_{[1,n]}\right)}{\partial x_{*0}},$$

for $\lambda^2 = 1$. By (4.34a), (4.34b), and (4.74) (or (7.29a) and (7.29b)), it follows that for each $n \in \omega_1$:

$$\frac{1}{x'_{[1,n]}}\left[(n-1)\eta^{\langle n\rangle}(x'_{[1,n]})+2\sum_{i=1}^{n}x'_i\frac{\partial \eta^{\langle n\rangle}(x'_{[1,n]})}{\partial x_i}\right]=-(n-3)\zeta^{\langle n\rangle}(x'_{[1,n]}). \tag{C.16}$$

Along with (C.11) and (C.16), identity (C.15) is equivalent to (C.6).

3) The variants of identities (C.13) and (C.14) with '$\zeta^{\langle n\rangle}$' and '$\psi_{x_0}^{\langle 1,n\rangle}$' in place of '$\eta^{\langle n\rangle}$' and '$\psi^{\langle 1,n\rangle}$' respectively remain valid, whereas alike variant of (C.15) has the form:

$$\left(\sum_{i=1}^{n}\frac{\partial^2}{\partial x_{*i}^{2}}-\frac{\partial^2}{\partial x_{*0}^{2}}\right)\!\left[\zeta^{\langle n\rangle}(x'_{[1,n]})\psi_{x_0}^{\langle 1,n\rangle}(x_{*0}+\lambda x'_{[1,n]},\overline{x}_{[1,n]})\right]$$
$$=\sum_{i=1}^{n}\frac{\partial^2 \zeta^{\langle n\rangle}(x'_{[1,n]})}{\partial x_{*i}^{2}}\psi_{x_0}^{\langle 1,n\rangle}(x_{*0}+\lambda x'_{[1,n]},\overline{x}_{[1,n]}) \tag{C.17}$$
$$+\frac{\lambda}{x'_{[1,n]}}\left[(n-1)\zeta^{\langle n\rangle}(x'_{[1,n]})+2\sum_{i=1}^{n}x'_i\frac{\partial \zeta^{\langle n\rangle}(x'_{[1,n]})}{\partial x_i}\right]\psi_{x_0 x_0}^{\langle 1,n\rangle}(x_{*0}+\lambda x'_{[1,n]},\overline{x}_{[1,n]}).$$

a) For each $n\in\omega_1\text{-}\{2\}$, it follows from (4.74a) or (7.29a) that

$$a_n^{-1}\frac{\partial \zeta^{\langle n\rangle}(x'_{[1,n]})}{\partial x_{*i}}=\frac{\partial x'_{[1,n]}{}^{1-n}}{\partial x_{*i}}=-\frac{\partial x'_{[1,n]}{}^{1-n}}{\partial x_i}=\frac{(n-1)x'_i}{x'_{[1,n]}{}^{n+1}},$$
$$a_n^{-1}\frac{\partial^2 \zeta^{\langle n\rangle}(x'_{[1,n]})}{\partial x_{*i}\partial x_{*j}}=\frac{\partial}{\partial x_{*j}}\frac{(n-1)x'_i}{x'_{[1,n]}{}^{n+1}}=-\frac{\partial}{\partial x_j}\frac{(n-1)x'_i}{x'_{[1,n]}{}^{n+1}} \tag{C.18}$$
$$=-(n-1)\frac{x'_{[1,n]}{}^2 \delta_{ij}-(n+1)x'_i x'_j}{x'_{[1,n]}{}^{n+3}}$$

for each $i\in\omega_{1,n}$ and each $j\in\omega_{1,n}$.

Hence,

$$a_n^{-1}\sum_{i=1}^{n}x'_i\frac{\partial \zeta^{\langle n\rangle}(x'_{[1,n]})}{\partial x_i}=\frac{n-1}{x'_{[1,n]}{}^{n+1}}\sum_{i=1}^{n}x'_i x'_i=\frac{n-1}{x'_{[1,n]}{}^{n-1}},$$
$$a_n^{-1}\sum_{i=1}^{n}\frac{\partial^2 \zeta^{\langle n\rangle}(x'_{[1,n]})}{\partial x_{*i}^{2}}=-\frac{n-1}{x'_{[1,n]}{}^{n+3}}\sum_{i=1}^{n}\left[x'_{[1,n]}{}^2\delta_{ii}-(n+1)x'_i x'_i\right] \tag{C.19}$$
$$=-\frac{(n-1)[n-(n+1)]}{x'_{[1,n]}{}^{n+1}}=\frac{n-1}{x'_{[1,n]}{}^{n+1}},$$

so that

$$(n-1)\zeta^{\langle n\rangle}(x'_{[1,n]})+2\sum_{i=1}^{n}x'_i\frac{\partial \zeta^{\langle n\rangle}(x'_{[1,n]})}{\partial x_i}=3(n-1)a_n x'_{[1,n]}{}^{1-n}. \tag{C.20}$$

Therefore, (C.17) becomes

$$D_{*[0,n]}\left[\zeta^{\langle n\rangle}(x'_{[1,n]})\psi_{x_0}^{\langle 1,n\rangle}(x_{*0}+\lambda x'_{[1,n]},\overline{x}_{[1,n]})\right]$$
$$=\frac{(n-1)a_n}{x'_{[1,n]}{}^{n+1}}\left[\psi_{x_0}^{\langle 1,n\rangle}(x_{*0}+\lambda x'_{[1,n]},\overline{x}_{[1,n]})+3\lambda x'_{[1,n]}{}^2\psi_{x_0 x_0}^{\langle 1,n\rangle}(x_{*0}+\lambda x'_{[1,n]},\overline{x}_{[1,n]})\right] \quad \text{(C.21)}$$

which proves (C.7a).

b) For $n=2$, it follows from (4.74b) or (7.29b) that

$$a_2^{-1}\frac{\partial\zeta^{\langle 2\rangle}(x'_{[1,2]})}{\partial x_{*i}}=-\frac{\partial}{\partial x_{*i}}\frac{\ln x'_{[1,2]}+2}{x'_{[1,2]}}=\frac{\partial}{\partial x_i}\frac{\ln x'_{[1,2]}+2}{x'_{[1,2]}}$$
$$=\frac{[1-(\ln x'_{[1,2]}+2)]x'_i}{x'_{[1,2]}{}^3}=-\frac{(\ln x'_{[1,2]}-1)x'_i}{x'_{[1,2]}{}^3},$$

$$a_n^{-1}\frac{\partial^2\zeta^{\langle n\rangle}(x'_{[1,n]})}{\partial x_{*i}\partial x_{*j}}=-\frac{\partial}{\partial x_{*j}}\frac{(\ln x'_{[1,2]}-1)x'_i}{x'_{[1,2]}{}^3}=\frac{\partial}{\partial x_j}\frac{(\ln x'_{[1,2]}-1)x'_i}{x'_{[1,2]}{}^3} \quad \text{(C.22)}$$
$$=\frac{(\ln x'_{[1,2]}-1)\delta_{ij}}{x'_{[1,2]}{}^3}+\frac{[1-3(\ln x'_{[1,2]}-1)]x'_i x'_j}{x'_{[1,2]}{}^5}$$
$$=\frac{x'_{[1,2]}{}^2(\ln x'_{[1,2]}-1)\delta_{ij}-(3\ln x'_{[1,2]}-4)x'_i x'_j}{x'_{[1,2]}{}^5}$$

for each $i\in\omega_{1,2}$ and each $j\in\omega_{1,2}$.

Hence,

$$a_2^{-1}\sum_{i=1}^{2}x'_i\frac{\partial\zeta^{\langle 2\rangle}(x'_{[1,2]})}{\partial x_i}=-\frac{\ln x'_{[1,2]}-1}{x'_{[1,2]}{}^3}\sum_{i=1}^{n}x'_i x'_i=-\frac{\ln x'_{[1,2]}-1}{x'_{[1,2]}},$$

$$a_2^{-1}\sum_{i=1}^{2}\frac{\partial^2\zeta^{\langle 2\rangle}(x'_{[1,2]})}{\partial x_{*i}{}^2}=\frac{1}{x'_{[1,n]}{}^5}\sum_{i=1}^{2}\left[x'_{[1,2]}{}^2(\ln x'_{[1,2]}-1)\delta_{ii}-(3\ln x'_{[1,2]}-4)x'_i x'_i\right] \quad \text{(C.23)}$$
$$=\frac{2(\ln x'_{[1,2]}-1)-(3\ln x'_{[1,2]}-4)}{x'_{[1,n]}{}^3}=-\frac{\ln x'_{[1,2]}-2}{x'_{[1,n]}{}^3},$$

so that, particularly,

$$(n-1)\zeta^{\langle n\rangle}(x'_{[1,n]})+2\sum_{i=1}^{n}x'_i\frac{\partial\zeta^{\langle n\rangle}(x'_{[1,n]})}{\partial x_i}$$
$$=-\frac{a_2}{x'_{[1,2]}}\left[\ln x'_{[1,2]}+2+2(\ln x'_{[1,2]}-1)\right]=-\frac{a_2}{x'_{[1,2]}}(3\ln x'_{[1,2]}-1). \quad \text{(C.24)}$$

Therefore, (C.17) becomes

$$D_{*[0,2]}\left[\zeta^{\langle 2\rangle}(x'_{[1,2]})\psi_{x_0}^{\langle 1,2\rangle}(x_{*0}+\lambda x'_{[1,2]},\overline{x}_{[1,2]})\right]$$
$$=-\frac{a_2(\ln x'_{[1,2]}-2)}{x'_{[1,n]}{}^3}\psi_{x_0}^{\langle 1,n\rangle}(x_{*0}+\lambda x'_{[1,n]},\overline{x}_{[1,n]}) \quad \text{(C.25)}$$
$$-\frac{a_2\lambda}{x'_{[1,n]}{}^2}(3\ln x'_{[1,2]}-1)\psi_{x_0 x_0}^{\langle 1,n\rangle}(x_{*0}+\lambda x'_{[1,n]},\overline{x}_{[1,n]}),$$

which proves (C.7b). QED.•

**Theorem C.1.** 1) By (C.6) it follows from (C.1) that

$$D_{*[0,n]}\Omega^{(0)}\left(\eta^{\langle n\rangle},\psi^{\langle 1,n\rangle},\tau_n,\overline{E}_n\right) = \int_{\overline{E}_n} D_{*[0,n]} K^{(0)}\left(\eta^{\langle n\rangle},\psi^{\langle 1,n\rangle},\tau_n,\overline{x}_{[1,n]}\right) dv_{\{n\}}$$

$$= f^{\langle 1,n\rangle}\left(x_{*0},\overline{x}_{*[1,n]}\right) - \lambda(n-3)\int_{\overline{E}_n} \zeta^{\langle n\rangle}\left(x'_{[1,n]}\right) \frac{\partial \mathcal{F}^{\langle 1,n\rangle}\left(x_{*0}+\lambda x'_{[1,n]},\overline{x}_{[1,n]}\right)}{\partial x_{*0}} dv_{\{n\}}. \quad (C.26)$$

In the case of $n \neq 3$, since $\partial \mathcal{F}^{\langle 1,n\rangle}\left(x_{*0}+\lambda x'_{[1,n]},\overline{x}_{[1,n]}\right)/\partial x_{*0}$ is supposedly bounded in a neighborhood of the point $\overline{x}_{*[1,n]}$, therefore the integral in the final expression of (C.26) converges as $x'_{[1,n]}$ in the integrand is contracted to 0.

2) Owing to (7.25) and (C.26), application of the operator $D_{*[0,3]}$ to both sides of identity (7.23) at $n=3$, i.e. of identity (7.30), results in the identity:

$$D_{*[0,3]}\psi^{\langle 1,3\rangle}\left(x_{*0},\overline{x}_{*[1,3]}\right) = -f^{\langle 1,3\rangle}\left(x_{*0},\overline{x}_{*[1,3]}\right) \quad (C.27)$$

ss expected.

3) By (C.7a), it follows from (C.2) at $n=1$ that

$$D_{*[0,1]}\Omega^{(0)}\left(\eta^{\langle 1\rangle},\psi^{\langle 1,1\rangle},\tau_1,\overline{E}_1\right) = \int_{\overline{E}_1} D_{*[0,1]} K^{(1)}\left(\eta^{\langle 1\rangle},\psi^{\langle 1,1\rangle},\tau_1,x_1\right) dx_1 = 0. \quad (C.28)$$

In all other cases, (C.7a) and (C.7b) evidence that $D_{*[0,n]} K^{(1)}\left(\eta^{\langle n\rangle},\psi^{\langle 1,n\rangle},\tau_n,\overline{x}_{[1,n]}\right)$ diverges as $x'_{[1,n]}$ is contracted to 0 so that the integral

$$\int_{\overline{E}_n} D_{*[0,n]} K^{(1)}\left(\eta^{\langle n\rangle},\psi^{\langle 1,n\rangle},\tau_n,\overline{x}_{[1,n]}\right) dv_{\{n\}}$$

diverges in a neighborhood of the point $\overline{x}_{[1,n]} = \overline{x}_{*[1,n]}$. Therefore, the operator $D_{*[0,n]}$ cannot be applied to the integral $\int_{\overline{E}_n} K^{(1)}\left(\eta^{\langle n\rangle},\psi^{\langle 1,n\rangle},\tau_n,\overline{x}_{[1,n]}\right) dv_{\{n\}}$ under the sign of integral.

## Acknowledgements

I am indebted to my son Gil Iosilevskii, Associate Professor with the Faculty of Aerospace Engineering at Technion, for reading the draft of this exposition, for verifying my key calculations, for useful criticisms and suggestions, and for his technical assistance in preparing its final version. I am also indebted to Lev P. Pitaevsky, an academician (full member) of the Russian Academy of Sciences and Professor of Theoretical Physics with the Department of Physics at the University of Trento (Italy) for his useful criticism, which stimulated my more adequate interpretation and discussion of the final result, expressed by equation (1.15).